\documentclass[12pt]{amsart}
\usepackage{amsmath, amsthm, amssymb,cite}
\usepackage{fullpage}
\usepackage[usenames,dvipsnames,svgnames,table]{xcolor}
\usepackage[usenames,dvipsnames]{xcolor}
\usepackage{hyperref}
\usepackage{enumitem}
\usepackage[english]{babel}
\usepackage[normalem]{ulem}
\renewcommand{\P}{\mathbf P}

\newtheorem{theorem}{Theorem}
\newtheorem{lemma}{Lemma}
\newtheorem{proposition}[lemma]{Proposition}

\numberwithin{lemma}{section}

\newcommand{\Ez}{{\mathcal E_0}}
\newcommand{\E}{{\mathcal E}}
\newcommand{\Elind}{E^{(2)}_{lin}}
\newcommand{\Elint}{E^{(3)}_{lin}}

\newcommand{\End}{E^{n,(2)}}
\newcommand{\Ent}{E^{n,(3)}}
\newcommand{\Enthigh}{E^{n,(3)}_{high}}
\newcommand{\Enthighc}{E^{n,(3)}_{high,c}}

\newcommand{\Ennf}{E^{n}_{NF}}
\newcommand{\Ennfhigh}{E^{n}_{NF,high}}
\newcommand{\Ennfhighc}{E^{n,c}_{NF,high}}
\newcommand{\Ennflow}{E^{n}_{NF,low}}
\newcommand{\M}{\mathfrak M}

\numberwithin{equation}{section}

\newcommand{\R}{{\mathbb R}}

\newcommand{\tW}{{\tilde W}}
\newcommand{\tQ}{{\tilde Q}}

\renewcommand{\H}{{\mathcal H }}
\renewcommand{\AA}{\mathbf A}

\newcommand{\bR}{\mathbf R}

\newcommand{\NN}{{n}}

\usepackage[makeroom]{cancel}

\newcommand{\tR}{{\tilde {\mathbf R}}}
\newcommand{\dH}{{\dot{\mathcal H} }}

\newcommand{\W}{{\mathbf W}}

\newcommand{\err}{\text{\bf err}}
\newcommand{\errw}{\text{\bf err}(L^2)}
\newcommand{\errr}{\text{\bf err}(\dot H^\frac12)}
\newcommand{\norm}{{\mathbf N}}
\newcommand{\ub}{\underline{b}}
\newcommand{\ua}{\underline{a}}
\newcommand{\uM}{\underline{M}}
\newcommand{\uA}{\underline{A}}
\newcommand{\uB}{\underline{B}}
\newcommand{\mfw}{{\mathfrak{w}}}
\newcommand{\mfr}{{\mathfrak{r}}}
\newcommand{\uG}{\underline{G}}
\newcommand{\uK}{\underline{K}}

\begin{document}

\title{Two dimensional gravity water waves with constant vorticity: I. Cubic lifespan}

\author{Mihaela Ifrim}
\address{Department of Mathematics, University of California at Berkeley}
\thanks{The first author was supported by the Simons Foundation}
\email{ifrim@math.berkeley.edu}

\author{ Daniel Tataru}
\address{Department of Mathematics, University of California at Berkeley}
 \thanks{The second author was partially supported by the NSF grant DMS-1266182
as well as by the Simons Foundation}
\email{tataru@math.berkeley.edu}

\begin{abstract}
  This article is concerned with the incompressible, infinite depth
  water wave equation in two space dimensions, with gravity and
  constant vorticity but with no surface tension.  We consider this
  problem expressed in position-velocity potential holomorphic
  coordinates, and prove local well-posedness for large data,
 as well as cubic  lifespan bounds for   small data solutions.  
\end{abstract}

\maketitle

\section{Introduction}
The motion of water in contact with air is well described by the incompressible Euler 
equations in the fluid domain, combined with two boundary conditions on the 
free surface, i.e., the interface with air. In special, but still physically relevant cases,
the equations of motion can be viewed as  evolution equations for the free 
surface. These equations are  commonly referred to as the water wave equations.
Most notably, this is the case for irrotational flows. However, in two space 
dimensions there is a natural extension of these equations to flows  
with constant vorticity.
 
In previous work \cite{HIT}, \cite{IT2} we have considered the local
and long time behaviour for the irrotational gravity wave equations
with infinite depth in two space dimensions. In this article we take a 
first step toward understanding the more difficult question of the
long time behavior of gravity waves with infinite
depth and constant vorticity, either in $\R \times \R$ or in the
periodic case $\R \times \mathbb T$.  We begin by establishing a local
well-posedness result. Then we consider the lifespan of small data
solutions, where, like in the zero vorticity case \cite{HIT}, we are
able to prove cubic lifespan bounds. To the best of our knowledge,
this is the first long time well-posedness result for this problem.

We remark that it is of further interest to consider small localized data,
and establish an almost global in time result, as it was done in the irrotational case
in \cite{HIT}, improving an earlier result of Wu~\cite{wu}. However, in the constant
vorticity case this problem presents some interesting new challenges. We hope to consider 
this in subsequent work.

The motivation to study the constant vorticity problem comes from multiple sources.
On one hand, from a mathematical perspective, it provides us with the possibility 
to consider vorticity  effects in a framework  where the equations of motion 
can still be described in terms of the water/air interface, while allowing for a larger
range of dynamic behavior, which is particularly interesting over large time scales.
On the other hand, from a practical perspective, constant vorticity flows are good models 
for the water motion in the presence of countercurrents. An interesting example of this 
type is provided by tidal effects.

The constant vorticity problem is a subset of the full vorticity problem, and as such, 
local well-posedness for regular enough data can be viewed as a consequence of
results for the general problem. For this we refer the reader to \cite{DL,HL,cs,sz}.
There has also been a considerable body of work devoted to 
the study of solitary waves in constant vorticity flows. A good source of information
in this direction is provided by several recent articles \cite{MR2754336}, \cite{MR3217710},
\cite{MR3351023}, as well as the 
survey article \cite{MR2721042}. Various ways of formulating the equations
have been described in the articles \cite{MR2309783}, \cite{MR2370681}, 
\cite{MR2864943}. 

The conformal formulation for two dimensional water waves, which we adopt here,
following our previous work \cite{HIT},  originates in early work on traveling waves,
see e.g.  Levi-Civita~\cite{LC}.  For the dynamical problem, to the best of our knowledge it 
first appears  in work of Ovsjannikov~\cite{ov}, but was better developed later in work of Wu~\cite{wu2}
and also Dyachenko-Kuznetsov-Spector-Zakharov~\cite{zakharov}.
It  has been widely used since then in order to study a variety of water wave problems. However, to the best of our 
our knowledge, this is the first article where this formulation is fully implemented
in the constant vorticity case.

\subsection{Equations of motion}
The water flow is governed by the  incompressible Euler equations in the fluid domain
$\Omega_t$,  with a dynamic and a kinematic boundary condition on the free surface
of the fluid $\Gamma_t$.   Denoting the fluid velocity by $\mathbf{u}(t,x,y)=(u(t,x,y),v(t,x,y))$, and the pressure by $p(t,x,y)$, 
the equations of motion are the equation of mass conservation
\begin{equation}
u_x+v_y=0 \mbox{ in } \Omega(t),
\end{equation}
and Euler's equations also in $\Omega(t)$ 
\begin{equation}
\label{euler}
\left\{
\begin{aligned}
& u_t + u  u_x +vu_y = - p_x 
\\
&v_t+ uv_x+vv_y=-p_y-{g},
\end{aligned}
\right.
\end{equation}
where  ${g}$ is the gravitational constant.
The boundary conditions for capillary-gravity waves are the dynamic boundary condition
\begin{equation}
p=p_0 \, \mbox{ on } \Gamma(t); 
\end{equation}
$p_0$ being the constant atmospheric pressure, and the kinematic
boundary condition, which asserts that the free boundary $\Gamma(t)$
is transported along the flow.  Since we are in the two dimensional
case, the vorticity will also be transported along the flow. This
makes it possible to study flows with a non-zero constant vorticity
field,
\[
\omega = u_y-v_x=-c, \mbox{ where } c \mbox{ is a constant}.
\]
Then the velocity field $\bf u$ can be represented as
\[
{\bf u}=(cy+\varphi_x, \varphi_y),
\]
where $\varphi (t, x,y)$ is called the (generalized) velocity
potential.  Here $\varphi$ it is defined up to an arbitrary function
of time and, by the incompressibility condition, it satisfies the
Laplace equation in $\Omega (t)$
\[
\bigtriangleup \varphi (t,x,y)=0.
\]
This brings us to our boundary condition on the bottom, which asserts 
that
\[ 
\lim_{y \to -\infty} \varphi(x,y) = 0, \qquad \text{uniformly in $x$}.
\]
Then $\varphi$ is uniquely determined by its values on the free surface $\Gamma(t)$.

Introducing its harmonic conjugate $\theta (t,x,y)$:
\[
\varphi_x=\theta_y, \quad \varphi_y=-\theta_x, 
\]
we can rewrite the equations in \eqref{euler} as a single scalar
equation in the fluid domain:
\begin{equation*}
\nabla\left( \varphi_t-c\theta+cy\varphi_x+\frac{1}{2}(\varphi_x^2+\varphi_y^2)+gy\right) = 0. 
\end{equation*}
 Since $\varphi$ is only defined up to an
arbitrary function of time, we can absorb a function of time into
$\varphi$. Using the fact that the pressure is constant on $\Gamma (t)$, 
we obtain the following analogue of the Bernoulli's equation:
\begin{equation}
\label{imp}
\varphi_t-c\theta+cy\varphi_x+\frac{1}{2}(\varphi_x^2+\varphi_y^2)+gy =0  \mbox{ on } \Gamma(t).
\end{equation}
This is the equation that makes possible reduction of dimensionality
of the problem, in the same manner as for purely potential flows. This
is in combination with the kinematic boundary condition, which is already 
restricted to $\Gamma (t)$. We remark
that expressing $\theta$ and the full gradient $\nabla \varphi$ on $\Gamma (t)$ 
in terms of the restriction of $\varphi$ to $\Gamma (t)$ requires using the 
Dirichlet to Neumann map associated to the fluid domain, and in turn makes 
our equations nonlocal.

\subsection{The equations in holomorphic coordinates}
The first issue we need to address is the choice of coordinates. Here we take our cue from \cite{HIT}, and work
in holomorphic coordinates. There are also other ways of expressing the equations, for instance in Cartesian coordinates using the Dirichlet to Neumann map associated to the water domain, see e.g. \cite{abz}.  However, we prefer the holomorphic coordinates due to the simpler form of the equations; in particular, in these coordinates the Dirichlet to Neumann map is diagonalized, and given in terms of the standard Hilbert transform.

We obtain a system which models the time evolution of the free
surface, described via a function $W$ as the graph of the function
$\alpha \to W(\alpha)+\alpha$, and that of the holomorphic velocity
potential $Q = \varphi+i\theta$ restricted to the free surface.  The derivation of the
equations is relegated to \emph{Appendix}~\ref{s:eq}, also using some
of the analysis from \emph{Appendix} in \cite{HIT}. Some minor changes
are needed for the periodic case; these are also described in
\cite{HIT}, \emph{Appendix B}. The outcome of this computation
is the following set of equations:
\begin{equation}
\label{ww2d1}
\left\{
\begin{aligned}
&W_t+ (W_{\alpha}+1)\underline{F}+i\frac{c}{2}W=0\\
&Q_t-igW+\underline{F}Q_{\alpha}+icQ+\mathbf{P}\left[\frac{\vert Q _{\alpha}\vert ^2}{J}\right]-i\frac{c}{2} T_1=0,
\end{aligned}
\right.
\end{equation}
where $ J := |1+W_\alpha|^2$, $\P$ is the projection onto negative wave numbers
\[
\P= \frac12(I-iH), \mbox{ with } H \mbox{ the Hilbert transform, }
\]
and 
\begin{equation*}
\begin{aligned}
 F :=& \P\left[\frac{Q_\alpha - \bar Q_\alpha}{J}\right],\quad F_1:=
\P\left[ \frac{W}{1+\bar{W}_{\alpha}}+ \frac{\bar{W} }{1+W_{\alpha}}\right],\\
\underline{F}:=& F-i\frac{c}{2}F_1, \qquad \quad  \ T_1:= \P\left[  \frac{W\bar{Q}_{\alpha}}{1+\bar{W}_{\alpha}} - \frac{\bar{W}Q_{\alpha}}{1+W_{\alpha}}\right]  .
\end{aligned}
\end{equation*}
These equations are considered either in $\mathbb{R} \times
\mathbb{R}$ or in $\mathbb{R} \times \mathbb{S}^1$. They model an
evolution in the space of functions which admit bounded holomorphic
extensions to the lower half-space; equivalently, their Fourier transform is 
supported on the negative real line. By a slight abuse, we call such functions 
holomorphic.

This is a Hamiltonian system, where the Hamiltonian has the form
\begin{equation}\label{ww-energy-re}
\begin{aligned}
\E(W,Q) =& \Re \int  g\vert W\vert^2(1+W_\alpha) - iQ\bar{Q}_{\alpha}
+ c  Q_{\alpha} (\Im W)^2 -
\frac{c^3}{2i} |W|^2W(1+W_\alpha) \, d\alpha . \\
\end{aligned}
\end{equation}
A second conserved quantity is the horizontal momentum,
\begin{equation}
\mathcal{P}(W,Q)=\int \left\lbrace\frac{1}{i}\left( \bar{Q}W_{\alpha}-Q\bar{W}_{\alpha}\right)-c\vert W\vert^2+\frac{c}{2}\left( W^2\bar{W}_{\alpha}+\bar{W}^2W_{\alpha}\right)\right\rbrace\, d\alpha ,
\label{horizontal-m-re}
\end{equation}
which is the Hamiltonian for the group of horizontal translations.

We remark that if $c = 0$ then the equations \eqref{ww2d1} are the
gravity water wave equations studied in \cite{HIT}. The sign of $c$ is
not important, as it can be switched via the flip $\alpha \to
-\alpha$. For convenience we assume $c > 0$.  The space-time scaling
\[
(W(t,\alpha),Q(t,\alpha)) \to \lambda^{-2} W(\lambda t, \lambda^2 x), 
\lambda^{-3} Q(\lambda t, \lambda^2 x)
\]
is a symmetry for gravity waves (thus it leaves $g$ unchanged). Here
it changes $c \to \lambda c$.  Finally, the purely spatial scaling
\[
(W(t,\alpha),Q(t,\alpha)) \to \lambda^{-2} W( t, \lambda^2 x), 
\lambda^{-3} Q( t, \lambda^2 x)
\]
has the effect of leaving $c$ unaffected, but it changes $g \to
\lambda^{-1}g$.  Thus one could use scaling considerations to set both
$c = 1$ and $g=1$.  However, we choose not to do that, and instead use
the coefficients $c$ and $g$ to keep better track of the expressions
arising in our analysis.  In this context, it is useful to observe
that one can interpret $c^2/g$ as an (inverse) semiclassical
parameter, so that all energy expressions can be viewed as
homogeneous.

We further remark that the terms involving $c$ are lower order, though they cannot be
viewed as bounded. Thus, it is natural to expect that the local theory for this problem 
is not very different from the gravity waves; indeed, our results in this regard
are similar to \cite{HIT}.  However, we will see that the long time behavior 
is quite different in the constant vorticity case.

To further motivate our expectations concerning this system, 
we note that the linearization of the system \eqref{ww2d1} around the zero solution is 
\begin{equation} \label{ww2d-0} \left\{
\begin{aligned}
 & w_{ t} +  q_\alpha   =  0
\\
& q_t +icq- ig w = 0,
\end{aligned}
\right.
\end{equation}
restricted to holomorphic functions (in our terminology, these are functions with negative spectrum).
It is not difficult to see that this is a dispersive equation. Expressed as a second order equation 
this becomes
\begin{equation}
w_{tt} + icw_t+igw_{\alpha}=0.
\end{equation}
From here we can deduce the associated dispersion relation,
\begin{equation} \label{dispersion}
\tau^2 + c \tau  + g \xi=0, \qquad \xi \leq 0.
\end{equation}
This is the intersection of a lateral parabola with the left half space. It has two branches,
intersecting the axis $\xi=0$ at $\tau = 0$ and $\tau = -c$.

 The conserved energy, respectively 
momentum for \eqref{ww2d-0} are 
\[
\mathcal E_0(w,r) = \int |w|^2 - i q \bar q_\alpha \ d\alpha = \|w\|_{L^2}^2 + \| q\|_{\dot H^\frac12}^2,
\]
\[
\mathcal P_0(w,r) = \int \frac{1}{i}\left( \bar{q}w_{\alpha}-q\bar{w}_{\alpha}\right)-c\vert w\vert^2\, d\alpha .
\]
The former motivates the functional setting where we will study the
equations \eqref{ww2d1}.  The system \eqref{ww2d-0} is a well-posed
linear evolution in the space $\dH_0$ of holomorphic functions endowed
with the $L^2 \times \dot H^{\frac12}$ norm. To measure higher regularity
we will use the spaces $\dH_n$ endowed with the norm
\[
\| (w,r) \|_{\dH_n}^2 := \sum_{k=0}^n
\| \partial^k_\alpha (w,r)\|_{ L^2 \times \dot H^\frac12}^2,
\]
where $n \geq 1$.

\subsection{ The differentiated equations and digonalization}
As the system \eqref{ww2d1} is fully nonlinear, a standard procedure
is to convert it into a quasilinear system by differentiating it. In
the case of gravity waves this yields a self-contained first order
quasilinear system for $(W_\alpha,Q_\alpha)$. This is no longer true
here, precisely due to the contributions from $F_1$ and $T_1$.

To write the differentiated system, as well as the linearized
system later on, we introduce several notations. First, as in \cite{HIT},  we set
\[
\W = W_\alpha, \qquad R = \frac{Q_\alpha}{1+\W}, \qquad Y = \frac{\W}{1+\W}.
\]
The expression $R$ has an intrinsic meaning, namely it is the complex
velocity on the water surface. $Y$, on the other hand, is introduced
for computational reasons only, in order to avoid rational expressions
above, and in many places in the sequel.

We also need two key auxiliary real functions.  The first is
$\underline{b}$, which we call the {\em advection velocity}, and is
given by
\begin{equation}
\begin{split}
\underline{b}:=&\ b-i\frac{c}{2}b_1,
\\
 b := \P \left[\frac{{Q}_\alpha}{J}\right] +  \bar \P\left[\frac{\bar{Q}_\alpha}{J}\right],&\qquad b_1:= \P\left[ \frac{W}{1+\bar{\W}}\right]-\bar{\P}\left[ \frac{\bar{W}}{1+\W}\right].
\label{defb}
\end{split}
\end{equation}
The second is the {\em frequency-shift} $a$, given by
\begin{equation}
\begin{split}
\underline{a}:=\ a+\frac{c}{2}a_1,\quad a := i\left(\bar \P \left[\bar{R} R_\alpha\right]- \P\left[R\bar{R}_\alpha\right]\right),  & \quad
a_1= R+\bar{R}- N,
\label{defa}
\end{split}
\end{equation}
where
\begin{equation}
N:=\P\left[W\bar{R}_{\alpha}-\bar \W R\left]+\bar{\P}\right[\bar{W}R_{\alpha}- \W \bar R\right].
\label{defN}
\end{equation}
The functions $\underline{a}$ and $\underline b$ are the leading order
coefficients in the linearized equation, and thus fully describe the
quasilinear nature of the problem.  The linearized equation is more
involved and is described in full later, but, as a baseline, the reader
should keep in mind the linear system
\begin{equation} \label{ww2d-0-model} 
\left\{
\begin{aligned}
 & w_{ t} + b w_\alpha +   q_\alpha   =  0
\\
& q_t +icq- i(g+a) w = 0.
\end{aligned}
\right.
\end{equation}

The real function $g+\ua$ has a physical interpretation as the normal
derivative of the pressure on the interface; this is proved in the
\emph{Appendix}. For more regular irrotational waves (i.e., $c=0$) this was proved
by Wu~\cite{wu2} to be positive; an alternate proof was provided in
\cite{HIT} in the context of holomorphic coordinates, assuming only 
scale invariant regularity $(\W,R) \in \dH_{\frac12}$. This positivity,
called the \emph{Taylor sign condition}, was crucial for the well-posedness
of the water wave system, both with gravity \cite{wu2}, \cite{HIT}
and with surface tension \cite{IT3}.

For the present problem we still need to know that $g+\ua$  
is nonnegative, which corresponds with the normal
derivative of the pressure being bounded away from zero. 
But  it is no longer the case that this comes for free. Thus, we will
impose  the positivity condition on the normal derivative of
the pressure  and solve the equations for as long as this condition
holds uniformly. We remark that this is always the case if we assume
that $R$ is small in $L^{\infty}$, which is the case for our small
data result.

Differentiating with respect to $\alpha$ yields a system for
$(W_\alpha,Q_\alpha)$, which turns out to be degenerate hyperbolic
with a double speed $\underline b$. This is explained in detail in the
context of irrotational gravity waves in \cite{HIT}, and easily carries
over here as the highest order terms in the equations are the same.
Then it is natural to diagonalize it. This is also done exactly as in
\cite{HIT}, using the operator
\begin{equation}
\AA(w,q) := (w,q - Rw), \qquad R := \frac{Q_\alpha}{1+W_\alpha}.
\label{defR}
\end{equation}
Noting that
\[
\AA(W_\alpha,Q_\alpha) = (\W,R), \qquad \W : = W_\alpha,
\]
it follows that the pair $(\W,R)$ diagonalizes the differentiated system. 
Thus, rather than repeating the more extensive computations in
\cite{HIT}, here we directly take advantage of this knowledge to
arrive more efficiently at the differentiated system for the diagonal variables 
$(\W,R)$.

We first introduce  $\underline{b}$ into the equations using the relations 
\begin{equation}\label{ftob}
F = b - \frac{\bar R}{1+\W}, \qquad F_1 = b_1 + \frac{\bar W}{1+\W},
\end{equation}
so that the system \eqref{ww2d1} is written in the form
\begin{equation*}
\left\{
\begin{aligned}
&W_t + \underline{b} (1+ W_\alpha) +i\frac{c}{2}W=  \bar R +i\frac{c}{2}\bar{W}
\\
&Q_t +  \underline{b} Q_\alpha - igW+icQ -i\frac{c}{2}\bar R W=  \bar \P\left[ |R|^2 \right]- i\frac{c}{2}\bar{\P}\left[W\bar R - \bar{W} R\right].
\end{aligned}
\right.
\end{equation*}
Here the terms on the right are antiholomorphic and disappear when the
equations are projected onto the holomorphic space.

Next we differentiate the equations. For the first equation we need
the expression for $\underline{b}_{\,\alpha}$, for which we introduce one last set of
quadratic auxiliary functions $M$, $M_1$ and $\underline{M}$ as
follows:
\begin{equation}\label{M-def}
\begin{aligned}
&\qquad \qquad \qquad \qquad  \underline{M}:=M-i\frac{c}{2}M_1, \\
&M :=  \frac{R_\alpha}{1+\bar \W}  + \frac{\bar R_\alpha}{1+ \W} -  b_\alpha =
\bar \P [\bar R Y_\alpha- R_\alpha \bar Y]  + \P[R \bar Y_\alpha - \bar R_\alpha Y],\\
&\\
& \qquad \qquad  M_1:=\W-\bar{\W}- b_{1,\alpha}=\P\left[ W\bar{Y}\right]_{\alpha}-\bar{\P}\left[ \bar{W}Y\right]_{\alpha}.
\end{aligned}
\end{equation}
Thus, we can substitute $\underline{b}_{\, \alpha}$ with 
\begin{equation}\label{M-def1}
\underline{b}_{\,\alpha} = \frac{R_\alpha}{1+\bar \W}  + \frac{\bar R_\alpha}{1+ \W} - i \frac{c}2  (\W-\bar{\W}) - \uM .
\end{equation}
For the second equation we switch directly to $R$, and then $\underline{b}_{\, \alpha}$ is no longer 
needed in view of the identity
\[
(\partial_t + \ub \partial_\alpha) R = \partial_\alpha (\partial_t + \ub \partial_\alpha) Q
- R \partial_\alpha ( W_t + \ub( 1+W_\alpha)).
\]
Taking the above discussion into account, after some straightforward computations, one arrives 
at the system
\begin{equation} \label{ww2d-diff}
\left\{
\begin{aligned}
 & \W_{ t} + \underline{b} \W_{ \alpha} + \frac{(1+\W) R_\alpha}{1+\bar \W}  =  (1+\W)\underline{M}+i\frac{c}{2}\W(\W-\bar{\W})
\\
&R_t  + \underline{b} R_{\alpha} +icR- i\frac{g\W-a}{1+\W} =i\frac{c}{2}\frac{R\W + \bar R \W +N}{1+\W}.
\end{aligned}
\right.
\end{equation}
This governs an evolution in the space of
holomorphic functions, and will be used both directly and in its
projected version. Obviously the two forms are algebraically equivalent.

We remark that while the transport coefficient $\ub$ appears
explicitly in these equations, the frequency shift $\ua$ does
not. This is due to the fact that the right hand side of the second
equation is still fully nonlinear in $\W$. To understand better the role
played by $\ua$ one needs to consider either the linearized equations,
which are discussed in Section~\ref{s:linearized}, or to further differentiate 
\eqref{ww2d-diff}, as in Section~\ref{s:ee}.

\subsection{The main results}
 To describe the lifespan of the solutions we begin with the control norms
in \cite{HIT}, namely
\begin{equation}\label{A-def}
A := \|\W\|_{L^\infty}+\| Y\|_{L^\infty} + \||D|^\frac12 R\|_{L^\infty \cap B^{0,\infty}_{2}},
\end{equation}
respectively
\begin{equation}\label{B-def}
B :=\||D|^\frac12 \W\|_{BMO} + \| R_\alpha\|_{BMO},
\end{equation}
where $|D|$ represents the multiplier with symbol $|\xi|$. In order to 
estimate lower order terms introduced in conjunction with $c$ we also 
introduce 
\begin{equation}\label{A1-def}
A_{-1/2} := \||D|^\frac12W\|_{L^\infty} + \| R\|_{L^\infty },
\end{equation}
and 
\begin{equation}\label{A2-def}
A_{-1} := \|W\|_{L^\infty}.
\end{equation}
It is also useful to introduce the notations
\begin{equation}
\uB := B + cA + c^2 A_{-1/2}, \qquad \uA := A + c A_{-1/2}+ c^2 A_{-1}.
\end{equation}
Here $A$ is a scale invariant quantity, while $B$ corresponds to the
homogeneous $\dH_1$ norm of $(\W, R)$ and $A_{-\frac12}$ corresponds
to the homogeneous $\dH_0$ norm of $(\W, R)$.  We note that $B$, and
all but the $Y$ component of $A$ are controlled by the energy and $\dH_1$ norm of
$(\W,R)$. 

Now we are ready to state our main results. We begin with the 
local well-posedness result:
\begin{theorem}
\label{baiatul}
 Let  $ n \geq 1$.  The system \eqref{ww2d1} is locally well-posed
for initial data $(W_0,Q_0)$ with the following regularity:
\[
(W_0,Q_0)  \in \dH_0, \qquad (\W_0,R_0) \in \dH_1,
\]
and satisfying the pointwise constraints
\begin{equation}\label{nocusp}
|\W(\alpha)+1| > \delta > 0, \qquad \text{(no interface singularities)}
\end{equation}
\begin{equation} \label{taylor}
  g+ \ua(\alpha) > \delta > 0 \qquad \text{(Taylor sign condition)}.
\end{equation}
Further, the solution can be continued for as long as $\uA$ and $\uB$ 
remain bounded and the pointwise conditions above hold uniformly.
The same result holds in the periodic setting.
\end{theorem}

The well-posedness above should be interpreted in the 
sense of Hadamard. To be more precise, it means that there
exists some time $T > 0$, depending only on the initial data size 
and on the constant $\delta$ in the above pointwise constraints,
so that the following properties hold:  
\begin{enumerate}
\item{(Regular data)} For each data $(W_0,Q_0)$, which is as above, but with additional regularity 
$(\W_0,R_0) \in \dH_n$, with $n \geq 2$, there exists a unique solution $(W,Q)$ in $[0,T]$,
with the property that 
\[
\| (\W,R)\|_{C[0,T;\dH_k]} \lesssim \| (\W_0,R_0)\|_{\dH_k}, \qquad 0 \leq k \leq n .
\]

\item{(Rough data)} For each data  $(W_0,Q_0)$ as above there exists a solution 
$(W,Q)$ in $[0,T]$, with the property that  
 \[
\| (\W,R)\|_{C[0,T;\dH_k]} \lesssim \| (\W_0,R_0)\|_{\dH_k}, \qquad k = 0,1.
\]
Further this solution is the unique limit of regular solutions, and
depends continuously on the initial data.

\item{(Weak Lipschitz dependence)} The solution $(W,Q)$ has a Lipschitz dependence 
on the initial data in the $\dH_1$ topology.
\end{enumerate}
The implicit constants in all the estimates above depend only on the
initial data size and on the constant $\delta$ in the above pointwise
constraints. Further, the last part of the Theorem~\ref{baiatul} asserts that in
effect these constants can alternatively be estimated purely in terms
of our uniform control parameters $A$ and $B$, rather than the full
Sobolev norm of the data.

Our second result in this paper is a cubic lifespan bound for the small data problem:

\begin{theorem}
\label{t:cubic}
Let $(W,Q)$  be a solution for the system \eqref{ww2d1} whose initial data 
satisfies 
\begin{equation}\label{small-data}
\|(W_0,Q_0)\|_{\dH_0} + \|(\W_0,R_0)\|_{\dH_1} \leq \epsilon \ll 1.
\end{equation}
Then the solution exists for a time $T_{\epsilon} \approx \epsilon^{-2}$,
with bounds 
\begin{equation}\label{small-data-out}
\|(W,Q)(t)\|_{\dH_0} + \|(\W,R)(t)\|_{\dH_1} \lesssim \epsilon, \qquad |t| < T_\epsilon .
\end{equation}
Further, higher regularity is also preserved,
\begin{equation}
 \|(\W,R)(t)\|_{\dH_n} \lesssim  \|(\W,R)(0)\|_{\dH_n} \qquad |t| < T_\epsilon ,
\end{equation}
whenever the norm on the right is finite.
\end{theorem}

To the best of our knowledge this is the first nontrivial lifespan bound for solutions 
to this problem. A similar result for the zero vorticity problem was proved in our 
earlier article \cite{HIT}. The problem here is considerably more difficult than the 
one in \cite{HIT}, both technically and conceptually. At the technical level, 
the normal form for the vorticity problem is much more involved, and quite 
nontrivial to compute (see the next section). Qualitatively, here we have stronger 
quadratic interactions at low frequency, which, unlike in \cite{HIT}, prevent us from 
obtaining cubic bounds for the linearized equation (and thus, for differences of solutions).

We further remark that  in \cite{HIT} we also provide a proof of an almost global result 
 for small localized data for gravity waves. Our aim is to also provide a similar result in 
this context. However, the ideas in \cite{HIT} do not directly carry over to this case,
due first, to the lack of a scaling symmetry, and secondly, to the lack of cubic estimates
for the linearized equation. We hope to be able to address these issues in subsequent 
work.

We note that the periodic case is almost identical and is not discussed separately. 
The only difference in the analysis is in how the constant functions are treated. This
is discussed in detail in the \emph{Appendix} to \cite{HIT}, and carries over to the present paper 
without any change.

\subsection{Outline of the paper}

There are three key steps in our analysis, which eventually provide
all the ingredients which are necessary in order to prove our main results.
These are as follows:
\bigskip

(i) \emph{The normal form analysis.} The constant vorticity water wave equation
has many quadratic interactions, yet we seek to prove small data lifespan bounds
as if the nonlinearities were cubic. At least formally the key to this is the normal form
analysis, which allows  us to replace quadratic nonlinearities with cubic ones.
While the normal form transformation for gravity waves is quite straightforward,
in the presence of constant vorticity, this is no longer the case. 

Indeed, the normal form turns out to be unbounded both at low
frequency and at high frequency.  This computation is fairly involved,
and is carried out in the next section.  Its redeeming feature is that
its outcome is also quite explicit.

The normal form we calculate here is not directly used in any of the estimates we derive
later on. However, it is crucially used in order to construct modified 
energies with cubic estimates, which is the  base of our \emph{quasilinear modified energy
method}. Even though the normal form is badly unbounded, it has 
enough of a ``null structure'', or antisymmetry, so that the cubic energy corrections
it generates are all of bounded type. 

\bigskip

(ii) \emph{The analysis of the linearized equation  in Section~\ref{s:linearized}}. This is a 
critical part of any local well-posedness result for a quasilinear problem.
The derivation of the equations is also interesting, as it clarifies
the quasilinear structure and the roles played by the advection coefficient $\ub$, 
and the frequency shift $\ua$.

As for the gravity waves in \cite{HIT}, we are able to prove that 
the linearized problem is well-posed in our base space $\dH_0$,
Further, the bounds we prove are in terms of our control parameters $\uA$ and $\uB$,
and not in terms the full Sobolev norm of the solution.

Unlike in \cite{HIT},
we are no longer able to prove cubic estimates for the linearized 
equation. This is due to the unbounded low frequency part of the 
normal form transformation, which looses its skew-adjoint structure
after linearization.   Because of this, we are able to use the 
bounds for  the linearized equation in the proof of local well-posedness,
but only partially in the proof of the cubic lifespan result.

\bigskip

 (iii) \emph{The cubic energy estimates in Section~\ref{s:ee}}. Since we already have the conserved 
Hamiltonian, which controls the $\dH_0$ norm of $(W,Q)$, our task here is to 
successively provide bounds for $(\W,R)$ in the $\dH_k$ spaces for $k = 0,1,\cdots$.
In all cases, we control the evolution of these norms using our pointwise control
parameters $\uA$ and $\uB$.

These bounds come in two flavors: (a) local bounds, which apply for large data,
and are needed for the local well-posedness result, and (b) cubic long time bounds
for small data, which are used for the cubic lifespan result. The former are obtained 
largely by differentiating the equation, and then by applying the bounds for the linearized equation.
The latter, however, requires computing cubic energy corrections, and there, we rely heavily 
on the normal form.     

The crucial step is the one for the $\dH_1$ norm of $(\W,R)$ (i.e., $k=1$), as this is the level where 
we have our well-posedness result; for this reason, we describe this case in detail. 
The case $k=0$ is simpler since $(\W,R)$ solves the linearized equation, so we already have
the local bound.  The case $k \geq 2$ is discussed last, without explicitly computing the 
modified energy.  
\bigskip

Once the bounds for the linearized equation and for the differentiated equation are established,
the remaining arguments in the proof of our main theorems are a fairly straightforward 
repetition of arguments in \cite{HIT}. We outline this in the last section of the paper.

Finally, the \emph{Appendix} plays two roles. In \emph{Appendix}~\ref{s:eq} we outline the derivation of the 
constant vorticity gravity wave equation in holomorphic coordinates. Last, but not least,
in \emph{Appendix}~\ref{s:est} we collect a number of bilinear Coifman-Meyer 
and nonlinear Moser type estimates, some from \cite{HIT}, and prove the ones which are new 
in this paper.

\section{ The normal form transformation}
\label{s:nf}

The nonlinear evolution \eqref{ww2d1} contains quadratic terms, yet for our problem
we seek to prove cubic lifespan bounds, as if the nonlinearity is at least cubic.
In this section we consider the question of finding a normal transformation, whose  
aim is to replace the original variables $(W,Q)$ with normal form variables $(\tilde W,\tilde Q)$,
of the form
\begin{equation}\label{nft}
\left\{
\begin{aligned}
\tilde{W} = & W + W_{[2]}
\\
\tilde{Q} = & Q  + Q_{[2]},
\end{aligned}
\right.
\end{equation}
where $W_{[2]}$ and $Q_{[2]}$ are quadratic forms in $(W,Q)$, so that 
the normal form variables satisfy an equation with only cubic and higher terms,
\begin{equation} \label{cubic-eq}
\left\{
\begin{aligned}
& \tilde{W}_t + \tilde Q_\alpha = \text{cubic and higher} \\
& \tilde Q_t -ig \tilde W+ic\tilde Q = \text{cubic and higher} .
\end{aligned} 
\right.
\end{equation}
We will indeed show that such a normal form transformation exists; 
as it turns out, it is highly unbounded, both at low and at high frequencies.
However, this is expected, and it does not cause any difficulties.
This is because we do not use the normal form directly, but only 
as a tool to help us construct modified energies in our  quasilinear modified 
energy method. As it turns, even though the normal  form transformation
is unbounded, it has some favorable structure, so that the modified 
cubic energies are nevertheless bounded.

\subsection{ The resonance analysis}
We begin by examining whether our evolution has quadratic resonant interactions.
Taking into account the possible complex conjugations and the dispersion relation 
\eqref{dispersion}, the  quadratic resonant interactions are associated with three
pairs of characteristic frequencies $(\tau_1,\xi_1)$, $(\tau_2,\xi_2)$ and $(\tau,\xi)$
so that $(\tau_1,\xi_1)+(\tau_2,\xi_2)=(\tau,\xi)$. Combining this with the dispersion relation
we obtain $\tau_1^2 + \tau_2^2 = \tau^2$, which leads to $\tau_1 \tau_2 = 0$.
Hence, we either have $(\tau_1,\xi_1) = (0,0)$ or  $(\tau_2,\xi_2) = (0,0)$.
Thus, resonant interactions occur only when either one of the inputs or the output
is at frequency zero. In terms of the normal form transformation, this indicates
that at most we will have singularities when either input is at frequency zero, or 
the output is at frequency zero. We will see that, due to the form of the quadratic terms
in the equation, the former scenario happens, but the latter does not. 

\subsection{The normal form computation}
 We begin with the quadratic and  expansion in the equation
\eqref{ww2d1}, and then we compute the normal form transformation which
eliminates the quadratic terms from the equation. 
Starting with 
\[
F \approx Q_\alpha - Q_\alpha W_\alpha + \P[\bar Q_\alpha W_\alpha - Q_\alpha \bar W_\alpha] 
+ \P[ (Q_\alpha -\bar Q_\alpha) (4|\Re W_\alpha|^2- |W_\alpha|^2)] ,
\]
we compute  the multilinear expansion:
\begin{equation} \label{ww-multi}
\left\{
\begin{aligned}
& W_t + Q_\alpha = G^{(2)} + G^{(3+)} \\
& Q_t -ig W+icQ = K^{(2)}  + K^{(3+)},
\end{aligned} 
\right.
\end{equation}
where the quadratic terms $(G^{(2)},K^{(2)})$ are given by
\[
\left\{
\begin{aligned}
\mathbf{G}^{(2)}&=-\P\left[ \bar{Q}_{\alpha}W_{\alpha}-Q_{\alpha}\bar{W}_{\alpha}\right]-i\frac{c}{2}\P\left[ W\bar{W}_{\alpha}+\bar{W}W_{\alpha}\right]+i\frac{c}{2}WW _{\alpha}\\ 
\mathbf{K}^{(2)}&=-Q_{\alpha}^2-\P\left[ \vert Q_{\alpha}\vert ^2\right]+i\frac{c}{2}WQ_{\alpha}+i\frac{c}{2}\P\left[ W\bar{Q}_{\alpha}-\bar{W}Q_{\alpha}\right].
\end{aligned}
\right.
\]
The role of the normal form transformation is to eliminate the quadratic terms $(G^{(2)},K^{(2)})$
from the equation \eqref{ww-multi}. We can divide the quadratic terms above into two classes:

a) holomorphic, i.e., those which are the product of two holomorphic functions,

b) mixed, i.e., those which are the (projected) product of one holomorphic function with 
another's conjugate.  
 
Thus, we expect the normal form to have a similar structure,
\begin{equation}\label{nft-ha}
 W_{[2]} =  W_{[2]}^h +  W_{[2]}^a, \qquad  Q_{[2]} =  Q_{[2]}^h +  Q_{[2]}^a,
\end{equation}
where, allowing for all possible combinations, the above components must have the form
\begin{equation}\label{nft-ha+}
\begin{split}
W_{[2]}^h = & \   B^h(W,W) + C^h(Q,Q)+ D^h(W,Q),
\\
W_{[2]}^a = & \ B^a(W,\bar W) + C^a(Q,\bar Q)+D^{a}(W, \bar Q)+ E^a ( Q, \bar{W}),
\\
 Q_{[2]}^h = & \ F^h(W,W)+H^{h}(Q,Q) +A^h(W,Q),
\\
  Q_{[2]}^a =& \ F^a(W,\bar W) + H^a(Q,\bar Q)+A^a(W,\bar Q) +G^a ( Q, \bar W).
\end{split}
\end{equation}
All the above expressions are translation invariant bilinear forms, whose 
symbols we need to compute. Our main result is as follows:

\begin{proposition}\label{p:nf}
The equation \eqref{ww-multi} admits  a normal form transformation \eqref{nft} 
as in \eqref{nft-ha}-\eqref{nft-ha+}, where the  symbols for the bilinear operators
in \eqref{nft-ha+} are given by \eqref{sol}, \eqref{sol1}, and \eqref{long-nft1}.
\end{proposition}

For later use, i.e., for computing the modified energy cubic corrections, we also translate 
the symbols \eqref{sol}, \eqref{sol1}, and \eqref{long-nft1} into  the spatial description.
Precisely, the holomorphic terms have the form:
\begin{equation}
\left\{
\begin{aligned}
&B^h(W,W)=-W\W+i\dfrac{c^2}{2g}(\W \partial^{-1}_{\alpha}W+W^2)+
\dfrac{c^4}{4g^2}W\partial^{-1}_{\alpha}W, 
\\
& C^h(Q,Q)=-\dfrac{c^2}{4g^2}QQ_{\alpha},
 \\
&  D^h(W,Q)=-\dfrac{c}{2g}(\W Q+WQ_{\alpha})+i\dfrac{c^3}{4g}(WQ+\partial^{-1}_{\alpha}W Q_{\alpha}), 
\\
& F^h(W,W)=i\dfrac{c}{4}W^2+\dfrac{c^3}{4g}W\partial^{-1}_{\alpha}W, 
\\
& H^{h}(Q,Q)=-\dfrac{c}{2g}QQ_{\alpha}, 
\\
& A^{h}(W,Q)=-WQ_{\alpha}+i\dfrac{c^2}{2g}\partial^{-1}_{\alpha}WQ_{\alpha}+i\dfrac{c^2}{4g}WQ,
\end{aligned}
\right.
\end{equation}
and the antiholomorphic counterparts:
\begin{equation}
\left\{
\begin{aligned}
&B^a(W,\bar{W})=-\W\bar{W}-i\dfrac{c^2}{2g}\W\partial^{-1}_{\alpha}\bar{W}+i\dfrac{c^2}{4g}\vert W\vert^2
-\dfrac{c^4}{4g^2}W \partial^{-1}_{\alpha}\bar{W}, 
\\
&  C^a(Q,\bar{Q})=-\dfrac{c^2}{4g^2}\bar{Q}Q_{\alpha},
 \\
&  D^a(W,\bar{Q})=-\dfrac{c}{2g}\W \bar{Q}+i\dfrac{c^3}{4g^2}W \bar{Q}, 
\\ 
& G^a(Q, \bar{W})=-Q_{\alpha}\bar{W} -i\dfrac{c^2}{2g}Q_{\alpha}\partial^{-1}_{\alpha}\bar{W}
\\
&E^a(Q,\bar{W})=-\dfrac{c}{2g}Q_{\alpha}\bar{W} -i\dfrac{c^3}{4g^2}Q_{\alpha}\partial_{\alpha}^{-1}\bar{W}, 
\\
 & H^{a}(Q,\bar{Q})=-\dfrac{c}{2g}\bar{Q}Q_{\alpha}, 
\\
& A^{a}(W,\bar{Q})=i\dfrac{c^2}{4g}W\bar{Q}, 
\\
& 
 F^a(W,\bar{W})=i\dfrac{c}{2}\vert W\vert^2-\dfrac{c^3}{4g}W\partial^{-1}_{\alpha}\bar{W}.
\end{aligned}
\right.
\end{equation}
We note that while for computational purposes it is convenient to separate the two
components of the normal form, in order to see the antisymmetric structure of the 
(low frequency) unbounded  part, one has to consider them together; see the computations in 
Section~\ref{s:ee}.  Toward that goal, we rewrite the quadratic normal form components in the following form: 
\begin{equation}\label{nf-eq}
\begin{split}
W^{[2]} = & \ - (W+\bar W) W_\alpha - \frac{c}{2g} \left[ (Q+\bar Q) W_\alpha + (W + \bar W) Q_\alpha\right] 
\\ & \ + \frac{ic^2}{2g} \left[ (\partial^{-1} W - \partial^{-1} \bar W) W_\alpha + W^2 +\frac12 |W|^2 \right]
- \frac{c^2}{4g^2} (Q+\bar Q) Q_\alpha 
\\ & \ + \frac{ic^3}{4g^2} \left[ (Q+\bar Q) W + (\partial^{-1} W - \partial^{-1} \bar W) Q_\alpha \right] + 
 \frac{c^4}{4g^2} (\partial^{-1} W - \partial^{-1} \bar W) W,
\\
Q^{[2]} = & \ - (W+\bar W) Q_\alpha - \frac{c}{2g} (Q+\bar Q) Q_\alpha + \frac{ic}{4} (W^2+ 2|W|^2)
\\ & \ + \frac{ic^2}{2g} \left[ (\partial^{-1} W - \partial^{-1} \bar W) Q_\alpha + \frac12  (Q+\bar Q) W \right]
+  \frac{c^3}{4g} (\partial^{-1} W - \partial^{-1} \bar W) W .
\end{split}
\end{equation}
\begin{proof}
A-priori, computing the normal form, i.e., all  symbols of the bilinear expressions above, might appear
quite involved. However, there are several observations which bring this analysis to a more 
manageable level:
\begin{itemize}
\item The analysis for the holomorphic products, and for the mixed terms is completely separate.
\item The system we obtain for the symbols has polynomial coefficients, so the solutions 
are rational functions. 
\item Counting $c$ as one half of a derivative; the problem is homogeneous. Thus, organizing 
symbols based on the powers of $c$,  each such term will have a specific homogeneity. 
Further, at each power of $c$, we will encounter only half the terms, as there is a half derivative 
difference between the scaling of $W$ and that of $Q$.
\item The terms without $c$ are already known from the gravity wave problem \cite{HIT}.
\end{itemize}
Given the above considerations, the natural strategy is to split the analysis into the holomorphic
and the mixed part, and in each of these cases to successively solve for increasing powers 
of $c$. The computation stops at $c^4$.
\smallskip

{\bf (i) Holomorphic terms:} Here we seek a normal form  for  the system
\begin{equation*}
\left\{
\begin{aligned}
& W_t + Q_\alpha =i\frac{c}{2}WW_{\alpha} \\
& Q_t-igW+icQ=  -  Q_\alpha^2 +i\frac{c}{2}WQ_\alpha
  + cubic.
\end{aligned} 
\right.
\end{equation*}
By checking parity, our normal form 
must be
\[
\left\{
\begin{aligned}
&\tW = W + B^h(W,W) + C^h(Q,Q)+ D^h(W,Q) \\
& \tQ = Q +F^h(W,W)+H^{h}(Q,Q) +A^h(W,Q) ,
\end{aligned}
\right.
\]
where $B^h, \ C^h, \ F^h$ and $H^h$ are symmetric bilinear forms with
symbols $B^h(\xi,\eta)$, $C^h(\xi,\eta)$, $F^h(\xi,\eta)$,
$H^h(\xi,\eta)$ and $A^h$ and $D^H$ are arbitrary. We compute
\[
\left\{
\begin{aligned}
&\tW_{t} +\tQ_{\alpha} =-2 B^h(Q_{\alpha},W) + 2igC^h(W,Q)-2icC^h(Q,Q)-D^h(Q_{\alpha},Q)+igD^h(W,W)\\
&\hspace{2cm}-ic D^h(W,Q)+\partial_{\alpha}\left( F^h(W,W)+H^{h}(Q,Q) +A^h(W,Q)\right) 
+i\frac{c}{2}WW_{\alpha}+\mbox{cubic}, \\
& \tQ_t-ig\tW+ic\tQ= -2F^h(Q_{\alpha},W)+2igH^h(W,Q)-icH^{h}(Q,Q) -A^h(Q_{\alpha},Q)+igA^h(W,W)\\
&\hspace{3.4cm} -igB^h(W,W)-igC^h(Q,Q)-igD^h(W,Q)+icF^h(W,W) -Q^2_{\alpha}\\
&\hspace{3.4cm}+i\frac{c}{2}WQ_{\alpha}+\mbox{cubic} .
\end{aligned}
\right.
\]
We denote the two input frequencies by $\xi$ and $\eta$, both of which are negative. 
Then, matching like terms, we obtain the following linear system for the symbols:
\[
\left\{
\begin{aligned}
&2\eta B^h -2gC^h+cD^h-(\xi +\eta)A^h=0\\
& 2cC^h+\left[ \xi D^h \right]_{sym}- (\xi +\eta)H^h=0  \\
&g\left[ D^h\right]_{sym}+(\xi +\eta)F^h=-i\frac{c}{4}(\xi+\eta) \\
&2\eta F^h-2gH^h+gD^h=i\frac{c}{2}\eta \\
&cH^h+\left[ \xi A^h \right] _{sym}+gC^h=-i\xi\eta\\
&g\left[ A^h\right]_{sym}-gB^h +cF^h=0,
\end{aligned}
\right.
\]
where $sym$ stands for symmetrization.
Using the first equation in the system above helps us to determine the symmetrized symbol of $A^h$ as follows
\[
A^h=\frac{1}{\xi +\eta}\left[ cD^h -2gC^h+2\eta B^h\right] ,
\]
which implies 
\begin{equation}
\label{A}
\begin{aligned}
&\left[ A^h\right] _{sym}=\frac{1}{\xi +\eta}\left[ c\left[ D^h \right] _{sym}-2gC^h+(\eta+\xi) B^h\right] \\
&\left[ \xi A^h\right] _{sym}=\frac{1}{\xi +\eta}\left[ c\left[ \xi D^h \right] _{sym}-g(\xi +\eta )C^h+2\xi\eta B^h\right].
\end{aligned}
\end{equation}
The forth equation provides an expression for the symmetrized symbol of $D^h$ as follows
\[
D^h=-\frac{2}{g}\eta F^h +2H^h+i\frac{c}{2g}\eta,
\]
which implies 
\begin{equation}
\label{D}
\begin{aligned}
&g\left[ D^h\right] _{sym}=-(\xi+\eta) F^h +2gH^h+i\frac{c}{4}(\xi +\eta),\\
& \left[ \xi D^h\right] _{sym}=-\frac{2}{g}\xi\eta F^h +(\xi +\eta)H^h+i\frac{c}{2g}\xi\eta,\\
\end{aligned}
\end{equation}
Using \eqref{D} in \eqref{A}  gives 
\begin{equation}
\label{nouA}
\begin{aligned}
&\left[ A^h\right] _{sym}=\frac{1}{\xi +\eta}\left[ -\frac{c}{g}(\xi+\eta) F^h +2cH^h+i\frac{c^2}{4g}(\xi +\eta)-2gC^h+(\eta+\xi) B^h\right] \\
&\left[ \xi A^h\right] _{sym}=\frac{1}{\xi +\eta}\left[ -\frac{2c}{g}\xi\eta F^h +(\xi +\eta)cH^h+i\frac{c^2}{2g}\xi \eta -g(\xi +\eta )C^h+2\xi\eta B^h\right].
\end{aligned}
\end{equation}

Thus, we return to the following system
\[
\left\{
\begin{aligned}
& 2cC^h+\left[ \xi D^h \right]_{sym}- (\xi +\eta)H^h=0  \\
&g\left[ D^h \right]_{sym}+(\xi +\eta)F^h=-i\frac{c}{4}(\xi+\eta) \\
&cH^h+\left[ \xi A^h\right] _{sym}+gC^h=-i\xi\eta\\
&g\left[ A^h\right]_{sym}-gB^h +cF^h=0,
\end{aligned}
\right.
\]
where we substitute the corresponding values from \eqref{D}, \eqref{nouA}, and obtain
\[
\left\{
\begin{aligned}
& 2cC^h-\frac{2}{g}\xi\eta F^h +(\xi +\eta)H^h+i\frac{c}{2g}\xi\eta- (\xi +\eta)H^h=0  \\
&-(\xi+\eta) F^h +2gH^h+i\frac{c}{4}(\xi +\eta)+(\xi +\eta)F^h=-i\frac{c}{4}(\xi+\eta) \\
&cH^h+\frac{1}{\xi +\eta}\left[ -\frac{2c}{g}\xi\eta F^h +(\xi +\eta)cH^h+i\frac{c^2}{2g}\xi \eta -g(\xi +\eta )C^h+2\xi\eta B^h\right]+gC^h=-i\xi\eta\\
&\frac{1}{\xi +\eta}\left[ -c(\xi+\eta) F^h +2cgH^h+i\frac{c^2}{4}(\xi +\eta)-2g^2C^h+(\eta+\xi)g B^h\right]-gB^h +cF^h=0,
\end{aligned}
\right.
\]
so
\[
\left\{
\begin{aligned}
& 2cC^h-\frac{2}{g}\xi\eta F^h =-i\frac{c}{2g}\xi\eta  \\
&H^h=-i\frac{c}{4g}(\xi+\eta) \\
&cH^h -\frac{c}{g}\frac{\xi\eta}{\xi +\eta} F^h +\frac{\xi\eta}{\xi +\eta} B^h=-i\frac{\xi\eta}{2}-i\frac{c^2}{4g}\frac{\xi \eta}{\xi +\eta}\\
&cH^h-gC^h=-i\frac{c^2}{8g}(\xi+\eta).
\end{aligned}
\right.
\]
The solution is
\begin{equation}
\label{sol}
\begin{aligned}
&B^h=-i\frac{1}{2}(\xi +\eta)+i\frac{c^2}{4g}\frac{(\xi +\eta)^2}{\xi\eta}-i\frac{c^4}{8g^2}\frac{\xi +\eta}{\xi\eta},\\
C^h=-i&\frac{c^2}{8g^2}(\xi +\eta), \quad F^h=i\frac{c}{4} -i\frac{c^3}{8g}\frac{\xi +\eta}{\xi \eta},\quad H^h=-i\frac{c}{4g}(\xi+\eta).\\
\end{aligned}
\end{equation}
It remains to find $A^h$ and $D^h$, which we obtain from
\[
\left\{
\begin{aligned}
&2\eta B^h -2gC^h+cD^h-(\xi +\eta)A^h=0\\
&2\eta F^h-2gH^h+gD^h=i\frac{c}{2}\eta .\\
\end{aligned}
\right.
\]
Thus, 
\begin{equation}
\label{sol1}
\begin{aligned}
A^h=-i\eta+i\frac{c^2}{2g}\eta \xi^{-1}+i\frac{c^2}{4g}, \quad D^h=i\frac{c^3}{4g^2} \frac{\xi +\eta}{\xi }-i\frac{c}{2g}(\xi +\eta).
\end{aligned}
\end{equation}

\bigskip 

{\bf (ii) Mixed terms:} Here we need a normal form for the system:
\begin{equation*}
\left\{
\begin{aligned}
& W_t + Q_\alpha = -\P\left[ \bar{Q}_{\alpha}W_{\alpha}-Q_{\alpha}\bar{W}_{\alpha}\right]-i\frac{c}{2}\P\left[ W\bar{W}_{\alpha}+\bar{W}W_{\alpha}\right] + cubic \\
&  Q_t-igW+icQ=   -  \P\left[ \vert Q_{\alpha}\vert ^2\right]+i\frac{c}{2}\P\left[ W\bar{Q}_{\alpha}-\bar{W}Q_{\alpha}\right]  + cubic.
\end{aligned} 
\right.
\end{equation*}
The general expression for our normal form is
\begin{equation}
\label{long-nft}
\left\{
\begin{aligned}
\tW &= W + B^a(W,\bar W) + C^a(Q,\bar Q)+D^{a}(W, \bar Q)+ E^a ( Q, \bar{W})\\
 \tQ &= Q + F^a(W,\bar W) + H^a(Q,\bar Q)+A^a(W,\bar Q) +G^a ( Q, \bar W),
\end{aligned}
\right.
\end{equation}
where no symmetry assumption is required.
We compute 
\[
\left\{
\begin{aligned}
&\tW_t +\tQ_{\alpha}= - B^a(Q_{\alpha},\bar W) -B^a(W,\bar Q_{\alpha})+ig C^a(W,\bar Q)-igC^a(Q,\bar W)-D^{a}(Q_{\alpha}, \bar Q)\\ 
& \hspace*{2cm}  -igD^a(W,\bar W) +icD^a(W, \bar{Q})+ igE^a ( W, \bar{W})-icE^a(Q,\bar W)-E^a(Q, \bar Q_{\alpha})\\ 
& \hspace*{2cm} -\P\left[ \bar{Q}_{\alpha}W_{\alpha}-Q_{\alpha}\bar{W}_{\alpha}\right] + \partial_{\alpha}\left[ F^a(W,\bar W) + H^a(Q,\bar Q)+A^a(W,\bar Q) +G^a ( Q, \bar W) \right] 
\\ 
&\hspace*{2cm}  -i\frac{c}{2}\P\left[ W\bar{W}_{\alpha}+\bar{W}W_{\alpha}\right] + cubic \\
& \tQ_t -ig\tW+ic\tQ = -F^a(Q_{\alpha},\bar W) -F^a(W, \bar Q_{\alpha})+ig H^a(W,\bar Q)-igH^a(Q,\bar W)-A^a(Q_{\alpha},\bar Q)\\
 &\hspace*{3.4cm} -igA^a(W,\bar W) +icA^a(W, \bar Q)-G^a(Q, \bar Q_{\alpha})+igG^a(W, \bar W)-igB^a(W, \bar W)\\
 &\hspace*{3.4cm}-igC^a(Q, \bar Q)-igD^a(W, \bar Q)-ig E^a (Q, \bar W)+icF^a(W, \bar W)+icH^a (Q, \bar{Q})\\
&\hspace*{3.4cm}+icA^a (W, \bar{Q}) -  \P\left[ \vert Q_{\alpha}\vert ^2\right]+i\frac{c}{2}\P\left[ W\bar{Q}_{\alpha}-\bar{W}Q_{\alpha}\right]  + cubic.
\end{aligned}
\right.
\]
Now we denote by $\xi$ the frequency of the holomorphic input and by $\eta$ the frequency of the 
conjugated input (both negative). Matching again like terms, it remains to solve the following system 
\begin{equation}
\label{sistem}
\left\{
\begin{aligned}
&\xi B^a+gC^a+cE^a-(\xi -\eta) G^a=-i\xi\eta\\
&\eta B^a+gC^a+(\xi -\eta) A^a+cD^a=-i\xi\eta\\
&\xi D^a-\eta E^a-(\xi -\eta) H^a=0\\
&g D^a-g E^a-(\xi -\eta) F^a=i\frac{c}{2}(\eta -\xi)\\
&\xi F^a+g H^a+gE^a=-i\frac{c}{2}\xi\\
&\eta F^a+g H^a+2cA^a-gD^a=i\frac{c}{2}\eta\\
&\xi A^a+g C^a-cH^a-\eta G^a=i\xi \eta \\
&g A^a+g B^a-cF^a-g G^a=0.\\
\end{aligned}
\right.
\end{equation}
To  avoid solving an $8\times 8$ system we expand the result in terms of powers of $c$.
A homogeneity analysis shows that $B^a$, $C^a$, $G^a$ and $A^a$ contain only even powers,
while the other four symbols contain only odd powers. We solve for the coefficients of increasing
powers of $c$, noting that at each stage we only need to solve a $4 \times 4$ system.  
The outcome of the first step (i.e., if $c=0$) is already known from gravity waves \cite{HIT}.
The computations are somewhat tedious, but quite elementary. We omit them and only
write the final result:
\begin{equation}
\label{long-nft1}
\left\{
\begin{aligned}
&C^a(\xi, \eta)=-i\frac{c^2}{4g^2}\xi, \qquad  \quad B^a(\xi, \eta)=-i\frac{c^4}{4g^2}\eta^{-1}+i\frac{c^2}{2g}\xi\eta^{-1}+i\frac{c^2}{4g}-i\xi ,  \\
&D^a(\xi, \eta)=i\frac{c^3}{4g^2} -i\frac{c}{2g}\xi,  \quad E^{a}(\xi, \eta)=i\frac{c^3}{4g^2} \xi\eta^{-1}-i\frac{c}{2g}\xi,\\
&A^a(\xi, \eta)= i\frac{c^2}{4g} ,\qquad  \quad \quad \ \  F^a(\xi, \eta)=-i\frac{c^3}{4g}\eta^{-1}+i\frac{c}{2}, \\
&H^a(\xi, \eta)=-i\frac{c}{2g}\xi,\quad  \qquad  G^a(\xi, \eta)=i\frac{c^2}{2g}\xi \eta^{-1}-i\xi .
\end{aligned}
\right.
\end{equation}
\end{proof}

\section{The linearized equation}
\label{s:linearized}

In this section we derive the linearized water wave equations,  and
prove energy estimates for  them.  We recall that  for the similar problem in \cite{HIT}
we were able to prove quadratic energy estimates in $\dH_0$, which apply 
for the large data problem, but also cubic energy estimates for the small data problem.
By contrast, here  we are only able to prove the quadratic energy estimates.
 This suffices for the local well-posedness theory, but is not so useful 
in order to establish improved lifespan bounds. It appears unlikely that 
cubic energy estimates hold in $\dH_0$  for the  linearized problem; in any case,
we leave this question open.

\subsection{Computing the linearization.} The solutions for the
linearized water wave equation around a profile $\left( W,Q\right) $
are denoted by $(w,q)$. However, it will be more convenient to
immediately switch to diagonal variables $(w,r)$, where
\[
r := q - Rw.
\]
We first recall the equations,
\begin{equation}
\label{ww2d1-re}
\left\{
\begin{aligned}
&W_t+ \underline{F}(W_{\alpha}+1) +i\frac{c}{2}W=0\\
&Q_t-igW+\underline{F} Q_{\alpha}+icQ+\mathbf{P}\left[|R|^2\right]-i\frac{c}{2} T_1=0,
\end{aligned}
\right.
\end{equation}
where $\underline{F} = F - i \dfrac{c}2 F_1$ with
\[
F = \P\left[\frac{R}{1+\bar \W} - \frac{\bar R}{1+ \W}\right],
\qquad 
F_1 = \P \left[ \frac{W}{1+\bar \W} +  \frac{\bar W}{1+\W}\right],
\qquad T_1 = \P[W \bar R - \bar W R].
\]
The linearization of $R$ is 
\[
\delta R =
\dfrac{q_{\alpha}- Rw_{\alpha}}{1+\W}
= \dfrac{r_{\alpha}+ R_\alpha w}{1+\W},
\]
and that of $|R|^2$ is
\[
\delta |R|^2 = n+ \bar n, \qquad n := \bar R \delta R =
 \frac{ \bar R(r_{\alpha}+R_\alpha w)}{1+\W}.
\]

The linearizations of $F$, $F_1$ and $T_1$ can be expressed in the form
\[
\delta F = \P[ m - \bar m], \qquad \delta F_1 = \P[ m_1 + \bar m_1], \qquad 
\delta T_1 = \P\left[m_2-\bar{m}_2 \right],
\]
where the auxiliary variables $m,m_1,m_2$ correspond to differentiating
$F$, $F_1$ and $T_1$ with respect to the holomorphic variables,
\[
\begin{split}
m := & \ \frac{q_\alpha - R w_\alpha}{J} + \frac{\bar R w_\alpha}{(1+\W)^2} =
 \frac{r_\alpha +R_\alpha w}{J} + \frac{\bar R w_\alpha}{(1+\W)^2},
\\
 m_1:=& \ \frac{1}{1+\bar{\W}}w-\frac{\bar{W}}{(1+\W)^2}w_{\alpha}, \ \ \ m_2:=  \ \bar{R}w-\dfrac{\bar{W} r_{\alpha}+ \bar{W} R_\alpha w}{1+\W}.
\end{split}
\]

Given all of the above, it follows that the linearized water wave equations take the form
\begin{equation*}
\left\{
\begin{aligned}
&w_{t}+ \underline{F} w_\alpha + (1+ \W) ( \P[ m-\bar m] - i\frac{c}{2} \P[m_1+\bar m_1]) +i\frac{c}{2}w  =0\\
&q_{t}+  \underline{F} q_\alpha + Q_\alpha( P[ m-\bar m] - i\frac{c}{2} \P[m_1+\bar m_1])   -ig w +icq+\P\left[n+\bar n\right]  -i\frac{c}{2} \P[m_2-\bar m_2]=0 .
\end{aligned}
\right.
\end{equation*}
Now we transition in both equations from $\underline F$ to the real advection coefficient 
$\underline b$ using the relations \eqref{ftob}.
We also move to the right all the terms which we expect to be
perturbative.  These are terms like $\bar P m$, $\bar Pn$, $\bar P
m_1$, $\bar Pm_2$, which are lower order since the differentiated holomorphic
variables have to be lower frequency. The same applies to their
conjugates.  Then, our equations are rewritten as
\begin{equation*}
\left\{
\begin{aligned}
  &\left(\partial_t + \underline{b} \partial_\alpha \right) w
  +\left[ i\frac{c}{2}w + (1+\W)(m-i\frac{c}2 m_1) - \dfrac{\bar R
    w_\alpha}{1+\W} - i\frac{c}2\frac{\bar{W}w_{\alpha}}{1+\W}\right] =
  \underline{\mathcal{G}}_0
  \\
  & \left(\partial_t + \underline{b} \partial_\alpha \right) q -ig
  w+icq +\left[ Q_\alpha (m-i\frac{c}2 m_1) + n - i\frac{c}{2} m_2 -
  \dfrac{\bar R q_\alpha}{1+\W}-
  i\frac{c}2\frac{\bar{W}q_{\alpha}}{1+\W} \right]
=\underline{\mathcal{K}}_0 ,
\end{aligned}
\right.
\end{equation*}
where 
\[
\left\{
\begin{aligned}
\underline{\mathcal{G}}_0 = &\ (1+\W)  \left((\P\bar  m+ \bar \P  m) + i\frac{c}{2}(\P \bar m_1-\bar \P {m}_1) \right) 
\\
\underline{\mathcal{K}}_0 = &\ Q_\alpha  \left((\P\bar  m+ \bar \P  m) + i\frac{c}{2}(\P \bar m_1-\bar \P {m}_1) \right) + (\bar \P n - \P \bar n) - i\frac{c}2 (\P \bar m_2 + \bar \P m_2)  .
\end{aligned}
\right.
\]
Taking advantage of algebraic cancellations in the square brackets above we are left with
\begin{equation*}
\left\{
\begin{aligned}
&\left(\partial_t +  \underline{b}  \partial_\alpha \right)  w+
\frac{1}{1+\bar{\W}} \left( r_{\alpha}+ R_{\alpha} w +i\frac{c}2(\bar{\W}-\W) w \right) 
 = \underline{ \mathcal{G}}_0\\
& \left(\partial_t + \underline{b} \partial_\alpha \right) q -ig
  w+icq +\frac{R}{1+\bar{\W}} \left( r_{\alpha}+ R_{\alpha} w +i\frac{c}2(\bar{\W}-\W) w \right) 
 - i \frac{c}2 (R + \bar R) w 
=\underline{\mathcal{K}}_0 .
\end{aligned}
\right.
\end{equation*}
Now  we can  switch from $q$ to $r= q-Rw$  and obtain a diagonalized system,
namely 
\begin{equation*}
\left\{
\begin{aligned}
&\left(\partial_t +  \underline{b}  \partial_\alpha \right)  w+
\frac{1}{1+\bar{\W}} \left( r_{\alpha}+ R_{\alpha} w +i\frac{c}2(\bar{\W}-\W) w \right) 
 = \underline{ \mathcal{G}}_0\\
& \left(\partial_t + \underline{b} \partial_\alpha \right) r +icr- i \left( g + \frac{c}2 (\bar R -R)  + i  (\partial_t +  \underline{b}  \partial_\alpha ) R \right) w
=\underline{\mathcal{K}}_0 - R \underline{\mathcal{G}}_0.
\end{aligned}
\right.
\end{equation*}
For the expression $ (\partial_t +  \underline{b}  \partial_\alpha ) R$ in the second equation
we use \eqref{ww2d-diff} and compute the coefficient of $w$ as follows:
\[
\begin{split}
  g  +  \frac{c}2 ( \bar R-R) -i  (\partial_t +  \underline{b}  \partial_\alpha ) R 
= & \ g  +  \frac{c}2 ( \bar R-R) + cR- \frac{g\W-a}{1+\W} - \frac{c}{2}\frac{R\W + 
\bar R \W  +N}{1+\W}
\\
= & \ \frac{ g + a + c/2(R+\bar R - N) }{1+\W}.
\end{split}
\]
This motivates the definition of $\underline{a}$ in \eqref{defa}. 
With this notation, we write the final form of the linearized equations as
\begin{equation}
\label{lin(wr)0}
\left\{
\begin{aligned}
&\left(\partial_t +  \underline{b}  \partial_\alpha \right)  w+\frac{r_\alpha}{1+\bar{\W}}+\frac{R_{\alpha} w}{1+\bar{\W}} = \underline{ \mathcal{G}}(w,r)\\
& \left(\partial_t + \underline{b} \partial_\alpha \right) r +icr  -i\frac{g+\underline{a}}{1+\W}w=\underline{\mathcal{K}}(w,r), 
\end{aligned}
\right.
\end{equation}
where $\underline{\mathcal{G}}(w,r)$, and $\underline{\mathcal{K}}(w,r)$ are given by
\[
\underline{\mathcal{G}}(w,r) = \mathcal{G}(w,r)-i\frac{c}{2}\mathcal{G}_1(w,r),
\qquad \underline{\mathcal{K}}(w,r) = \mathcal{K}(w,r)-i\frac{c}{2}\mathcal{K}_1(w,r),
\]
where 
\begin{equation}
\label{defG1}
\mathcal{G}(w,r) = (1+\W) (\P\bar m+ \bar \P m),  \qquad
\mathcal{G}_1(w,r) = - (1+\W) (\P \bar m_1-\bar \P {m}_1) +
\frac{(\bar{\W}-\W)w}{1+\bar{\W}},
\end{equation}
\begin{equation}
\label{defK1}
\mathcal{K}(w,r) = \bar \P n - \P \bar n,  \qquad
\mathcal{K}_1(w,r) =  \P \bar m_2 + \bar \P m_2.
\end{equation}
If $c = 0$, then these equations coincide with those in
in \cite{HIT}. The fact that $g+ \underline{a}$ is real and 
positive is crucial for the well-posedness of the linearized system.
 
We remark that while $(w,r)$ are holomorphic, it is not directly obvious
that the above evolution preserves the space of holomorphic states. To
remedy this one can also project the linearized equations onto the
space of holomorphic functions via the projection $\P$.  Then we obtain
the equations
\begin{equation}
\label{lin(wr)}
\left\{
\begin{aligned}
&\left(\partial_t +  \M_{\underline{b}}  \partial_\alpha \right)  w+\P\left[ \frac{1}{1+\bar{\W}}r_{\alpha}\right] +\P\left[ \frac{R_{\alpha}}{1+\bar{\W}}w\right]  =  \P\underline{\mathcal{G}}(w,r)\\
& \left(\partial_t + \M_{\underline{b}} \partial_\alpha \right) r +icr-i\P\left[ \frac{g+\underline{a}\, }{1+\W}w\right]=\P\underline{\mathcal{K}}(w,r).
\end{aligned}
\right.
\end{equation}
Since the original set of equations \eqref{ww2d1} is fully
holomorphic, it follows that the two sets of equations,
\eqref{lin(wr)0} and \eqref{lin(wr)}, are algebraically equivalent.

In order to investigate the possibility of cubic linearized energy
estimates it is also of interest to separate the quadratic parts $
\underline{\mathcal{G}}^{2}$ and $\underline{\mathcal{K}}^{2}$ of
$\underline{\mathcal{G}}$ and $\underline{\mathcal{K}}$.
The holomorphic quadratic parts of $ \mathcal{G}$ and $\mathcal K$, 
which also appear in \cite{HIT}, are given by 
\begin{equation*}
 \P \mathcal{G}^{(2)}(w,r) =-P \left[ \W \bar r_\alpha
  \right] + P\left[ R \bar w_\alpha\right] , \qquad  \P \mathcal{K}^{(2)}(w,r)  = -\P[R \bar r_\alpha] .
\end{equation*}
Next we compute the similar decomposition for $\mathcal{G}^{2}_1$ and
$\mathcal{K}^{2}_1$, since the rest was done in \cite{HIT}.  These are
split into quadratic and higher terms as shown below
\[
\begin{split}
   \mathcal{G}_1 = &  \,  \mathcal{G}^{(2)}_1+  \mathcal{G}^{(3+)}_1,
  \qquad   \mathcal{K}_1 = \  \mathcal{K}^{(2)}_1+  \mathcal{K}^{(3+)}_1.
\end{split}
\]
For the quadratic parts we have the holomorphic components
\begin{equation}
\label{quadratic}
\begin{aligned}
 \P \mathcal{G}^{(2)}_1(w,r) = & \ \P \left[\mathbf{\W} \bar{w}\right] + \P\left[W \bar{w}_{\alpha} \right] + \P\left[\bar \W {w} \right] - \P\left[\W {w} \right], 
\\ 
  \P \mathcal{K}^{(2)}_1(w,r) =  &\ \P \left[ W \bar r_\alpha\right]-\P \left[ R\bar{w}\right] ,
\end{aligned}
\end{equation}
and the antiholomorphic components
\[
\begin{aligned}
  \bar \P \mathcal{G}^{(2)}_1(w,r) = & \ -\bar \P \left[\bar \W
    {w}\right] - \bar \P\left[\bar W {w}_{\alpha} \right] +\bar
  \P\left[\bar \W {w} \right],
  \\
  \bar \P \mathcal{K}^{(2)}_1(w,r) = &\ \bar \P \left[ \bar W
    r_\alpha\right]- \bar \P \left[ \bar R {w}\right].
\end{aligned}
\]

The cubic terms have the form
\[
\begin{split}
  \mathcal{G}_1^{(3+)}(w,r) = & - \W ( \P \bar m_1 - \bar \P m_1)  - 
( \P \bar m_1^{(3+)}  - \bar \P m_1^{(3+)})  - \bar Y( \bar \W - \W) w ,
\\
  \mathcal{K}_1^{(3+)}(w,r) = &  \   \bar \P m_2^{(3+)} + \P \bar m_2^{(3+)}.
\end{split}
\]
For the purpose of simplifying nonlinear estimates, it is convenient
to express $\mathcal G^{(3)}_1$ and $\mathcal K^{(3)}_1$ in a polynomial
fashion. This is done using the variable $Y=\dfrac{\W}{1+\W}$. Then we have
\[
\begin{aligned}
  \bar \P m_1 =  \ - \bar \P [ \bar{Y}w +(Y^2-2Y+1)\bar{W}w_{\alpha}],\quad \bar \P m^{(3+)} =  \ \bar \P [ (2Y-Y^2)\bar{W}w_{\alpha}],
  \end{aligned}
  \]
  \[
  \begin{aligned}
 \bar \P m_2^{(3+)} = \ \bar \P [\bar{W}Yr_{\alpha}  -\bar{W} R_\alpha  (1-Y) w].
\end{aligned}
\]

\subsection{Quadratic estimates for large data.}  Our goal here is to
study the well-posedness of the system \eqref{lin(wr)} in $L^2 \times
\dot H^\frac12$.  We begin with a more general version of the system
\eqref{lin(wr)}, namely
\begin{equation}\label{lin(wr)inhom}
\left\{
\begin{aligned}
&\left(\partial_t +  \M_{\underline{b}}  \partial_\alpha \right)  w+\P\left[ \frac{1}{1+\bar{\W}}r_{\alpha}\right] +\P\left[ \frac{R_{\alpha}}{1+\bar{\W}}  w\right]  = G
 \\
 & \left(\partial_t + \M_{\underline{b}} \partial_\alpha \right) r +icr-i\P\left[ \frac{g+\underline{a}\, }{1+\W}w\right] = K,
\end{aligned}
\right.
\end{equation}
and  define the associated linear energy
\[
\Elind(w,r) = \int_{\R} (g+\underline{a}) |w|^2 + \Im ( r \bar  r_\alpha)\, d\alpha.
\]
We note that the positivity of the energy is closely related to 
the Taylor sign condition \eqref{taylor}.  Therefore, when $\underline{a}$ is positive, we have
\[
\Elind(w,r) \approx_A \Ez(w,r).
\]

Our first result uses the control parameters $A_{-1/2}$, $A$ and $B$
defined in \eqref{A1-def}, \eqref{A-def}, and\eqref{B-def},
respectively in order to establish (nearly) cubic bounds for the system
\eqref{lin(wr)inhom}:

\begin{proposition}\label{plin-short} The linear equation \eqref{lin(wr)inhom} is well-posed in 
$\dH_0$,  and the following estimate holds:
\begin{equation}\label{lin-gen2}
  \frac{d}{dt} \Elind(w,r)  =  2 \Re \int_{\R} (g+\underline{a}) \bar w \,  G -
i \bar r_\alpha \,  K  + c^2   \Im R |w|^2 \, d\alpha +  O_A(  \uA \uB)  \Elind(w,r).
\end{equation}
\end{proposition}

We will also need a weighted version of this:
\begin{lemma} \label{weighted-en}
Let $f$ be  a real function and  $(w,r)$ solutions to \eqref{lin(wr)inhom}. Then for the 
difference
\[
I: = \frac{d}{dt} \Re \! \int f ((g+ \ua) |w|^2 - i \bar r_\alpha r)\, d\alpha 
- \Re \int   f ((g+ \ua) \bar w F  - i \bar r_\alpha  G)  + (\partial_t + b \partial_\alpha)f  ((g+ \ua) |w|^2 - i \bar r_\alpha r) \, d\alpha 
\]
we have the estimate
\begin{equation}
\begin{split}
|I| \lesssim_A  & \ ( \uB \|f\|_{L^\infty} + \uA \|\vert D\vert^\frac12 f\|_{BMO} ) 
\| (w,r)\|_{L^2 \times \dot H^\frac12}^2.
\end{split}
\end{equation}
\end{lemma}

Our main use for the result in Proposition~\ref{plin-short} is to
apply it to the linearized equation \eqref{lin(wr)}:
\begin{proposition} \label{plin-short+} 
The linearized equation \eqref{lin(wr)} is well-posed in $L^2 \times \dot H^\frac12$,
and the following estimate holds:
\begin{equation}\label{lin2}
\frac{d}{dt} \Elind(w,r)  \lesssim_A (\uB + c \uA)   \Elind(w,r).
\end{equation}
\end{proposition}

\begin{proof}[Proof of Proposition~\ref{plin-short}]
  A direct computation yields
\[
\begin{aligned}
\frac{d}{dt} \int (g+\underline{a})|w|^2 d \alpha = & \ 2 \Re \!\! \int  (g+\underline{a}) \bar w (\partial_t
+ \M_{\underline{b}} \partial_\alpha  ) w +
\underline{a} \bar w   [\underline{b},\P] w_\alpha\,  d\alpha \\ & \ + 
\int \left[ \underline{a}_t+((g+\underline{a}\,)\underline{b})_\alpha \right]|w|^2 \, d \alpha .
\end{aligned}
\]
A similar computation shows that
\[
\frac{d}{dt} \int \Im (  r \partial_\alpha \bar r)  \, d \alpha = 2 \Im  \int    (\partial_t
+ \M_{\underline{b}\, } \partial_\alpha) r  \, \partial_\alpha \bar r \, d \alpha. \ \
\]
Adding the two and using the equations \eqref{lin(wr)inhom}, the quadratic
$\Re (w \bar r_\alpha)$ term cancels modulo another commutator term,
and we obtain
\begin{equation}\label{dE2}
\frac{d}{dt}  \Elind(w,r) = 2\Re \int (g+\underline{a}) \bar w\,   G- i \bar r_\alpha \,  K\, d\alpha 
+ c^2  \int \Im R |w|^2 \, d\alpha
+     \underline{err}_1,
\end{equation}
where
\[
\begin{aligned}
  \underline{err}_1 = & \int \left[ \underline{a}_t+((g+\underline{a}\, )\underline{b}\, )_\alpha- c^2 \Im  R \right] |w|^2 \, d\alpha -
  2\Re \int (g+\underline{a}\, )\frac{R_{\alpha}}{1+\bar{\W}}\vert w\vert ^2\,  d\alpha
  \\ & -2 \Re \int \underline{a}\, \bar w \, (\left[\bar Y, \P\right] (r_\alpha+R_\alpha w) \, d\alpha -2\Re\int \underline{a} \bar{w}[\P,\underline{b}\, ] w_\alpha\, d\alpha .
\end{aligned}
\]

Using the auxiliary function $\underline{M}$ in \eqref{M-def},
we  rewrite it as
\[
\begin{aligned}
  \underline{err}_1 =& \ \int \left( \underline{a}_t+\underline{b}\,
    \underline{a}_{\alpha}- c^2 \Im R\right) |w|^2 -(g+\underline{a}\,
  )( i \frac{c}2 (\W -\bar\W) +\underline{M})\, \vert w\vert ^2 \, d\alpha 
\\  & \ - 2 \Re \int
  \underline{a} \bar w \, (\left[\bar Y,\P \right] (r_\alpha+R_\alpha
  w)
  +   [\P,\underline{b}\, ] w_\alpha) \, d\alpha .\\
\end{aligned}
\]

At this point we need to re-express the error in terms of $a$, $a_1$,
$b$, $b_1$, $M$, and $M_1$. The terms which do not have any $c$ factors
were already estimated in \cite{HIT}, so that we only need to worry about the remaining
components. There is only one exception: the first term in
$\underline{err}_1$, whose counterpart was estimated separately in \cite{HIT},
and  which we also estimate separately here. Thus, the error we want
to bound is composed of $err_1$ (the same as in \cite{HIT}) and
additional terms, which we call $err_{1,1}$:
\[
\begin{aligned}
\underline{err}_1=&err_1+ err_{1,1},
\end{aligned}
\]
where $err_{1,1}$ is further separated into terms
\begin{equation}
\label{err11}
err_{1,1}:=err_{1,1}^1+err_{1,1}^2+err_{1,1}^3-err_{1,1}^4+err_{1,1}^5+err_{1,1}^6,
\end{equation}
which will be estimated separately, and are listed below:
\[
err_{1,1}^1:=\int \left( {a}_t+\underline{b}\, {a}_{\alpha}\right)  |w|^2\, d\alpha ,
\]
\[
err_{1,1}^2:=\frac{c}2 \int \left( {a}_{1,t}+\underline{b}\, {a}_{1,\alpha} - i g(\W-\bar \W)  - 2c \Im R
\right)  |w|^2\, d\alpha,
\]
\[
err_{1,1}^3:=\frac{c}{2} \int  \left( i(g+a)M_1-a_1M - i a(\W-\bar \W)\right) |w|^2 \, d\alpha  , 
\]
\[
err_{1,1}^4:=2\Re\int i\frac{c}{2}a\bar{w} \left[ \P, b_1 \right] w_{\alpha} \, d\alpha,
\]
\[
err_{1,1}^5:=2\Re\int  \frac{c}{2}a_1\bar{w}\left\lbrace \left[ \bar{Y}, \P\right](r_{\alpha}+R_{\alpha}w)+\left[ \P, b\right]w_{\alpha}  \right\rbrace \, d\alpha ,
\]
\[
err_{1,1}^6:=\Re\int i\frac{c^2}{4}a_1 \bar{w}  ( 2\left[P,b_1 \right]w_{\alpha}+
 (M_1-\W+\bar \W) w)   \,d\alpha .
\]

 To conclude the proof of
\eqref{lin-gen2}  we want to show that
\begin{equation}\label{err12}
 |\underline{err}_{1}| \lesssim_{A}  (A+cA_{-1/2}) (B +cA)  \Elind(w,r) .
\end{equation}
From \cite{HIT} we have the bound 
\[
|err_{1}| \lesssim_{A}  AB  \Elind(w,r) ,
\]
so it remains to establish \eqref{err12} for each of the $err_{1,1}$
components. The positive constant $c$ has the role of a scaling
parameter, therefore it makes sense to group terms accordingly. 
\bigskip

{\em The bound for $err_{1,1}^1$}.
This follows by Lemma~\ref{l:a}  in
Section~\ref{s:est}. 

\medskip

{\em The bound for $err_{1,1}^2$}.
This follows by Lemma~\ref{l:a1}, which we included in
Section~\ref{s:est}. 

\medskip

{\em The bound for $err_{1,1}^3$}. Here we use the pointwise 
estimates
\[
\| a \|_{L^\infty} \lesssim A^2, \quad \|a_1\|_{L^\infty} \lesssim_A A_{-1/2}, \quad 
\| M\|_{L^\infty} \lesssim_A AB, \quad \|M_1\|_{L^\infty} \lesssim_A A^2.
\]
The first and the third are from \cite{HIT}, while for the 
second and fourth we use Lemma~\ref{l:N} and Proposition~\ref{l:a1}.
 This yields
\[
 \vert err^2_{1,1}\vert \lesssim c\left( A^2+ ABA_{-1/2}+A^3+A^4\right)\Elind(w,r) .
\]
\medskip

{\em The bound for $err_{1,1}^4$}. For the terms in $err_{1,1}^3$ we
use the pointwise bounds of $a$ and $b_{1,\alpha}$ (recall that
$b_{1,\alpha}=\W-\bar{\W} -M_1$), which were obtained in \cite{HIT}
(see Lemma~2.6 in the \emph{Appendix}), and in Lemma~\ref{l:N},
respectively:
 \[
 \vert a\vert\lesssim A^2 , \quad   \|b_{1,\alpha}\|_{L^{\infty}} \lesssim A+ A^2.
 \] 
 To estimate the commutator in $L^2$ we use a Coifman-Meyer bound, see
 e.g.  Lemma~2.1 in \cite{HIT}:
 \[
\|\left[\P, b_1\right] w_{\alpha}\|_{L^2}\lesssim \Vert b_{1,\alpha}\Vert_{L^{\infty}}\Vert w\Vert_{L^2}.
\]
 Combining the results above leads to
 \[
\vert  err_{1,1}^3\vert \lesssim c \left( A^3+A^4\right)\Elind(w,r) .
 \]
 \medskip

{\em The bound for $err_{1,1}^5$}. After using the above pointwise bound on $a_1$
it remains to  estimate the commutators in $L^2$ (as above). 
The only difficulty here is that we need to move half of derivative
from $r_{\alpha}$ onto $Y$. Such  estimates  were already considered  in \cite{HIT}:
\[
\|\left[\bar Y ,\P\right]  r_\alpha\|_{L^2} \lesssim \| |D|^\frac12 Y\|_{BMO}
\| r \|_{\dot H^\frac12}, \qquad    \| [\P,b] w_\alpha\|_{L^2} \lesssim 
\|b_\alpha\|_{BMO} \|w\|_{L^2},
\]
and suffice due to the bounds for $b$, $a_1$ and $Y$
in Lemma~2.5, and Lemma~ 2.7 (in \cite{HIT}), and  Lemma~\ref{l:a1}.  
For the remaining term in $err_{1,1}^4$ we write
$[\bar Y, \P] (R_\alpha w)= [\bar \P, \bar \P [\bar Y R_\alpha]] w]$ and 
use the Coifman-Meyer Lemma, \cite{cm}, (also discussed in \emph{Appendix}~B in \cite{HIT}) to estimate
\[
\|  \bar \P[\bar \P [\bar Y R_\alpha] w]\|_{L^2} \lesssim 
\|w\|_{L^2} \|\bar \P [\bar Y R_\alpha]\|_{BMO} \lesssim 
\|w\|_{L^2} \||D|^\frac12 Y\|_{BMO} \| |D|^\frac12 R\|_{BMO},
\]
where the bilinear bound in the second step follows after a bilinear
Littlewood-Paley decomposition (again, see \cite{HIT}).
Hence,
\[ 
\vert err_{1,1}^4\vert  \lesssim_A c A_{-1/2}AB \,  \Elind(w,r).
\]
\medskip 

{\em The bound for $err_{1,1}^6$.}  Here we use the pointwise  
bounds on $a_1$, $M_1$ and $b_{1,\alpha}$, along with the Coifman-Meyer commutator
estimate.

\medskip
This concludes the proof of \eqref{err12}. 
\end{proof}

\begin{proof}{Proof of Lemma~\ref{weighted-en}}
We rewrite $I$ in the form 
\[
I: = D_1 + D_2 + D_3 + D_4,
\]
where
\[
\begin{aligned}
D_1 = & \ \int   f  (\partial_t + \ub \partial_\alpha) \ua |w|^2 \, d \alpha,\quad D_3 = \  \int f   \left(b_\alpha |w|^2 - 2 \Re ( \bar w \P\left[ \frac{R_\alpha}{1+ \bar \W} w \right]\right) \, d \alpha,
\\
 \quad D_2 = & \ \Re \int f  \bar \P[\ub w_\alpha] \bar w \,  d \alpha,\quad 
D_4 =  \   \int i f \left( \bar r_\alpha \P\left[ \frac{(g+a) w}{1+ \W}\right] - (g+a) \bar w \P\left[\frac{r_\alpha}{1+\W}\right]\right) \, d \alpha .
\end{aligned}
\]
For the first term we  use the pointwise bounds in Lemma~\ref{l:a} and Lemma~\ref{l:a1}.

For the second term we directly use a Coifman-Meyer estimate to bound the middle factor in $L^2$.

For the third term we use the pointwise bound on $b_1$, and then harmlessly replace $b_\alpha$
by $\P\left[ \dfrac{R}{1+\bar \W}\right]_\alpha$. Then it remains to estimate in $L^2$ the difference 
\[
\P\left[ \frac{R}{1+\bar \W}\right]_\alpha w - \P\left[ \frac{R_\alpha}{1+ \bar \W} w \right] .
\]
If $w$ is the low frequency factor in either term, then we only need its cofactor in $BMO$
and we win. Else we can drop the first projection, and we are left with estimating the difference
\[
\P\left[ \frac{R}{1+\bar \W}\right]_\alpha - \frac{R_\alpha}{1+ \bar \W} 
\]
in $L^\infty$, which was done in \cite{HIT}.

Finally, the last term cancels if we drop the projections. Hence we are left with estimating 
in $L^2$ the expression
\[
\bar \P \left[ \frac{r_\alpha}{1+\W}\right] ,
\]
which is done by a Coifman-Meyer estimate.
In the first term we bound the expression 
\[
\bar \P\left[ \frac{a w}{1+ \W}\right] 
\]
in $L^2$ and then move a half-derivative on $f$.
\end{proof}

\begin{proof}[Proof of Proposition~\ref{plin-short+}]
  To estimate the terms involving $ \mathcal{G}$ and $ \mathcal{K}$ we
  separate the quadratic and cubic parts, but more importantly we
  group these expressions keeping track of the scaling parameter $c$.
In our previous paper \cite{HIT} we have already established the bounds
for the components without $c$, namely
\begin{equation*}
\| \P \mathcal{G}^{(2)}(w, r)\|_{L^2} + \| \P \mathcal{K}^{(2)}(w, r)\|_{\dot H^\frac12}
\lesssim B(\|w\|_{L^2} + \|r\|_{\dot H^\frac12}),
\end{equation*}
while  the cubic and higher terms satisfy
\begin{equation*}
\| \P \mathcal{G}^{(3+)}(w, r)\|_{L^2} + \| \P \mathcal{K}^{(3+)}(w, r)\|_{\dot H^\frac12}
\lesssim_A AB(\|w\|_{L^2} + \|r\|_{\dot H^\frac12}).
\end{equation*} 

Hence  it suffices to estimate the terms with the $c$ factor,
and show that the quadratic terms satisfy
\begin{equation}\label{gk2-est}
\| \P \mathcal{G}_1^{(2)}(w, r)\|_{L^2} + \|\P \mathcal{K}_1^{(2)}(w, r)\|_{\dot H^\frac12}
\lesssim A(\|w\|_{L^2} + \|r\|_{\dot H^\frac12}),
\end{equation}
while  the cubic and higher terms satisfy
\begin{equation}\label{gk3-est}
\| \P \mathcal{G}_1^{(3+)}(w, r)\|_{L^2} + \| \P \mathcal{K}_1^{(3+)}(w, r)\|_{\dot H^\frac12}
\lesssim_A A^2(\|w\|_{L^2} + \|r\|_{\dot H^\frac12}).
\end{equation}

In order to obtain the estimates claimed in
\eqref{gk2-est},\eqref{gk3-est} we use the Coifman-Meyer \cite{cm}
type commutator estimates described in the \emph{Appendix}~B,
Lemma~2.1. 

The bounds for all  terms in  $\P\mathcal{G}_1^{(2)}(w,r)$ are immediate, except for the second, where  we use $(B.10)$ with  $s = \frac12$
and $\sigma = \frac12$, to write
\[
 \| \left[ \P,W\right] \bar{w}_{\alpha} \|_{L^2} \lesssim \| \W \|_{L^{\infty}} \|w\|_{L^2}.
\]

For $\P\mathcal{K}_1^{(2)}(w,r)$ we use again $(B.10)$ with $s = \frac12$,
and $\sigma = \frac12$, and conclude that
\[
\| \left[\P, W\right]  \bar r_\alpha\|_{\dot H^\frac12} \lesssim \| |D|^{\frac12}W\|_{L^{\infty}}  \|r\|_{\dot H^\frac12}, \quad \|  \left[ \P,R\right] \bar{w}\|_{\dot H^\frac12} \lesssim \| \vert D\vert ^{\frac{1}{2}}R\|_{L^{\infty}} \|w\|_{L^2}.
\]
 Thus \eqref{gk2-est} follows.

For the cubic and higher parts of $\mathcal{G}_1$ and
$\mathcal{K}_1$ we apply  the same type of commutator
estimates, as well as the $BMO$ bounds in \cite{HIT}, Proposition ~$2.2$.
Precisely, in  $\mathcal{G}_1^{(3+)}$ there are three nontrivial terms to be estimated 
in $L^2$, namely 
\[
\P[ \bar w_\alpha (\bar Y^2 - 2\bar Y) W], \qquad W \P[ \bar w_\alpha ( \bar Y^2 -2\bar Y+1) W], \qquad  \P [W \bar \P[  w_\alpha ( Y^2 - 2Y +1) \bar W]].
\]
For the first two we use $(B.12)$ and $(B.14)$ from \cite{HIT} as follows:
\begin{equation*}
\|  \P [\bar w_\alpha (\bar{Y}^2-2\bar{Y})W]\|_{L^2} \lesssim
\| w\|_{L^{2}} \| \partial_{\alpha}\P\left[ ( \bar{Y}^2-2\bar{Y}) W\right] \|_{BMO}
\lesssim_A \| \W\|_{BMO} \|Y\|_{L^\infty}  \| w\|_{L^2} .
\end{equation*}
Similarly, for the remaining terms we have 
\[
\|  \P [\bar w_\alpha (\bar{Y}^2-2\bar{Y}+1)W]\|_{L^2} \lesssim_A  \| \W\|_{BMO} \| w\|_{L^2}. 
\]

Finally, we estimate the cubic component of $\mathcal{K}_1$; we again
use $(B.12)$ and $(B.15)$ from \cite{HIT}, and obtain
\[
\| |D|^\frac12 \P[\bar r_\alpha \bar Y  W]\|_{L^2}
\lesssim \| r\|_{\dot H^\frac12} \|  \partial_\alpha  \P[ \bar Y W]\|_{BMO}
\lesssim_A \| r\|_{\dot H^\frac12} \|Y\|_{L^\infty} \| \W\|_{L^{\infty}},
\]
 and also
\[
\begin{split}
\| |D|^\frac12 \P[ \bar w (1-\bar Y) \bar R_\alpha  W]\|_{L^2}
 & \ \lesssim \| w(1-Y)\|_{L^2} \| |D|^\frac12  \P[\bar R_\alpha  W]\|_{BMO}
\\
 & \ \lesssim_A   \|w\|_{L^2} \| |D|^\frac12 R\|_{BMO} \|\W\|_{L^{\infty}}.
\end{split}
\]
This concludes the proof of \eqref{gk3-est}, and thus of the proposition.

\end{proof}

\subsection{Nearly cubic estimates for small data.}  Our goal here is
to investigate the possibility of obtaining cubic estimates for the
system \eqref{lin(wr)} in $L^2 \times \dot H^\frac12$. Unlike in
\cite{HIT}, this is no longer possible. Instead, our more limited goal
in this section is to identify a main portion of the linearized
equation for which cubic estimates are valid, precisely up to $c^2$
terms.  This will come in very handy later on in the proof of cubic
estimates for the differentiated equation.

Our model problem this time is the following subset of \eqref{lin(wr)}:
\begin{equation}\label{lin(wr)inhom+}
\left\{
\begin{aligned}
&\left(\partial_t +  \M_{\underline{b}}  \partial_\alpha \right)  w+\P\left[ \frac{1}{1+\bar{\W}}r_{\alpha}\right] +\P\left[ \frac{R_{\alpha}}{1+\bar{\W}}  w\right]  =  -\P \left[ \W \bar r_\alpha
  \right] + \P\left[ R \bar w_\alpha\right]+   G
 \\
 & \left(\partial_t + \M_{\underline{b}} \partial_\alpha \right) r +icr-i\P\left[ \frac{g+\underline{a}\, }{1+\W}w\right] =-\P[R \bar r_\alpha]  + K.
\end{aligned}
\right.
\end{equation}

In our previous work \cite{HIT}, for the case $c=0$, we have identified
a cubic correction to $\Elind(w,r)$ for which cubic estimates
hold for solutions to the linearized equation, namely
\[
\Elint(w,r):=\int _{\mathbb{R}}(g+\underline{a})\vert w\vert ^2 
+\Im(r\bar{r}_{\alpha})+2\Im(\bar{R}wr_{\alpha})-2\Re(\bar{\W}w^2)\, d\alpha.
\]
Our next result asserts that in our case the time derivative of
$\Elint$ will be quartic at the leading order, but will have a cubic
terms with a coefficient of  $c^2$.  

\begin{proposition}\label{plin-short3} 
Assume that $\uA$ is small. Then we have the energy equivalence
\begin{equation}
\Elint(w,r) \approx \Elind(w,r).
\end{equation}
Further, the linear equation \eqref{lin(wr)inhom} is well-posed in 
$\dH_0$,  and the following estimate holds:
\begin{equation}\label{lin-gen2+}
\begin{aligned}
  \frac{d}{dt} \Elint(w,r)  &=  2\Re \int_{\mathbf{R}}\left[ (g+\underline{a})\bar{w}-i\bar{R}_{\alpha}r_{\alpha}-2\bar{\W}w\right] G
  -i\left[ \bar{r}_{\alpha} -(\bar{R}w)_{\alpha}\right] K\, d\alpha  \\
  &\qquad \qquad + c^2   \Im R |w|^2 \ d\alpha +  O_A(  \uA \, \uB)  \Elind(w,r).
\end{aligned}
\end{equation}
\end{proposition}

Applying this to the linearized equation \eqref{lin(wr)} we obtain:
\begin{proposition} \label{plin-short3+} 
The linearized equation \eqref{lin(wr)} is well-posed in $L^2 \times \dot H^\frac12$,
and the following estimate holds:
\begin{equation}\label{lin2+}
\frac{d}{dt} \Elint(w,r)  =  2 \Re \int_{\R} c\left[ g w\P \mathcal{G}^{(2)}_1(w,r) -i\bar{r}_{\alpha}\P \mathcal{K}^{(2)}_1(w,r)\right] -i\frac{c^2}{2}R |w|^2 \, d\alpha 
 + O_A(  \uA \uB )  \Elind(w,r),
\end{equation}
where  $\P \mathcal{G}^{(2)}_1(w,r)$, $\P \mathcal{K}^{(2)}_1(w,r)$  are as in \eqref{quadratic}.
\end{proposition}

\begin{proof}[Proof of Proposition~\ref{plin-short3} ]
The energy equivalence is immediate due to the estimates in \cite{HIT} for the added
cubic terms, and the pointwise bounds on $\ua$ in Lemmas~\ref{l:a}, \ref{l:a1}.
 
To prove the estimate in \eqref{lin-gen2+} we compute the time
derivative of the cubic component of the energy $\Elint(w,r)$,
using the projected equations for $w$ and $r$ and the unprojected
equations for $R$ and $\W$:
\[
\begin{aligned}
\frac{d}{dt} \! \left( \! \Im\! \int   \bar{R} w r_\alpha \, d\alpha -\Re\! \int \bar{\W}w^2d\alpha \!\right)\! =&
 \Im \! \int \!   - i g\bar \W w r_\alpha     -  \bar{R}  r_\alpha r_\alpha+ i g\bar{R} w  w_\alpha
+ \bar R r_\alpha G + \bar R w K_\alpha \, d\alpha \\
 & \hspace{.7cm}+\Re\int \bar{R}_{\alpha}w^2+2\bar{\W} wr_{\alpha}- 2 \bar \W w G\, d\alpha+ \underline{err_2},
 \end{aligned}
\]
where
\begin{equation}
\label{underline err2}
\begin{aligned}
\underline{err_2} =  & \ \Im \int\!\! \left\lbrace   \left( i\left( \frac{g\bar{\W}^2+a}{1+\bar{\W}} \right)-\underline{b}\bar{R}_{\alpha} -i\frac{c}{2}\frac{\bar{R}\bar{\W}+R\bar{\W}+\bar{N} }{1+\bar{\W}}\right) wr_{\alpha}\right. \\
    & \left. \hspace{1.5cm}   -  \bar{R}w\partial_{\alpha}\left( \M_{\underline{b}\,}r_{\alpha}-i\P\left[ \frac{\underline{a}-\W}{1+\W}w\right] + \P[R\bar{r}_\alpha]\right)\right.
  \\
   &\left.   \hspace{1.5cm}  -\bar{R}r_{\alpha} \left( \M_{\underline{b}}w_{\alpha}-\P\left[ \frac{\bar{\W}}{1+\bar{\W}}r_{\alpha}\right]
   +\P\left[ \frac{R_{\alpha}}{1+\bar{\W}}w\right] + \P[\W \bar{r}_\alpha - R\bar{w}_\alpha] \right) \right\rbrace \, d\alpha \\
&+ \Re\int \left\{ w^2\left( \underline{b} \bar \W_{\alpha}+ \frac{\bar{\W}-\W}{1+\W}\bar{R}_{\alpha} - (1+\bar{\W})\underline{\bar{M}} +i\frac{c}{2}\bar{\W}(\bar{\W}-\W) \right) \right.
\\ &  \hspace{1.5cm}  \left. + 2\bar{\W} w \left( \M_{\underline{b}}w_{\alpha}-
\P\left[ \frac{\bar{\W}}{1+\bar{\W}}r_{\alpha}\right] +\P\left[ \frac{R_{\alpha}}{1+\bar{\W}}w \right]
+ \P[\W \bar{r}_\alpha - R\bar{w}_\alpha]\right) \right\} \, d\alpha .
\end{aligned}
\end{equation}
We now separate \eqref{underline err2} in two parts, namely $err_2$, which was already 
estimated in  \cite{HIT},  and in $err_{2,1}$, 
\[ 
\underline{err}_2:=err_2+i\frac{c}{2}err_{2,1},
\]
where 
\begin{equation}
\label{err2,1}
\begin{aligned}
err_{2,1} :=  & \ \Im \int\!\! \left\lbrace   \left( b_1\bar{R}_{\alpha} -\frac{\bar{R}\bar{\W}+R\bar{\W}+\bar{N} }{1+\bar{\W}}\right) wr_{\alpha}+
  \bar{R}w\partial_{\alpha}\P\left[ b_1r_{\alpha}\right]\right. \\
  &\left.\ \ \ \  \qquad -  \bar{R}w\P\left[ \frac{a_1}{1+\W}w\right] +\bar{R}r_{\alpha}\P\left[ b_1w_{\alpha}\right] \right\rbrace  d\alpha \\
&+ \Re\int \left\lbrace  w^2\left[  -b_1 \bar \W_{\alpha} - (1+\bar{\W})\bar{M}_1 +\bar{\W}(\bar{\W}-\W) \right]   - 2\bar{\W} w \P \left[ b_1w_{\alpha}\right] \right\rbrace \, d\alpha.
\end{aligned}
\end{equation}
We still need to estimate this last error.
 
 Adding $\underline{err}_2$ to \eqref{lin-gen2+} (but applied to solutions to \eqref{lin(wr)inhom+}) gives us
 \begin{equation*}
 \begin{aligned}
 \frac{d}{dt}\Elint(w, r)&=\frac{d}{dt}E^{(2)}_{lin}(w,r)+2 \frac{d}{dt}  \left\lbrace \left( \Im \int   \bar{R} w r_\alpha \, d\alpha -\Re\int \bar{\W}w^2d\alpha\right) \right\rbrace \\
 &=2\Re \int_{\mathbf{R}} (g+\underline{a})\bar{w}\left\lbrace G-\P\left[ \W\bar{r}_{\alpha} \right]+\P\left[ R\bar{w}_{\alpha}\right]\right\rbrace-i\bar{r}_{\alpha}\left\lbrace  K-\P\left[ R\bar{r}_{\alpha}\right]   \right\rbrace \,d\alpha \\
 &\hspace*{1.5cm}+c^{2}\Im \int_{\mathbf{R}} R\vert w\vert^2 \, d \alpha +\underline{err_1}\\
 &\hspace*{1.5cm}+ 2\Im \! \int \!   - i g\bar \W w r_\alpha     -  \bar{R}  r_\alpha r_\alpha+ i g\bar{R} w  w_\alpha
+ \bar R r_\alpha G + \bar R w K_\alpha \, d\alpha \\
 &\hspace*{1.5cm} +2\Re\int \bar{R}_{\alpha}w^2+2\bar{\W} wr_{\alpha}- 2 \bar \W w G\, d\alpha+ 2\underline{err}_2.
 \end{aligned}
 \end{equation*}
 Rewriting the above expression leads to
  \begin{equation*}
 \begin{aligned}
 \frac{d}{dt}\Elint(w, r)&=\frac{d}{dt}E^{(2)}_{lin}(w,r)+2 \frac{d}{dt}  \left\lbrace \left( \Im \int   \bar{R} w r_\alpha \, d\alpha -\Re\int \bar{\W}w^2d\alpha\right) \right\rbrace \\
 &=2\Re \int_{\mathbf{R}} (g+\underline{a})\bar{w}G-i\bar{r}_{\alpha}  K \,d\alpha +c^{2}\Im \int_{\mathbf{R}} R\vert w\vert^2 \, d \alpha +\underline{err_1}\\
 &\quad +  2\Re \int_{\mathbf{R}} (g+\underline{a}) \bar{w}\left\lbrace \P\left[ R\bar{w}_{\alpha}\right] -\P\left[ \W\bar{r}_{\alpha} \right]\right\rbrace+i\bar{r}_{\alpha}  \P\left[ R\bar{r}_{\alpha}\right]   \,d\alpha \\
 &\quad + 2\Im \! \int \!   - i g\bar \W w r_\alpha     -  \bar{R}  r_\alpha r_\alpha+ i g\bar{R} w  w_\alpha
+ \bar R r_\alpha G + \bar R w K_\alpha \, d\alpha \\
 &\quad  +2\Re\int \bar{R}_{\alpha}w^2+2\bar{\W} wr_{\alpha}- 2 \bar \W w G\, d\alpha+ 2\underline{err}_2.
 \end{aligned}
 \end{equation*}
 We obtain
 \begin{equation*}
 \begin{aligned}
  \frac{d}{dt}\Elint(w, r)&= 2\Re \int_{\mathbf{R}}\left[ (g+\underline{a})\bar{w}-i\bar{R}_{\alpha}r_{\alpha}-2\bar{\W}w\right] G
  -i\left[ \bar{r}_{\alpha} -(\bar{R}w)_{\alpha}\right] K\, d\alpha  \\
  &\quad +c^2\Im \int_{\mathbf{R}}R\vert w\vert^2\, d\alpha  +2\Re \int_{\mathbf{R}}\underline{a}\bar{w}\left\lbrace \P\left[ R\bar{w}_{\alpha}\right] -\P\left[ \W\bar{r}_{\alpha} \right]\right\rbrace \, d\alpha +\underline{err}_1+ 2\underline{err}_2.
   \end{aligned}
   \end{equation*}
   We introduce 
   \[
   \underline{err}_3:=2\underline{err}_2-2\Re \int_{\mathbf{R}}\underline{a}\bar{w}\left\lbrace \P\left[ \W\bar{r}_{\alpha} \right]-\P\left[ R\bar{w}_{\alpha}\right] \right\rbrace \, d\alpha,
   \]
   which implies
   \begin{equation*}
 \begin{aligned}
  \frac{d}{dt}\Elint(w, r)&= 2\Re \int_{\mathbf{R}}\left[ (g+\underline{a})\bar{w}-i\bar{R}_{\alpha}r_{\alpha}-2\bar{\W}w\right] G
  -i\left[ \bar{r}-(\bar{R}w)\right]_{\alpha} K\, d\alpha  \\
  &\quad +c^2\Im \int_{\mathbf{R}}R\vert w\vert^2\, d\alpha  +\underline{err}_1+ \underline{err}_3.
   \end{aligned}
   \end{equation*}
Given the bound \eqref{err12} for $\underline{err}_1$, the proof of \eqref{lin-gen2+} is concluded if we show that 
\begin{equation}
\label{err3}
\vert \underline{err}_3\vert \lesssim \underline{A}\, \underline{B}E^{(2)}_{lin}(w,r).
\end{equation}
 Due to the pointwise bound for $\ua$, proved in the \emph{Appendix}, and the Coifman-Meyer type 
$L^2$ bound for  $\P\left[ \W\bar{r}_{\alpha} \right]-\P\left[ R\bar{w}_{\alpha}\right]$ from \cite{HIT},
it suffices to estimate $\underline{err}_2$,  which in turn reduces to estimate $err_{2,1}$,
  \[
  \vert err_{2,1}\vert\lesssim \uA \, \uB E^{(2)}_{lin}(w,r).
  \]
 For the remainder of the proof we separately estimate  several types of terms in
  $err_{2,1}$:
 
\medskip 

\textbf{Terms involving $b_1$.} Here, we use the bounds for $b_1$ from in Lemma~\ref{l:N}:
 \[
 \Vert b_{1, \alpha}\Vert_{L^{\infty}}\lesssim A+A^2, \quad \Vert \vert D\vert^{\frac{1}{2}}b\Vert_{L^{\infty}}\lesssim A_{-1/2}+A_{-1/2}A.
 \]
We first collect all the terms that are contained in the first integral of $err_{2,1}$, and include $b_1$:
\[
I_1:= \int_{\mathbf{R}}b_1\bar{R}_{\alpha}wr_{\alpha}+\bar{R}w\partial_{\alpha}\P\left[ b_1r_{\alpha}\right]  +\bar{R}r_{\alpha}\P\left[ b_1w_{\alpha}\right]\, d\alpha .
\]
After integrating by parts, we cancel two of the terms in $I_1$ and obtain 
\[
I_1=I_2+I_3, 
\]
where 
\[
\begin{aligned}
I_2:&= \int_{\mathbf{R}}\bar{R}_{\alpha}w \bar{\P}\left[ b_1r_{\alpha}\right]\, d\alpha , \mbox{ and } I_3:&= \int_{\mathbf{R}}\bar{R}r_{\alpha}\P\left[ b_1w_{\alpha}\right]  -\bar{R}w_{\alpha}\P\left[ b_1r_{\alpha}\right]\, d\alpha .
\end{aligned}
\]
The bounds for $I_2$ follow easily because it has a commutator structure:
\[
\begin{aligned}
I_2:= \int_{\mathbf{R}}\bar{R}_{\alpha}w \bar{\P}\left[ b_1r_{\alpha}\right]\, d\alpha=\int_{\mathbf{R}}\P \left[ \bar{R}_{\alpha}w\right]  \bar{\P}\left[ b_1r_{\alpha}\right]\, d\alpha =\int_{\mathbf{R}}\left[ \P , w\right] \bar{R}_{\alpha}\, \cdot   \left[\bar{\P}, b_1\right]r_{\alpha}\, d\alpha,
\end{aligned}
\]
where we can estimate both factors in $L^2$ using Coifman-Meyer estimates,
\[
\| \left[ \P , w\right] \bar{R}_{\alpha} \|_{L^2} \lesssim \|
w\|_{L^2} \|R_\alpha \|_{BMO}, \qquad \| \left[\bar{\P},
  b_1\right]r_{\alpha}\|_{L^2} \lesssim \| |D|^\frac12 b_{1} \|_{BMO}
\|r\|_{\dot H^\frac12}.
\]
The bound for $I_3$ follows from Lemma 2.9 in \cite{HIT} (see the \textit{Appendix}).

We next collect all the terms that are contained in the second integral appearing in the expression of $err_{2,1}$, and which include $b_1$:
\[
I_4:=\int_{\mathbf{R}} -b_1\bar{\W}_{\alpha}w^2-2\bar{\W}w\P\left[ b_1w_{\alpha} \right]\, d\alpha .
\]
As before, we integrate by parts and rewrite the expression for $I_4$ as:
\[
I_4:=\int_{\mathbf{R}} b_{1, \alpha}\bar{\W}w^2 +2b_1\bar{\W}ww_{\alpha}-2\bar{\W}w\P\left[ b_1w_{\alpha} \right]\, d\alpha =
\int_{\mathbf{R}} b_{1, \alpha}\bar{\W}w^2 +2\bar{\W}w\bar{\P}\left[ b_1w_{\alpha} \right]\, d\alpha.
\]
The first integral on the RHS is easy to bound since we know that that $b_{1,\alpha}$ is in $L^{\infty}$,
\[
\big| \int_{\mathbf{R}} b_{1, \alpha}\bar{\W}w^2\, d\alpha\big| \lesssim \Vert b_{1,\alpha}\Vert_{L^{\infty}}\Vert \W\Vert_{L^{\infty}}\Vert w\Vert^2_{L^2}.
\]
For the last integral in $I_4$ we use the Coifman-Meyer type estimate (established first in \cite{HIT}) to obtain
\[
\big| \int_{\mathbf{R}}2\bar{\W}w\bar{\P}\left[ b_1w_{\alpha} \right]\, d\alpha\big| \lesssim \Vert \W\Vert_{L^{\infty}}\Vert  w\Vert _{L^2}\Vert  \left[\bar{\P} ,b_1\right] w_{\alpha}  \Vert_{L^2}\lesssim \Vert b_{1,\alpha}\Vert_{L^{\infty}}\Vert \W\Vert_{L^{\infty}}\Vert w\Vert^2_{L^2}.
\]
Thus,
\[
\big| I_4\big|\lesssim \Vert b_{1,\alpha}\Vert_{L^{\infty}}\Vert \W\Vert_{L^{\infty}}\Vert w\Vert^2_{L^2}.
\]
\newpage 

\textbf{Quadrilinear terms bounded via  $L^2 \cdot L^2$ paring.}  Some of the he remaining terms in 
$err_{2,1}$ have straightforward bounds:
\[
\begin{aligned}
\big| I_7\big| :&= \big| \int_{\mathbf{R}} (1+\bar{\W})\bar{M}_1w^2\, d\alpha\big| \lesssim_A (1+A)A^2\Vert w\Vert^2_{L^2}, \\
\big| I_8\big| :&=\big| \int_{\mathbf{R}} \bar{\W}(\bar{\W}-\W)w^2\, d\alpha \big| \lesssim_A A^2 \Vert w\Vert_{L^2}^2,  \\
\end{aligned}
\]
but other require a little bit of work. This includes the following expressions:
\[
\begin{aligned}
I_5:&=\int_{\mathbf{R}} \frac{\bar{R}\bar{\W}+R\bar{\W}+\bar{N}}{1+\bar{\W}}wr_{\alpha}\, d\alpha ,\\
I_6:&=\int_{\mathbf{R}}\bar{R}w\partial_{\alpha}\P \left[ \frac{a_1}{1+\W}w\right]\, d\alpha =-\int_{\mathbf{R}}\partial_{\alpha}\bar{\P}\left[ \bar{R}w\right] \P \left[ a_1(1-Y)w\right]\, d\alpha .\\
\end{aligned}
\] 
To obtain the bound for $I_6$ we use the Cauchy-Schwartz inequality and  Lemma ~2.1 from \cite{HIT}
\[
\big| I_6\big| \lesssim \Vert \partial_{\alpha}\bar{\P}\left[ \bar{R}w\right]\Vert_{L^2}\Vert  \P \left[ a_1(1-Y)w\right]\Vert_{L^2}\lesssim \Vert R_{\alpha}\Vert_{BMO}
\Vert a_1\Vert_{L^{\infty}}(1+\Vert Y\Vert_{L^{\infty}})\Vert w \Vert^2_{L^2}.
\]
Finally, the bound for $I_5$ is also a consequence of commutator estimates; to see this  we rewrite it as follows
\[
I_5:=\int_{\mathbf{R}}f wr_{\alpha}\, d\alpha=\int_{\mathbf{R}} \bar{\P}\left[ fr_{\alpha}\right] w\, d\alpha , \mbox{ where }
f:=\frac{\bar{R}\bar{\W}+R\bar{\W}+\bar{N} }{1+\bar{\W}}.
\]
Thus, 
\[
\big| I_5\big| =\big| \int_{\mathbf{R}} \bar{\P}\left[ fr_{\alpha}\right] w\, d\alpha\big| \lesssim \Vert \left[\bar{\P}, f\right] r_{\alpha}\Vert_{L^2}\Vert w\Vert_{L^2}
\lesssim \Vert \vert D\vert^{\frac{1}{2}}f\Vert_{BMO} \Vert r\Vert_{\dot{H}^{\frac{1}{2}}}\Vert w\Vert_{L^2},
\]
where 
\[
\Vert \vert D\vert^{\frac{1}{2}}f\Vert_{BMO}\lesssim\underline{A}\, \underline{B}.
\]
\end{proof}

\begin{proof}[Proof of Proposition~\ref{plin-short3+} ] To prove the bound in \eqref{lin2+} it suffices to apply the estimate in \eqref{lin-gen2+} with
\[
G=\P\underline{\mathcal{G}} , \quad K=\P\underline{\mathcal{K}}.
\]
In fact, we only have to resume our work in finding the new terms introduced by the vorticity assumption, which are the 
 components carrying the $c$ when expanding the RHS of \eqref{lin(wr)}, $\P\mathcal{G}_1$ and $\P \mathcal{K}_1$:
\[
\underline{\mathcal{G}}=\mathcal{G}-i\frac{c}{2}\mathcal{G}_1, \quad \underline{\mathcal{K}}=\mathcal{K}-i\frac{c}{2}\mathcal{K}_1.
\]
The terms we want to single out are the $c$-cubic terms appearing in the cubic part of the energy $E^3(w,r)$, \eqref{lin-gen2}. 

For this we need to recall the quadratic components of $ \P \mathcal{G}_1$ and  $\P \mathcal{K}_1$
\[
\begin{aligned}
 \P \mathcal{G}^{(2)}_1(w,r) = & \ \P \left[\mathbf{\W} \bar{w}\right] + \P\left[W \bar{w}_{\alpha} \right] + \P\left[\bar \W {w} \right] - \P\left[\W {w} \right] ,
\\ 
  \P \mathcal{K}^{(2)}_1(w,r) =  &\ \P \left[ W \bar r_\alpha\right]-\P \left[ R\bar{w}\right] ,
\end{aligned}
\]
and the antiholomorphic components (which will not matter in this computation anyways, since we work with the projected equations)
\[
\begin{aligned}
  \bar \P \mathcal{G}^{(2)}_1(w,r) = & \ -\bar \P \left[\bar \W
    {w}\right] - \bar \P\left[\bar W {w}_{\alpha} \right] +\bar
  \P\left[\bar \W {w} \right],
  \\
  \bar \P \mathcal{K}^{(2)}_1(w,r) = &\ \bar \P \left[ \bar W
    r_\alpha\right]- \bar \P \left[ \bar R {w}\right].
\end{aligned}
\]
 
 Thus, we need to bound the following terms
 \[
 \int_{\mathbf{R}} g w\P \mathcal{G}^{(2)}_1(w,r) \, d\alpha, \quad  \int_{\mathbf{R}} \bar{r}_{\alpha}  \P \mathcal{K}^{(2)}_1(w,r)\, d\alpha , \quad \int_{\mathbf{R}} c^2   \Im R |w|^2 \ d\alpha.
 \]
 These bounds are all obtained using Lemma~2.1 in \cite{HIT}, for example, for the first integral we obtain
 \[
  \big| \int_{\mathbf{R}} g w\left\lbrace   \left[\P ,\mathbf{\W}\right] \bar{w} + \left[\P, W\right] \bar{w}_{\alpha}  - \left[\P, \W\right] w  \right\rbrace  \, d\alpha\big|
  \lesssim  A\Vert w\Vert_{L^2}^2.
 \]
 For the second integral we need to move half of derivative off of $\bar{r}_{\alpha}$
 \[
\big|  \int_{\mathbf{R}} \bar{r}_{\alpha} \left\lbrace  \left[ \P,  W\right] \bar r_\alpha-\left[\P,  R\right] \bar{w}\right\rbrace \, d\alpha\big|
\lesssim B\Vert r\Vert_{\dot{H}^{\frac{1}{2}}}^2+ A\Vert r\Vert_{\dot{H}^{\frac{1}{2}}}\Vert w\Vert_{L^2}.
 \]
 The last one, is trivial
 \[
 \big|  \int_{\mathbf{R}} R |w|^2 \ d\alpha \big|\lesssim A_{-1/2} \Vert w\Vert^2_{L^2}.
 \]

\end{proof}

\section{ Higher order energy estimates}
\label{s:ee}
The main goal of this section is to prove energy bounds for the
differentiated equations. These are the main ingredient for the
lifespan part of Theorem~\ref{baiatul}, as well as for the cubic result in
Theorem~\ref{t:cubic}. Precisely, we will establish two types energy
bounds for $(\W,R)$ and their higher derivatives. The first one is a
quadratic bound which applies independently of the size of the initial
data; this yields the last part of Theorem~\ref{baiatul}. The cubic
energy bound applies in the small data case, and yields the cubic
lifespan bound in our small data result in Theorem~\ref{t:cubic}.
The large data result is as follows:

\begin{proposition}\label{t:en=large}
 For any $n \geq 0$ there exists an energy functional $\End$ with  the
following properties whenever the conditions \eqref{nocusp} and \eqref{taylor}
hold uniformly:

(i) Norm equivalence:
\begin{equation*}
\End(\W,R) \approx_{\uA} \Ez(\partial^{n} \W, \partial^{n} R). 
\end{equation*}

(ii) Quadratic energy estimates for solutions to \eqref{ww2d-diff}:
\begin{equation}\label{en-quad}
\frac{d}{dt} \End(\W,R)  \lesssim_{\uA} \uB \Ez(\partial^{n} \W, \partial^{n} R).
\end{equation}
\end{proposition}

The small data result is as follows:

\begin{proposition}\label{t:en=small}
a) For any $n \geq 0$ there exists an energy functional $\Ent$ which
  has the following properties as long as  $\uA \ll 1$:

(i) Norm equivalence:
\begin{equation*}
\Ent (\W,R)= (1+ O(\uA)) \Ez (\partial^{n} \W, \partial^{n} R) +  O(c^4\uA)
 \Ez (\partial^{n-1} \W, \partial^{n-1} R),
\end{equation*}

(ii) Cubic energy estimates:
\begin{equation}\label{en-cubic}
\frac{d}{dt} \Ent (\W,R)  \lesssim_{\uA} \uB\, \uA \left( \Ez (\partial^{n} \W, \partial^{n} R) +  c^4
 \Ez (\partial^{n-1} \W, \partial^{n-1} R)\right).
\end{equation}

Here if $n=0$ then $ \Ez (\partial^{-1} \W, \partial^{-1} R)$ is naturally replaced by 
$\E(W,Q)$.
\end{proposition}

The case $n = 0$ of Proposition~\ref{t:en=large} corresponds to
$\dH_0$ bounds for $(\W_{\alpha},R)$. But these functions solve the linearized
equation, so the desired bounds  are  a consequence of Proposition~\ref{plin-short}.
Hence, we will begin with the proof of the cubic bounds for the $n=0$ 
case.

Next, we compute the differentiated equations, first for $n=1$ and then for 
$n \geq 2$, and show that, up to a certain class of bounded errors,
suitable modifications of $(\W^{(n)},R^{(n)})$ solve a linear 
system which is quite similar to the linearized equation. There are two reasons
why we separate the case $n=1$. First, this corresponds to $\dH_0$ bounds for
$(\W_{\alpha},R_\alpha)$, which play a special role as it is our
threshold for local well-posedness. Secondly, there is a subtle difference
in the choice of the modifications alluded above, due to the fact that 
certain terms which are different for  $n\geq 2$ coincide at $n= 1$.

At this point, the  bounds for the linearized equation already yield 
the large data result. What remains is to establish the cubic small data result,
which is where we implement our \emph{quasilinear modified modified energy method}, using the normal 
form calculated in Section~\ref{s:nf}.

As part of the argument we need to express various nonlinear
expressions in terms of their homogeneous expansions. To describe the
decomposition of a nonlinear analytic expression $F$ into homogeneous
components we will use the notation $\Lambda^k F$ to denote the
component of $F$ of homogeneity $k$. We similarly introduce the
operators $\Lambda^{\leq k}$ and $\Lambda^{\geq k}$. We carefully note 
that all of our multilinear expansions are with respect to the diagonal variables 
$W$ and $R$, and not with respect to $(W,Q)$. 

 Compared to  \cite{HIT} here we lack scaling, but we can still introduce 
a notion of order of a multilinear expression. We begin with single terms, 
for which we assign orders as follows:
\begin{itemize}
\item The order of $W^{(k)}$ is $k-1$.
\item The order of $R^{(k)}$ is $k - \frac12$. 
\item The order of $c$ is $\frac12$.
\end{itemize}
For a multilinear form involving products of such terms we 
define the total order as the sum of the orders of all factors.

While not all expressions arising in the $(\W^{(n)},R^{(n)})$ are
multilinear in $(W,R)$, they can be still viewed as multilinear in
$(W,R)$ and undifferentiated $Y$. Since $Y$ scales like $\W$, it is
natural to assign to it the homogeneity zero.  According to this
definition, all terms in the $\W^{(n)}$ equation have order
$n+\frac12$, and all terms in the $R$ equation have order $n$.
Moving on to integral multilinear forms, all $n$-th energies have order  
$2n$, and their time derivatives have order $2n+\frac12$. 

A second useful bookkeeping device will be needed when we deal with
integral multilinear forms. There it makes a difference how
derivatives and also complex conjugations are distributed among
factors. To account for this we define the leading order of a
multilinear form to be the largest sum of the orders of two factors
with opposite conjugations.  Since we only allow nonnegative orders,
for the $n$-th order energy this is at most $2n$. According to our
definition, all the terms in our $n$-th order energy will have order
$2n$, and all terms in its time derivative will have order
$2n+\frac12$. We remark that the half-integers in the definition of
the orders impose a parity constraint in the terms associated to each
power of $c$.

If all factors in a multilinear form have nonnegative orders, this imposes
a constraint on the order of each factor. Unfortunately this does not 
appear to be the case here, as our multilinear forms will also contain factors of $W$ and $R$,
which have negative order. This  is quite inconvenient. Fortunately, there is a simple way
 to avoid negative orders altogether. Precisely, we will never consider such factors
alone, but in combination with $c$;  thus, the allowed factors 
will be $cR$ and $c^2 W$, both of which have order $0$.
We carefully note last remark applies only partially in the special case $n=0$.

\subsection{The case \texorpdfstring{$n=0$}{0}  }

The goal of this subsection is to obtain  cubic energy estimates in $\dH_0$ for the system for 
 diagonal variables $(\W, R)$.  For convenience we recall the system here:
\begin{equation*}
\left\{
\begin{aligned}
 & \W_{ t} + \underline{b} \W_{ \alpha} + \frac{(1+\W) R_\alpha}{1+\bar \W}  =  
\underline{\mathcal{G}}_0
\\
&R_t  + \underline{b} R_{\alpha} +icR- i\frac{g\W-a}{1+\W} =\underline{\mathcal{K}}_0,
\end{aligned}
\right.
\end{equation*}
where 
\begin{equation*}
 \underline{\mathcal{G}}_0  =  (1+\W)\underline{M}+i\frac{c}{2}\W(\W-\bar{\W}), \quad  
\underline{\mathcal{K}}_0 =i\frac{c}{2}\frac{R\W + \bar R \W +N}{1+\W}.
\end{equation*}
We want to be able to apply the quadratic and cubic bounds  for the
``modified'' model \eqref{lin(wr)inhom+}, respectively. For that, we
rewrite the above systems as follows
\begin{equation}
\label{underlineG0}
\left\{
\begin{aligned}
 & \W_{ t} + \underline{b} \W_{ \alpha} + \frac{ R_\alpha}{1+\bar \W} + \frac{ R_\alpha}{1+\bar \W} \W= -\P\left[ \bar{R}_{\alpha}\W\right]+\P \left[  R\bar{\W}_{\alpha}\right] +  \underline{G}_0
\\
&R_t  + \underline{b} R_{\alpha} +icR- i\frac{(g+\underline{a})\W}{1+\W} = -\P\left[R\bar{R}_{\alpha} \right]+\underline{K}_0,
\end{aligned}
\right.
\end{equation}
where, the expressions for $\underline{G}_0$ and $\underline{K}_0$, for the purpose of this section, are
\begin{equation*}
\left\{
\begin{aligned}
 &\underline{G}_0  :=  \W \left( \P \left[ R\bar{\W}_{\alpha}\right]-\P\left[ \bar{R}_{\alpha}\W\right]   \right)+i\frac{c}{2}(1+\W)M_1 +i\frac{c}{2}\W(\W-\bar{\W})\\
 &\qquad +(1+\W)\left\lbrace  \bar{\P} \left[ \bar{R}Y_{\alpha}-R_{\alpha}\bar{Y} \right]-\P\left[ R\partial_{\alpha}(\bar{\W}\bar{Y}) \right] +\P\left[ \bar{R}_{\alpha}\W    Y\right]  \right\rbrace \\
&\underline{K}_0:= \bar{\P}\left[\bar{R}R_{\alpha} \right] -  i\frac{c}{2}N.
 \end{aligned}
\right.
\end{equation*}
In \eqref{underlineG0} we have identified the leading part of the  equation.  We want to
interpret the terms $(\uG_0,\uK_0)$ on the right as mostly perturbative, but also pay
attention to the holomorphic quadratic part, given by 
\begin{equation}
\label{quad-pgk0}
\begin{split}
\P  G_0^{(2)} = & \   i\frac{c}{2} \left( \P\left[ W\bar{Y}\right]_{\alpha}  +
\W^2-\P\left[ \W\bar{\W} \right]\right) , \quad \P  K_0^{(2)} =  \  i\frac{c}{2} \P \left[ W\bar{R}_{\alpha}-\bar{\W}R \right] .
\end{split}
\end{equation}
 Precisely, we first claim that the quadratic and cubic parts  
of $(\uG_0,\uK_0)$ satisfy the bounds 
\begin{equation}\label{quad-gk0}
\| (\uG_0^{(2)},\uK_0^{(2)} )\|_{\dH_0} \lesssim_{A} \uB \norm_0, 
\end{equation}
respectively 
\begin{equation}\label{cube-gk0}
\| (\uG_0^{(3)},\uK_0^{(3)} )\|_{\dH_0} \lesssim_{A} A \uB \norm_0,
\end{equation}
where
\begin{equation*}
\norm_0 = \| (\W, R)\|_{L^2 \times \dot H^\frac12}.
\end{equation*}
The bounds for the components of $M_1$ and $N$ are discussed in
Lemma~\ref{l:N}. The $\W$ prefactors in $\underline{G}_0$ are harmless, as they
are bounded by $A$ in $L^\infty$. For the remaining terms it suffices
to use Coifman-Meyer type estimates discussed in \emph{Appendix}~\ref{s:est}. For instance we have
 \[
\|  \P \left[ R\bar{\W}_{\alpha}\right]\Vert_{L^2}\lesssim  \Vert \left[  \P,  R\right] \bar{\W}_{\alpha}\Vert_{L^2}
\lesssim \Vert \vert D\vert^{\frac{1}{2}} \W\Vert_{BMO}\Vert R\Vert_{\dot{H}^{\frac{1}{2}}}\lesssim_A B\norm_0,
\] 
and all other terms in $\uG_0$ are similar. Finally, for the first term in $\uK_0$ we have 
\[
\| \P \left[\bar{R}R_{\alpha} \right] \|_{\dot H^\frac12} \lesssim \| |D|^\frac12 R\|_{L^2} \|R_\alpha\|_{BMO}
\lesssim B \norm_0.
\]

We are now ready to look at the cubic energies.  We start by
constructing the cubic normal form energy by selectioning the quadratic
and cubic terms from the corresponding linear energy for the normal
form variables. Precisely, we have
 \[
 \mathcal{E}_0(\tW_{\alpha}, \tQ_{\alpha})=\mathcal{E}_0(\W,
 Q_{\alpha})+2\int_{\mathbb{R}}\Re \bar{\W}\partial_{\alpha}
 W^{[2]}-2\Im \bar{Q}_{\alpha\alpha}\partial_{\alpha}Q^{[2]}\, d\alpha
 + \text{quartic}.
 \]
 In the first term we substitute $Q_{\alpha}=R(1+\W)$. In the integral
 we use the expressions \eqref{nf-eq}, integrate by parts to eliminate
 the $\partial^{-1}W$ and $Q$ factors, and then substitute
 $Q_{\alpha}$ by $R$. Separating the outcome of this computations into
 a leading part and a lower order part, we write it in the form
\[
E_{NF}^{(3)}:=E_{NF, high}^{(3)}+E_{NF, low}^{(3)},
\]
where 
\[
E^3_{NF, high}(\W,R):=\int _{\mathbb{R}}(g+c\Re R)\vert \W \vert ^2  +\Im(R\bar{R}_{\alpha})+2\Im(\bar{R}\W R_{\alpha})-2\Re(\bar{\W}\W^2)\, d\alpha ,
\]
and 
\[
\begin{aligned}
E_{NF, low}^{(3)}(\W, R):=& -c  \int_{\mathbb{R}} 2\Re R \left\lbrace \vert \W\vert^2 -\Im ( \bar{R}_{\alpha}R) \right\rbrace -\W\bar{W}R_{\alpha}+\W^2\bar{R}\, d\alpha \\
                   &-\frac{c^2}{g}\int_{\mathbb{R}}\frac{5}{2} \Im W\left\lbrace \vert \W\vert^2 -\Im (\bar{R}_{\alpha} R)\right\rbrace -\frac{1}{2}\bar{W}\W^2- \frac{1}{2}\W\bar{R}^2\, d\alpha\\
                   &-3\frac{c^3}{2g}\Im \int_{\mathbb{R}} R W\bar{\W}\, d\alpha -3\frac{c^4}{2g}\Re \int_{\mathbb{R}}\W \vert W\vert^2\, d\alpha.
\end{aligned}
\]
For the leading order part $E^{(3)}_{NF, high}(\W,R)$ we consider the appropriate quasilinear correction 
\[
E^{(3)}_{high}(\W,R):=\int _{\mathbb{R}}(g+\underline{a})\vert \W \vert ^2  +\Im(R\bar{R}_{\alpha})+2\Im(\bar{R}\W R_{\alpha})-2\Re(\bar{\W}\W^2)\, d\alpha,
\]
and the remainder $E_{NF, low}^3(\W, R)$ remains unchanged. Hence, we define the quasilinear cubic energy
\[
E^{0,(3)}:=E_{high}^3+E_{NF, low}^3.
\]
It remains to show that this energy has all the right properties in Proposition~\ref{t:en=small}.

 We begin with the energy equivalence. For the leading part this had already been in the context of the linearized equation (see Proposition \ref{plin-short3}), 
\[
E^{(3)}_{high}(\W, R)\approx(1+O(A))E^{(2)}_{lin}(\W, R).
\]
So it remains to show that
\[
E^{(3)}_{NF, low}(\W, R)\approx O(\uA)\left( \mathcal E_0(\W, R)+c^4\mathcal E(W, Q)\right) .
\]
The bound is straightforward for all terms not containing $R_{\alpha}$. So we now consider those.
For the first one we have
\[
\left| \int_{\mathbb{R}}R_{\alpha}\vert R\vert^2\, d\alpha \right| \lesssim \Vert R\Vert_{\dot{H}^{\frac{1}{2}}} \Vert \vert D\vert ^{\frac{1}{2}}(R\Re R)\Vert_{L^2}
\lesssim  \Vert R\Vert_{\dot{H}^{\frac{1}{2}}} ^2\Vert R\Vert_{L^{\infty}} \lesssim A_{-1/2}\norm_0^2,
\]
which suffices. For the next term we have
\[
\big| \int_{\mathbb{R}} R_{\alpha} \bar{W}\W\big| \lesssim \Vert R\Vert_{\dot{H}^{\frac{1}{2}}} \Vert \vert D\vert^{\frac{1}{2}}\P \left[ W\bar{\W}\right]\Vert_{L^2}
\lesssim  \Vert R\Vert_{\dot{H}^{\frac{1}{2}}}  \Vert \vert D\vert^{\frac{1}{2}}W\Vert_{BMO}\Vert \W\Vert_{L^2}\lesssim A_{-1/2}\norm_0^2.
\]
Finally, 
\[
\begin{aligned}
\Vert R_{\alpha}\bar{R} W \Vert _{L^2}\lesssim \Vert R\Vert_{\dot{H}^{\frac{1}{2}}} \Vert \vert D\vert^{\frac{1}{2}} (\bar{R}W)\Vert_{L^2}
\lesssim (\Vert \vert D\vert^{\frac{1}{2}}R \Vert_{BMO} \Vert W\Vert_{L^2} +\Vert R\Vert_{L^{\infty}} \Vert \vert D\vert^{\frac{1}{2}}W\Vert_{L^2}) \Vert R\Vert_{\dot{H}^{\frac{1}{2}}} .
\end{aligned}
\]

Now we consider the time derivative of modified quasilinear energy $E^{0,(3)}$ in order to prove the bound \eqref{en-cubic}. By construction we know that 
\[
\Lambda^{\leq 3}\frac{d}{dt}E^{0,(3)}=0.
\]
We note that this is an algebraic property, which follows from the normal form based construction
even though the normal form itself is unbounded.
Therefore it remains to estimate 
\[
\Lambda^{\geq 4}\frac{d}{dt}E^3=\Lambda^{\geq 4}\frac{d}{dt}E_{NF, high}^3+\Lambda^{\geq 4}\frac{d}{dt}E_{NF, low}^3.
\]
Due to \eqref{quad-gk0} and \eqref{cube-gk0} the estimate for the first
term on the RHS follows directly from the bound \eqref{lin-gen2+} in
Proposition~\eqref{plin-short3} for the leading part of the
linearized equation. Hence, it remains to consider the last term.

The $E^3_{NF, low}(\W, R)$ is a trilinear form of order zero and leading order  zero. We compute its time time derivative using the relation 
\begin{equation} \label{dt-tri}
\frac{d}{dt} \int f_1 f_2 f_3 d \alpha = \int (\partial_t + \ub \partial_\alpha)f_1 f_2 f_3
+  f_1 (\partial_t + \ub \partial_\alpha) f_2 f_3 +  f_1 f_2 (\partial_t + \ub \partial_\alpha) f_3
- \ub_\alpha f_1 f_2 f_3\, d \alpha .
\end{equation}
Then its time derivative will be a multilinear form of order
$\frac{1}{2}$, and also of leading order $\frac{1}{2}$. By inspection,
we see that in this time derivative we can associate each $W$ with a
$c^2$ factor and each $R$ with a $c$ factor, so that all each of the
factors in all of the multilinear monomials have degree at least zero.
Then, each multilinear monomial in $\Lambda^{\geq 4}\left(\dfrac{d}{dt}E_{NF,
  low}^3\right)$ contains exactly one factor of order $\frac{1}{2}$, and the
rest are all factors of order zero.  The factor of order $\frac{1}{2}$
can be either $R_{\alpha}$ or $c$, and the factors of order zero could
be $\W$, $Y$, $c^2W$, or $cR$. We have two cases to consider:

a) If the order $\frac{1}{2}$ factor is $c$, then we simply bound the two of the remaining factors in $L^2$ and the others in $L^{\infty}$.

b) If the order $\frac{1}{2}$ factor is $R_{\alpha}$, we use the
$\dot{H}^{\frac{1}{2}}$ bound for $R$, and we are left with estimating
an expression of the form $\vert D\vert^{\frac{1}{2}} (f_1\cdots
f_k)$, where $k\geq 3$, and each $f_k$ has order zero. For this we use
Lemma~\ref{l:com} to estimate
\[
\Vert \vert D\vert^{\frac{1}{2}} (f_1\cdots f_k)\Vert_{L^2}\lesssim \sum_{j}\Vert \vert D\vert^{\frac{1}{2}} f_j\Vert_{L^{\infty}} \Vert \prod_{l\neq j} f_{l}\Vert_{L^2},
\]
and conclude as above.

\subsection{The differentiated equations for 
\texorpdfstring{$n=1$}.}
We begin by differentiating \eqref{ww2d-diff}, in order to obtain a system for $(\W_{\alpha}, R_{\alpha})$:
\[
\left\{
\begin{aligned}
 & \! \W_{ \alpha t} \! + \! \underline{b} \W_{ \alpha \alpha}  \! + \! \frac{[(1\!+\!\W) R_\alpha]_\alpha}{1+\bar \W}   = \!
(\underline M  \! - \! \underline{b}\, _\alpha) \W_\alpha \!+  (1+\W) (R_\alpha \bar Y_\alpha\! + \! \underline{M}\, _\alpha) 
  +i\frac{c}{2}\left[\W (\W-\bar{\W}) \right] _{\alpha}
 \\
 & \! R_{\alpha t} \!+ \underline{b}\, R_{\alpha\alpha}  +icR_{\alpha}-i \left( \frac{(g+a)\W_{\alpha}}{(1+\W)^2} -\frac{a_{\alpha}}{1+\W}\right) =- \underline{b}\, _\alpha R_\alpha+  
 i\frac{c}{2}\left(\frac{\W(R+ \bar{R})+N}{1+\W}\right)_{\alpha}.
\end{aligned}
\right.
\]
We can expand the last term in the second equation, putting together
all terms which involve $\W_\alpha$. The reason for this is that the $\W_\alpha$ terms 
will be unbounded in a suitable sense, and need to be treated as part of the leading 
order linear operator. This will lead us to the coefficient $a_1$, which will be moved to the left, 
completing the coupling coefficient $a$ to $\underline{a}$. We have 
\[
\begin{split}
\left(\frac{\W(R+ \bar{R})+N}{1+\W}\right)_{\alpha} = & \ \frac{\W_\alpha}{(1+\W)^2}{\left\{ (R +\bar R)   -N  \right\}  } + \frac{\W   (R_\alpha +\bar R_\alpha) +N_\alpha    }{1+\W}
\\
= & \ \frac{{a}_1 \W_\alpha }{(1+\W)^2}   + \frac{\W   (R_\alpha +\bar R_\alpha) +N_\alpha    }{1+\W}.
\end{split}
\]  
Then our system becomes
\[
\left\{
\begin{aligned}
  & \! \W_{ \alpha t} \! + \! \underline{b} \W_{ \alpha \alpha}  \! + \! \frac{[(1\!+\!\W) R_\alpha]_\alpha}{1+\bar \W}   = \!
(\underline M  \! - \! \underline{b}\, _\alpha) \W_\alpha \!+  (1+\W) (R_\alpha \bar Y_\alpha\! + \! \underline{M}\, _\alpha) 
  +i\frac{c}{2}\left[\W (\W-\bar{\W}) \right] _{\alpha}
 \\
 & R_{t\alpha} + \underline{b}\, R_{\alpha\alpha}  +icR_{\alpha}-i \frac{(g+\underline{a})\W_{\alpha}}{(1+\W)^2}  =- \underline{b}\, _\alpha R_\alpha -i  \frac{a_{\alpha}}{1+\W}
+ i \frac{c}2 \frac{\W   (R_\alpha +\bar R_\alpha) +N_\alpha    }{1+\W}.
\end{aligned}
\right.
\]

In order to better compare this with the linearized system we introduce the
modified variable $\bR := R_\alpha(1+\W)$ to further obtain
\[
\left\{
\begin{aligned}
  & \W_{ \alpha t} + \underline{b}\, \W_{ \alpha \alpha} +
  \frac{\bR_\alpha}{1+\bar \W} =
  (\underline M - \underline{b}\, _\alpha) \W_\alpha  +   \bR \bar Y_\alpha +
(1+\W) \underline{M}\, _\alpha +i\frac{c}{2}\left[\W (\W-\bar{\W}) \right]
  _{\alpha}
  \\
  & \bR_{t} + \underline{b}\, \bR_{\alpha} +ic\bR  -  i\frac{(g+\underline{a})\W_\alpha}{1+\W}
= \left(\underline M - \underline{b}\,
    _\alpha - \frac{R_\alpha}{1+\bar \W}+ i\frac{c}{2}  (\W-\bar{\W})\right) \bR 
    -  i a_\alpha  \\
  &\hspace*{6.5cm}+i\frac{c}{2}\left[ \bar \W R_\alpha + \W \bar R_\alpha + N_\alpha\right].
\end{aligned}
\right.
\]

Expanding the $\underline{b}\, _\alpha$  terms via \eqref{M-def1} this yields
\begin{equation}
\left\{
\begin{aligned}
&  \W_{ \alpha t} + \underline{b}\,  \W_{ \alpha \alpha}  + \frac{\bR_\alpha}{1+\bar \W} + \frac{R_\alpha}{1+\bar \W} \W_\alpha  = \underline{G}\,_1
\\
&  \bR_{t} + \underline{b}\, \bR_{\alpha}+ic\bR -   i\frac{(g+a)\W_\alpha}{1+\W}  =\underline{K}\, _1,
\end{aligned}
\right.
\label{WR-diff}
\end{equation}
where on the right we have placed all terms which should be thought of as ``bounded''.
Precisely, $( \underline{G}\,_1,\underline{K}\, _1)$ have the form
\[
\underline{G}\, _1:=G_1-i\frac{c}2 G_{1,1}, \qquad  \underline{K}\, _1:=K_1-i\frac{c}2 K _{1,1},
\]
where the leading components have the same form as in \cite{HIT},
\begin{equation*}
\left\{
\begin{aligned}
 G_{1}  = & \ \bR \bar{ Y}_{\alpha} - \frac{\bar{R}_{\alpha}}{1+\W} \W_\alpha + 2 M   \W_{\alpha}  +  (1+\W)  {M}\, _{\alpha} 
\\
 {K}_{1} =& \ -2\left(   \frac{\bar{ R}_{\alpha} }{1+\W }+ \frac{R_{\alpha} }{1+\bar{\W}} \right)  \bR  + 2 {M}\, \bR + 
( R_{\alpha} \bar{R} _{\alpha} -i a_{\alpha} ),
\end{aligned}
\right.
\end{equation*}
while the extra terms containing the vorticity $c$ are as follows:
\[
\left\{
\begin{aligned}
&  G _{1,1}  =  2M_1   \W_{\alpha}  +  (1+\W)  M _{1,\alpha} -2\W_{\alpha}(\W-\bar{\W}) - \W (\W_{\alpha} -\bar{\W}_{\alpha}),
\\
& K _{1,1}=  2M_1\, \bR -2\bR (\W -\bar{\W} )-  \P [ 2 \W \bar R_\alpha\! + W \bar R_{\alpha \alpha} \! -\! \bar \W_\alpha R]-
 \bar \P [ 2 \bar \W  R_\alpha + \bar W  R_{\alpha \alpha} \! - \! \W_\alpha \bar R] .
\end{aligned} \right.
\]

Here on the left we have again the leading part of the linearized
equation.  Following the same approach as in \cite{HIT}, we will
interpret the terms on the right as mostly perturbative, but also keep track
of their quadratic part.  Thus, for bookkeeping purposes, we
introduce two types of error terms, denoted $\errw$ and $\errr$, which
correspond to the two equations. The bounds for these errors are in
terms of the control variables $\uA,\uB$, as well as the $L^2$ type norm
\begin{equation*}
\norm_1 = \| (\W_\alpha,R_\alpha)\|_{L^2 \times \dot H^\frac12}.
\end{equation*}

In \cite{HIT}, by $\errw$ we denote terms $G_1$ which  satisfy the estimates
\begin{equation*}
\|P G_1\|_{L^2} \lesssim_A AB \norm_1,
\end{equation*}
as well as either one one the following two:
\begin{equation*}
\|\bar P G_1\|_{L^2} \lesssim_A B \norm_1 \quad \text{or} \quad \|\bar P G_1\|_{\dot H^{-\frac12}} 
\lesssim_A A \norm_1.
\end{equation*}
Here, in order to manage the new $c$-terms, we also include in $\errw$
the terms in $G_{1,1}$ for which
\begin{equation*}
\|P G_{1,1}\|_{L^2} \lesssim_A A^2 \norm_1,
\end{equation*}
and either of the following two holds:
\begin{equation*}
\|\bar P G_{1,1}\|_{L^2} \lesssim_A A \norm_1 \quad \text{or} \quad 
 \|\bar P G_{1,1}\|_{\dot H^{-\frac12}} \lesssim_A A_{-1/2} \norm_1.
\end{equation*}

Similarly,  by $\errr$ we have denoted in \cite{HIT} the  terms $K_1$  which satisfy the estimates
\begin{equation*}
\|P K_1\|_{\dot H^\frac12} \lesssim_A AB  \norm_1, \qquad \| P K_1\|_{L^2} \lesssim_A A^2  \norm_1,
\end{equation*}
and
\begin{equation*}
 \|\bar P K_1\|_{L^2}  \lesssim_A A  \norm_1.
\end{equation*}
To that,  for the $c$-terms we add expressions $K_{1,1}$ which satisfy the estimates
\begin{equation*}
\|P K_{1,1}\|_{\dot H^\frac12} \lesssim_A A^2  \norm_1, \qquad \| P K_{2,1}\|_{L^2} \lesssim_A AA_{-\frac{1}{2}}  \norm_1,
\end{equation*}
and
\begin{equation*}
\mbox{ either } \|\bar P K_{1,1}\|_{L^2}  \lesssim_A A_{-1/2}  \norm_1 \ \ \  \mbox{  or  }   \ \ \ 
\|\bar P K_{1,1}\|_{\dot{H}^{\frac12}}  \lesssim_A A  \norm_1.
\end{equation*}
We will rely on \cite{HIT} for the analysis of the expressions $G_1$
and $K_1$, and handle just the new entries i.e., the new terms accompanied
by the vorticity factor $c$. 

We recall that the use of the more relaxed quadratic
control on the antiholomorphic terms, as opposed to the cubic control
on the holomorphic terms, is motivated by the fact that the equations
will eventually get projected on the holomorphic space, so the
antiholomorphic components will have less of an impact. A key property
of the space of errors is contained in Lemma~$3.3$ in \cite{HIT},
which we will refer to as \emph{the multiplicative bounds lemma}. For
convenience we recall it here:
\begin{lemma}\label{l:err}
Let $\Phi$ be a function which satisfies
\begin{equation}\label{Phi-est}
\| \Phi\|_{L^\infty} \lesssim A, \qquad \||D|^\frac12 \Phi\|_{BMO} \lesssim B.
\end{equation}
Then, we have the multiplicative bounds
\begin{equation}
 \Phi \cdot  \errw = \errw, \qquad \Phi \cdot \errr = \errr,
\end{equation}
\begin{equation}
 \Phi \cdot  P \errw = A\, \errw, \qquad \Phi \cdot P \errr = A\, \errr.
\end{equation}
\end{lemma}

We now expand some of the terms in the above system. For this we
will use the following bounds for $M_1$ (see Lemma~\ref{l:N}):
\begin{equation}\label{M-bd}
\| M_1\|_{L^\infty} \lesssim A^2, \qquad \|M_1\|_{\dot H^{1}} \lesssim A \norm_1.
\end{equation}
First we note that
\begin{equation}\label{Mwr}
M_1\W_\alpha = \errw, \qquad M_1 \bR = \errr.
\end{equation}
The first is straightforward in view of pointwise bound for $M_1$.
For the second, by the \emph{multiplicative bounds lemma}, we can
replace $M_1\bR$ by $M_1 R_\alpha$.  After a Littlewood-Paley
decomposition, the $\dot H^\frac12$ estimate for $M_1 R_\alpha$ is a
consequence of the pointwise bound in \eqref{M-bd} for low-high and
balanced interactions, and of the $\dot H^\frac12$ bound in
\eqref{M-bd} combined with Lemma~$2.1$ from \cite{HIT} (see
\emph{Appendix}), provided that we also move half of derivative from
$R_{\alpha}$ onto $M_1$ for the high-low interactions, arriving at
\[
\Vert M_{1}R_{\alpha}\Vert_{\dot{H}^{\frac{1}{2}}}\lesssim \Vert M_1\Vert_{L^{\infty}}\Vert R_{\alpha}\Vert_{\dot{H}^{\frac{1}{2}}} +\Vert \vert D\vert ^{\frac{1}{2}}R\Vert_{L^{\infty}}\Vert M_1\Vert_{\dot{H}^1}\lesssim A^2\norm_1 .
\]

Next we consider $(1+\W) M_{\, 1,\alpha} $, for which we claim that 
\begin{equation}\label{Malpha}
\begin{aligned}
&M_{1,\alpha}= \W_\alpha \bar \W + 2 \W \bar \W_\alpha + \P[ W \W_{\alpha \alpha}]
+\errw, \\
&\P \left[  W\bar{\W}_{\alpha \alpha}+ \bar{\W}_{\alpha}\W\right]=A^{-1}\errw.
 \end{aligned}
\end{equation}

By Lemma~\ref{l:err}, this shows that
\[
(1+\W) M_{1,\alpha} = \W_\alpha \bar \W + 2 \W \bar \W_\alpha + \P[ W \W_{\alpha \alpha}] +\errw.
\]
To prove \eqref{Malpha} we discard cubic terms to  rewrite
\[
\begin{aligned}
M_{1,\alpha}&= \W_{\alpha}\bar{Y}+2\W\bar{Y}_{\alpha}+\P\left[W\bar{Y}_{\alpha \alpha} \right]
-\bar{\P}\left[   \W_{\alpha}\bar{Y} + 2\W\bar{Y}_{\alpha} + \bar W_\alpha Y  + 2 \bar{\W}Y_{\alpha}
+ \bar W Y_{\alpha \alpha}   \right]\\ 
&=\W_{\alpha}\bar{\W} + 2\W\bar{\W}_{\alpha}+\P\left[ W\bar{\W}_{\alpha \alpha}\right] 
-\P\left[ f\right]  -\bar{\P}\left[ g \right] + \errw ,
\end{aligned} 
\]
where 
\[
f=W(\bar{\W}\bar{Y})_{\alpha \alpha}, \qquad  g =\W_{\alpha}\bar{Y} + 2\W\bar{Y}_{\alpha} + \bar W_\alpha Y  + 2 \bar{\W}Y_{\alpha}
+ \bar W Y_{\alpha \alpha}  .
\]
Finally,  for $f$ and $g$ we have $L^2$ bounds
\[
\Vert \P f\Vert_{L^2}\lesssim A^2\norm_1, \ \ \ \Vert \bar{\P}g\Vert_{L^2}\lesssim A\norm_1, 
\]
which follow from  a commutator type bound 
\begin{equation}\label{com2w}
\| P[W \bar \Phi_{\alpha \alpha}]\|_{L^2} \lesssim \|W_\alpha\|_{BMO} \|\Phi_\alpha\|_{L^2}
\end{equation}
derived from Lemma~2.1 in \cite{HIT}. 

The last term in $K_{1,1}$ is antiholomorphic, and also easily placed in $\errr$ by Coifman-Meyer type 
bounds.

Taking into account all of the above expansions, it follows that our system can be rewritten in the form
\[
\left\{
\begin{aligned}
&  \W_{ \alpha t} + \underline{b}\,  \W_{ \alpha \alpha}  + \frac{\bR_\alpha}{1+\bar \W} + \frac{R_\alpha}{1+\bar \W} \W_\alpha  =Princ \, (G_1)-i\frac{c}{2}Princ \, (G_{1,1})+\errw  
\\
&  \bR_{t} + \underline{b}\, \bR_{\alpha}+ic\bR -   i\frac{(g+\underline{a}\,)\W_\alpha}{1+\W}  =Princ \, (K_1)-i\frac{c}{2}Princ \, (K_{1,1}) + \errr, 
\end{aligned}
\right.
\]
where $Princ \, (G)$ refers to the terms in $G$ that cannot be treated as error; they are quadratic 
and higher order terms. We list their expressions below
\begin{equation}
\label{PrincG}
\begin{cases}
&Princ \, (G_1):= 2 \bR \bar Y_\alpha - \dfrac{2 \bar R_\alpha\W_\alpha}{1+\W}  + P[ R \bar \W_{\alpha \alpha}- \bar R_{\alpha \alpha}\W] \\
&\\
&Princ \, (G_{1,1}):= -3 \W  \W_{ \alpha} + 3 \bar \W \W_{ \alpha} + 3\W  \bar 
W_{\alpha} + P[W \bar \W_{\alpha\alpha}],
\\
\end{cases}
\end{equation}
respectively
\begin{equation}
\label{PrincK}
\begin{cases}
&Princ \,(K_1):=-2\left(\dfrac{\bar R_\alpha}{1+\W}+ \dfrac{R_\alpha}{1+\bar \W}\right)  \bR  - P[\bar R_{\alpha \alpha} R]   \\            
&\\                          
&Princ \, (K_{1,1}):=-2\bR (\W -\bar{\W}) -  \P [ 2 \W \bar R_\alpha + W \bar R_{\alpha \alpha} - \bar \W_\alpha R].
\\
\end{cases}
\end{equation}

One might wish to compare this system with the linearized system which
was studied before. However, of the terms on the right, i.e., those
in $Princ \, (G_1)$ and $Princ \, (K_1)$  cannot be all bounded in $L^2 \times \dot
H^\frac12$, even after applying the projection operator $\P$.
Precisely, the terms on the right which cannot be bounded directly in
$L^2 \times \dot H^\frac12$ are $ \displaystyle - 2\frac{\bar
  R_\alpha}{1+\W} \W_\alpha$, respectively $\displaystyle -2
\left(\frac{\bar R_\alpha}{1+\W}+ \frac{R_\alpha}{1+\bar \W}\right)
\bR$. There are no such problematic terms in the $ Princ \, (G_{1,1})$
or in $Princ\, (K_{1,1})$.

As in \cite{HIT}, we eliminate these terms  by conjugation with respect to
a real exponential weight $e^{2\phi}$,  where $\phi = - 2\Re \log(1+\W)$.  Then
\[
\phi_\alpha = - 2 \Re \frac{\W_\alpha}{1+\W}, \qquad
(\partial_t + \underline{b} \partial_\alpha) \phi =  2 \Re \frac{R_\alpha}{1+\bar \W} - 2 \underline{M}-ic(\W -\bar{\W})\Re \frac{\W}{1+\W}\, .
 \]
We denote the weighted variables by
\[
\mfw = e^{2\phi} \W_\alpha, \qquad \mfr = e^{2\phi} \bR.
\]

Using \eqref{Mwr} and Lemma~\ref{l:err} it follows that
 $\underline{M}w = \errw$, $\underline{M}r = \errr$. Similarly 
we have 
\[
c \W_{\alpha}(\W-\bar{\W})\Re Y = \errw, \qquad  c \bR(\W -\bar{\W})\Re Y = \errr .
\]
Then the only significant effect of the conjugation is to eliminate the above problem terms,
and we are left with 
 \[
\left\{
\begin{aligned}
&  \mfw_{ t} + \underline{b} \mfw_{\alpha}  + \frac{\mfr_\alpha}{1+\bar \W}  + \frac{R_\alpha}{1+\bar \W} \mfw =
   \P [ R \bar \W_{\alpha \alpha}- \bar R_{\alpha \alpha} \W]-i\frac{c}{2}Princ \, (G_{1,1}) +  \errw \\
\\
&  \mfr_{t} + \underline{b}\mfr_{\alpha}  +ic\mfr- i \frac{(g+\underline{a}\,) \mfw}{1+\W} =  -  \P [\bar R_{\alpha \alpha} R]-i\frac{c}{2} Princ \, (K_{1,1})   + \errr .\\
\end{aligned}
\right.
\]
Finally, as in \cite{HIT}, we can also harmlessly  replace $\mfw$ and $\mfr$ into the leading 
error terms on the right to rewrite the above equation as
 \[
\left\{
\begin{aligned}
&  \mfw_{ t} + \underline{b} \mfw_{\alpha}  + \frac{\mfr_\alpha}{1+\bar \W}  + \frac{R_\alpha}{1+\bar \W} \mfw =
   \P [ R \bar \mfw_{ \alpha}- \bar \mfr_{\alpha} \W]-i\frac{c}{2}Princ \, (G_{1,1}) +  \errw \\
\\
&  \mfr_{t} + \underline{b}\mfr_{\alpha}  +ic\mfr- i \frac{(g+\underline{a}\,) \mfw}{1+\W} =  -  \P [\bar \mfr_{\alpha} R]-i\frac{c}{2} Princ \, (K_{1,1})   + \errr .\\
\end{aligned}
\right.
\]
One downside to the conjugation is that the new variables $(\mfw,\mfr)$ are no longer 
holomorphic. To remedy this we project onto the holomorphic space to write a system for the
variables $(w,r)= (\P \mfw,\P \mfr)$. At this point one may legitimately be
concerned that restricting to the holomorphic part might remove a good
portion of our variables. However, this is not the case, as one can
verify the claim by reviewing the discussion in Lemma~3.4 in
\cite{HIT}.  Following again \cite{HIT}, after estimating some 
Coifman-Meyer type commutators, the system for $(w,r)$ is:
\begin{equation}
\label{Pw,Pr}
\left\{
\begin{aligned}
&  w_{ t} + \M_{\underline{b}}  w_{\alpha}  + \P\left[\frac{r_\alpha}{1+\bar \W}\right]  +
\P \left[\frac{R_\alpha}{1+\bar \W}  w \right] =  G 
\\
&\\
&   r_{t} + \M_{\underline{b}} r_{\alpha}  +ic r- i \P \left[ \frac{(g+\underline{a}\,)  w}{1+\W}\right] =K,
\end{aligned}
\right.
\end{equation}
where
\[
\begin{cases}
G:=&  \P [ R \bar w_{ \alpha}- \bar r_{ \alpha} \W]-i\dfrac{c}{2}\P\left[ Princ \, (G_{1,1})\right] + \errw
\\ \\
K:=&  -  \P [\bar r_{\alpha} R]-i\dfrac{c}{2} \P \left[ Princ \, (K_{1,1}) \right] + \errr,
\end{cases}
\]
This is our main system for the (slightly renormalized) differentiated variables
$(\W_\alpha, R_\alpha)$. In order to use it we need to properly relate $(w,r)$ to $ (\W_\alpha, R_\alpha)$,
and to estimate the terms in $G$ and $K$. This is done in the following lemma
\begin{lemma}
\label{l:erori}
a) The energy of $(w,r)$ above is equivalent to that of $(\W_\alpha,R_\alpha)$,
\begin{equation}
\| (w,r)\|_{L^2 \times \dot H^\frac12} \approx_{\uA} \| ( \W_\alpha,R_\alpha  )\|_{L^2 \times \dot H^\frac12}  =  \norm_1,
\end{equation}
and the difference is estimated by 
\begin{equation}
\| (w,r) - ( \W_\alpha,R_\alpha  ) \|_{L^2 \times \dot H^\frac12} \lesssim_{A} \uA \norm_1.
\end{equation}

b) The  error terms on the right in \eqref{Pw,Pr} are bounded,
\begin{equation}\label{bierror-n}
\| ( \P [R  \bar w_\alpha -  \W \bar r_\alpha]  , \P [R  \bar r_\alpha])\|_{L^2 \times \dot H^\frac12} \lesssim_A B
\norm_1, \qquad \| P [R  \bar r_\alpha] \|_{L^2} \lesssim A \norm_1 ,
\end{equation}
respectively
\begin{equation}\label{bierror-nc}
\| ( \P Princ \, (G_{1,1}) , \P Princ \, (K_{1,1}) )\|_{L^2 \times \dot H^\frac12} \lesssim_A A \norm_1, \quad \| \P Princ \, (K_{1,1})  \|_{L^2}  \lesssim_A A_{-1/2}
\norm_1.
\end{equation}
\end{lemma}
\begin{proof}
Part (a) and the first estimate in part (b) are from \cite{HIT}, 
while part (c) follows either directly, for some terms,
or by Coifman-Meyer estimates for the rest. 
\end{proof}

Given the above lemma, the conclusion of Proposition~\ref{t:en=large}
for $n=1$ follows from the energy estimates for the linearized
equation, namely part (a) of Propositions~\ref{plin-short}~\eqref{lin-gen2}; further, if $n=1$ then
we can take
\[
\End(\W,R) = \Elind( w, r).
\]

\subsection{The differentiated equations for \texorpdfstring{$n \geq 2$}\ }

The first step is to derive a set of equations for $(\W^{(n)},
R^{(n)})$.  We start again with the differentiated equations
\eqref{ww2d-diff} and differentiate $n$ times. 

Compared with the case $n=1$, we obtain many more terms. To separate
them into leading order and lower order, we call lower order terms
($lot$) any terms which do not involve $\W^{(n)}$, $ R^{(n)} $ or
derivatives thereof.  In the computation below we take care to
separate all the leading order terms. Toward that end we define again
the notion of \emph{error term}. Unlike in the case $n=1$, here we
also include lower order quadratic terms into the error.  As before,
we describe the error bounds in terms of the parameters $\uA$, $\uB$
and
\begin{equation}
\norm_n = \| (\W^{(n)},R^{(n)})\|_{L^2 \times \dot H^\frac12}.
\end{equation}
The acceptable errors in the $\W^{(n)}$ equation are denoted by
$\errw$ and are of two types, $\errw^{[2]}$ and
$\errw$.  The lower order quadratic holomorphic terms are placed $\errw^{[2]}$,
which is defined to be a linear combination of expressions of the form
\[
\P[\W^{(j)} R^{(n+1-j)}], \quad \P[\bar \W^{(j)} R^{(n+1-j)}], \quad  \P[\W^{(j)} \bar R^{(n+1-j)}], \qquad
2 \leq j \leq n-1,
\]
as well as terms involving the vorticity $c$, namely
\[
\begin{aligned}
&c \P[\W^{(j)} \W^{(n-j)}], \quad c \P[\bar \W^{(j)} \W^{(n-j)}], \quad 1 \leq j \leq n-1,\\
& c \P[R^{(j)} \bar R^{(n+1-j)}], \quad c \P[R^{(j)} R^{(n+1-j)}], \quad 2 \leq j \leq n-1.
\end{aligned}
\]
By interpolation and H\"older's inequality,  terms $G$ in $\errw^{[2]}$ satisfy the bound
\begin{equation*}
\| G\|_{L^2} \lesssim \uB \norm_n.
\end{equation*}

By $\errw^{[3]}$ we denote terms $G$ which  satisfy the estimates the same estimates 
as $\errw$ in the case $n=1$, but with $\norm_1$ replaced by $\norm_n$.

The acceptable errors in the $R^{(n)}$ equation are denoted by
$\errr$ and are also of two types, $\errr^{[2]}$ and
$\errr^{[3]}$. $\errr^{[2]}$ consists of holomorphic quadratic lower
order terms of the form
\[
\P[R^{(j)} R^{(n+1-j)}], \quad \P[\bar R^{(j)} R^{(n+1-j)}],   \qquad
2 \leq j \leq n-1,
\]
and
\[
\P[\W^{(j)} \W^{(n+1-j)}], \quad \P[\bar \W^{(j)} \W^{(n-j)}], \qquad
1 \leq j \leq n-1,
\]
as well as the $c$ terms 
\[
c\P[\W^{(j)} R^{(n-j)}], \quad c\P[\bar \W^{(j)} R^{(n-j)}], \quad  c\P[\W^{(j)} \bar R^{(n-j)}], \qquad
1 \leq j \leq n-1.
\]
By interpolation and H\"older's inequality,  terms $K$ in $\errr^{[2]}$ satisfy the bound
\begin{equation*}
\| K\|_{\dot H^\frac12} \lesssim \uB \norm_n, \qquad \| K\|_{L^2} \lesssim A \norm_n.
\end{equation*}
By $\errr^{[3]}$ we denote terms $K$ which satisfy the same estimates
as $\errr$ in the case $n=1$, but with $\norm_1$ replaced by
$\norm_n$.

We now proceed to differentiate $n$ times the
equation~\eqref{ww2d-diff}.  Our task is simplified due to \cite{HIT},
where this analysis has already been carried out for the terms without
$c$.  Hence, we only concentrate here on
the $c$ terms.

In addition, we remark that, as all terms in the $\W^{(n)}$ equation
have the same homogeneity, whenever all the Sobolev exponents are
within the lower order range, we are guaranteed to get the correct
$L^2$ estimate after interpolation and H\"older's inequality. The same
applies to all cubic terms with at most $n$ derivatives on any single
factor.  The same observation applies to all the lower order terms in
the $R^{(n)}$ equation, as well as to all cubic terms containing
$R^{(n)}$.  However, the terms containing $\W^{(n)}$ are unbounded and
belong to the principal part of the $R^{(n)}$ equation.  Because of
these considerations, the computation below is largely of an algebraic
nature.

We begin with  the $\W$ equation. For the $b_1$ term we have
\[
\begin{split}
\partial^{n} (b_1 \W_\alpha) = & \ b_1 \W^{(n)}_\alpha  + n b_{1,\alpha} \W^{(n)}
+ b_1^{(n)} \W_\alpha 
+ \errw
\\
= & \ \ b_1 \W^{(n)}_\alpha + n (\W-\bar \W)  \W^{(n)} + \errw.
\end{split}
\]
Here at the last step we use the simple observation that the linear part 
of the expression $b_1^{(n)} \W_\alpha$ contributes only lower order terms,
and the rest contributes cubic terms,  which can be estimated either directly
or using Coifman-Meyer bounds.

Next we consider 
 \[
\begin{split}
\partial^{n} ((1+\W) M_1) = & (1+\W) \partial^{n+1}
\left(\P\left[ W\bar{Y}\right]-\bar{\P}\left[ \bar{W}Y\right]\right)
 + \errw
\\
= &\W^{(n)} \W  + (n+1) \W \bar \W^{(n)}+  \P[  W \bar \W^{(n+1)}]     + \errw.
\end{split}
\]
Here all cubic terms involving $\bar \W^{(n+1)}$ can be directly bounded 
using Coifman-Meyer estimates. Finally, for the last $c$ term in the $\W^{(n)}$ 
equation we have
 \[
\begin{split}
\partial^{n} (\W(\W-\bar \W)) = & 2 \W^{(n)}  \W  -  \W^{(n)} \bar \W - \W \bar \W^{(n)} + \errw.
\end{split}
\]

Next we turn our attention to the $c$ terms in the $R^{(n)}$ equation, 
 We begin with
\[
\begin{split}
\partial^{n-1} (b_1 R_\alpha) = & \ b_1 R^{(n)}_\alpha  + n b_{1,\alpha} R^{(n)}
+ b^{(n)} R_\alpha + \errr
\\
= & \ \ b_1 R^{(n)}_\alpha  +   n(\W - \bar \W) R^{(n)} + \errr, 
\end{split}
\]
where again the quadratic terms in $b^{(n)} R_\alpha$ are lower order, and the cubic terms 
with $\bar \W^{(n)}$ are bounded as error terms via Coifman-Meyer estimates.

For the next term in the $R$ equation we write
\[
\begin{split}
\partial^{n} \frac{\W(R+\bar R)}{1+\W} = & \ \W(R^{(n)} + \bar R^{(n)}) +
 \frac{\W^{(n)}(R + \bar R)}{1+\W} + \errr ,
\end{split}
\]
where the cubic $R^{(n)}$ terms are estimated directly as error terms.
Finally, we need to consider
\[
\begin{split}
\partial^{n} \frac{N}{1+\W} = & \partial^{n} \left(\frac{  \P\left[W\bar{R}_{\alpha}-\bar \W R\left]+\bar{\P}\right[\bar{W}R_{\alpha}- \W \bar R\right] }{1+\W}\right) 
\\
=  & \ 
- \frac{N \W^{(n)}}{(1+\W)^2}+ 
\P[W\bar R_\alpha^{(n)} + n  \W\bar R^{(n)} - \bar \W^{(n)} R - \bar \W R^{(n)}] + \errr . 
\end{split}
\]

Combining the above computations for the $c$-terms in \eqref{ww2d-diff} with the prior 
computations for the non-$c$ terms  in \cite{HIT} we obtain the differentiated system
\[
\left\{
\begin{aligned}
 &  \W^{(n)}_{ t} + \ub \W^{(n)}_{ \alpha}  + \frac{((1+\W) R^{(n)})_\alpha}{1+\bar \W}
 + \frac{R_\alpha}{1+\W} \W^{(n)}  =  \uG_n + \errw
\\ &
 R^{(n)}_t + \ub R^{(n)}_\alpha + i c R^{(n)} -  i\left(\frac{(g+ \ua)\W^{(n)}}{(1+\W)^2}\right) =
\uK_n+ \errr,
\end{aligned}
\right.
\]
where 
\[
\uG_n = G_n - i \frac{c}2 G_{n,1}, \qquad \uK_n = K_n -i \frac{c}2 K_{n,1}.
\]
From \cite{HIT} we have
\[
\left\{
\begin{aligned}
G_n = & -  n \frac{\bar R_\alpha}{1+\W}  \W^{(n-1)}
-  (n-1) \frac{ R_\alpha}{1+\bar \W}  \W^{(n-1)}
 + P[ R \bar \W^{(n-1)}_\alpha - \W \bar R^{(n-1)}_\alpha]\\ &
+ R^{(n-1)}(n \bar \W_\alpha - (n-1)\W_\alpha)
+  n(R_\alpha \bar \W^{(n-1)} - \W_\alpha\bar R^{(n-1)})
\\
K_n = & - n\left( \frac{ R_\alpha}{1+\bar \W}+ \frac{\bar R_\alpha}{1+\W}\right)  R^{(n-1)}
 -  \left(P[R \bar R^{(n-1)}_\alpha] - nR_\alpha  \bar R^{(n-1)}\right),
\end{aligned}
\right.
\]
and from the above computations,
\[
\left\{
\begin{aligned}
G_{n,1} =& \   (n+2) [\W \bar \W^{(n)} - (\W-\bar \W) \W^{(n)}] + \P[ W \bar \W^{(n+1)}] 
\\
K_{n,1} = & -(n+1) [(\W - \bar \W) R^{(n)} + \W \bar R^{(n)}] + \P [ R \bar \W^{(n)} - W \bar R_\alpha^{(n)}] . 
\end{aligned}
\right.
\]
To bring this equation to a form closer to the linearized equation we
follow the lead of \cite{HIT} and perform several algebraic
holomorphic substitutions for the $R^{(n)}$ variable, beginning with 
\[
\bR =(1+\W)R^{(n)},
\]
and followed by
\[
\tR = \bR - R_\alpha \W^{(n-1)} +(2n+1) \W_\alpha R^{(n-1)}.
\]
Finally, we conclude with the exponential conjugation by $e^{n\phi}$
where we use the same $\phi$  as before, namely $\phi = -2 \Re \log(1+\W)$, in order to eliminate
the unbounded terms on the right. At the end we obtain an equation for
$(\mfw := e^{n\phi} \W^{(n-1)}, \mfr := e^{n \phi} \tR)$ where the
leading part is exactly as in the linearized equation. We do not
repeat this computation, as it primarily affects the part of the
equation without $c$, which is fully described in \cite{HIT}. All the
additional $c$ contributions are either cubic and bounded or quadratic
and lower order, so they are placed in the error. The outcome
is the equation
\[
\left\{
\begin{aligned}
 & \mfw_{ t} + \ub \mfw_{ \alpha} + \frac{ \mfr_\alpha}{1+\bar \W}
 + \frac{R_\alpha}{1+\W} \mfw  =    \P[ R  \bar \W^{(n)}_\alpha - \W  \bar R^{(n)}_\alpha]
+  (n+1) (R_\alpha  \bar \mfw -  \W_\alpha\bar \mfr)  
\\ &\hspace{3in}\, \ \ \ \ - \frac{ic}2  G_{n,1}
 + \err(L^2)
\\
& \mfr_t + \ub \mfr_\alpha +ic \mfr -  i\left(\frac{(g+ \ua) \mfw}{1+\W}\right) =
-  \P[ R  \bar R^{(n)}_\alpha]  -  (n+1) R_\alpha  \bar \mfr - \frac{ic}2 K_{n,1}  + \errr.
\end{aligned}
\right.
\]
As $(\mfw,\mfr)$ are no longer holomorphic, we project and work with
the projected variables $(w,r) = (\P \mfw,\P \mfr)$.  After some
additional commutator estimates, which are identical to those in the
$n=1$ case, we finally obtain
\begin{equation}\label{energy(n,3)}
\!\!\!\!\! \left\{
\begin{aligned}
 & w_{ t} + \M_{\ub} w_{ \alpha} + \P\! \left[\frac{  r_\alpha}{1+\bar \W}\right]
\! + \P\! \left[\frac{R_\alpha w}{1+\W} \right] \! =   \P[ R  \bar w_\alpha - \! \W  \bar r_\alpha]+ \! (n+1) \P [R_\alpha   \bar w -\!  \W_\alpha\bar r] \!\!\!\!
\\ &\hspace{3in}\, \ \ \ \ - \frac{ic}2 \P G_{n,1}
 + \err(L^2)
\\
& r_t + \M_{\ub} r_\alpha + icr -  i \P\! \left[\frac{(g+ \ua) w}{1+\W}\right] \! =
-  \P[ R  \bar r_\alpha]  -\!  (n+1)\P [R_\alpha  \bar r] - \frac{ic}2 \P K_{n,1}\! + \err(\dot H^\frac12).\!\!\!\!\!
\end{aligned}
\right.
\end{equation}
This is our main system for the (slightly modified) differentiated variables $(\W^{(n)},R^{(n)})$.
 In order to use it we again need to relate $(w,r)$ to $ (\W^{(n)}, R^{(n)})$,
and also to estimate the terms in $G_{n,1}$ and $K_{n,1}$. This is done in the following 

\begin{lemma}
\label{l:erori+}
a) The energy of $(w,r)$ above is equivalent to that of $(\W^{(n)}, R^{(n)})$,
\begin{equation}
\| (w,r)\|_{L^2 \times \dot H^\frac12} \approx_{\uA} \| (\W^{(n)}, R^{(n)})\|_{L^2 \times \dot H^\frac12}  =  \norm_n,
\end{equation}
and the difference is estimated by 
\begin{equation}
\| (w,r) - (\W^{(n)}, R^{(n)}) \|_{L^2 \times \dot H^\frac12} \lesssim_{A} \uA \norm_n.
\end{equation}

b) The  error terms on the right in \eqref{Pw,Pr} are bounded,
\begin{equation}\label{bierror-n+}
\| ( \P [R  \bar w_\alpha -  \W \bar r_\alpha]  , \P [R  \bar r_\alpha])\|_{L^2 \times \dot H^\frac12} \lesssim_A B
\norm_n, \qquad \| \P [R  \bar r_\alpha] \|_{L^2} \lesssim A \norm_n ,
\end{equation}
respectively
\begin{equation}\label{bierror-nn+}
\| ( \P [R_\alpha  \bar w -  \W_\alpha \bar r]  , \P [R_\alpha \bar r])\|_{L^2 \times \dot H^\frac12} \lesssim_A B
\norm_n, \qquad \|\P [R_\alpha \bar r]\|_{L^2} \lesssim_A A
\norm_n,
\end{equation}
and
\begin{equation}\label{bierror-nc+}
\| ( \P G_{1,n} , \P K_{1,n})\|_{L^2 \times \dot H^\frac12} \lesssim_A A, \qquad \| \P K_{1,n}\|_{L^2}  \lesssim_A A_{-1/2}
\norm_n.
\end{equation}
\end{lemma}
\begin{proof}
Part (a) and the first two bounds in part (b) are from \cite{HIT}, while the last bound in part  
(b) follows either directly, for some terms, or by Coifman-Meyer estimates, for the rest. 
\end{proof}

Given the above Lemma~\ref{l:erori}, the $n \geq 2$ case of the result
in Proposition~\ref{t:en=large} is a direct consequence of our
quadratic estimates for the linearized equation in
Proposition~\ref{plin-short}(a).

Comparing the linear $(w,r)$ equation obtained above for the case $n
\geq 2$ to the corresponding equation previously obtained for the case
$n=1$, we note that here we have two extra terms, namely the ones
estimated in \eqref{bierror-nn+}.  On one hand this is good, because
the bound \eqref{bierror-nn+} is no longer true if $n=1$. But on the
other hand, it means that we will no longer be able to use directly
the cubic energy $ \Elint(w,r)$ introduced in the context of the
linearized equation as the full high frequency part of our cubic
energy. To account for these extra terms, we will add a further
correction to the energy $\Elint(w,r)$, and define
\begin{equation}
\label{ent-high}
\Enthigh  ( w,r) :=
\int (1+a)  |w|^2 + \Im (\bar r r_\alpha)
 + 2n \Im (R_{\alpha} \bar w  \bar r) +
 2(  \Im[\bar{R} w r_\alpha] - \Re[\bar {W}_\alpha w^2])
\, d\alpha.
\end{equation}
Then we will prove a counterpart to the small data (nearly) cubic
energy estimates in Proposition~\ref{t:en=small}.  Precisely, we claim
that the evolution of this energy is governed by the following

\begin{lemma}\label{l:hf-cubic}
Let $(w,r)$ be defined as above. Then

a) Assuming that $A \ll 1$,  we have
\begin{equation}\label{en3-equiv}
 \Enthigh ( w, r) \approx \Ez  ( w, r)
\approx  \norm_n,
\end{equation}

b) The solutions $(w,r)$ of \eqref{energy(n,3)} satisfy
\begin{equation}\label{en3-evolve}
\begin{split}
\frac{d}{dt} \Enthigh( w, r) =  & 
2 \int  \! - \Re (\bar w \cdot ( \frac{ic}2 G_{n,1}\! -  \! \err(L^2)^{\left[ 2\right] } ) + \Im(\bar r_\alpha \cdot ( \frac{ic}2 K_{n,1} \! - \! \err(\dot H^\frac12)^{\left[ 2\right] })  
\, d\alpha \\
& \ \ \  + \int  \frac{c^2}{2}\Im R |w|^2  \, d\alpha   +  O_A(A \uB \norm_n).
\end{split}
\end{equation}
Further, the same relation holds if $(\bar w, \bar r)$ on the right are replaced
by $(\bar \W^{(n)},\bar R^{(n)})$.
\end{lemma}
We note that in the proof of this lemma we crucially use the fact that $(w,r)$ are the ones
associated to the differentiated equation, and not arbitrary functions in the energy space.
This is also the reason why this lemma is presented here, rather then in Section~\ref{s:linearized}.

\begin{proof}
Part (a) was already proved in \cite{HIT}, so we proceed to part (b).
 Here, we begin with the cubic linearized energy, $\Elint$. According to
the bound \eqref{lin-gen2+} in Proposition~\ref{plin-short3}, we have
\[
\begin{split}
\frac{d}{dt} \Elint( w, r) = & \
 \int  2\Re \left( (n \P [R_\alpha  \bar w -  \W_\alpha \bar r] -\frac{ic}2 G_{n,1} + P\err(L^2))
\cdot ( \bar w - \bar \P[\bar R r_\alpha] - \bar \P[\bar W_\alpha w])  \right)
\\ & \ \  - 2 \Im\left( (-  n\P [R_\alpha  \bar r] - \frac{ic}2 K_{n,1} + \P \err(\dot H^\frac12)) \cdot (\bar r_\alpha
+ \bar \P[\bar R w]_\alpha)\right)
\, d\alpha \\
& \ \ +  O_A\left( \uA\uB \|(w,r)\|_{L^2 \times \dot H^\frac12}^2\right) .
\end{split}
\]
By the Coifman-Meyer type estimates in Lemma~\ref{l:com} the following bounds hold:
\begin{equation}
\|  \bar P[\bar R r_\alpha]\|_{L^2} + \|   \bar P[\bar W_\alpha w]\|_{L^2}
+\|  \bar P[\bar R w]\|_{\dot H^\frac12}
   \lesssim A \|(w,r)\|_{L^2 \times \dot H^\frac12}.
\end{equation}
Combining this with \eqref{bierror-n} and with the bounds for $G_{n,1}$, $K_{n,1}$
 and the error terms we get
\[
\begin{split}
\frac{d}{dt} \Elint( w,r) \leq & \
 \int  2\Re \left( (n \P [R_\alpha   \bar w -  \W_\alpha\bar r] -\frac{ic}2 G_{n,1}+ \P\err(L^2)^{\left[ 2\right] })
\cdot \bar w \right)
\\ & \ \ - 2\Im\left( (-  n\P [R_\alpha \bar r] - \frac{ic}2 K_{n,1}+ \P \err(\dot H^\frac12)^{\left[ 2\right] })
\cdot  \bar r_\alpha \right) \, d\alpha \\
& \ \ +  O_A\left( \uA\uB\right)  \|(w,r)\|_{L^2 \times \dot H^\frac12}^2 ,
\end{split}
\]
where the output from all error terms which are cubic and higher error terms
is all included in the last RHS term.

It remains to consider the contribution of the extra term
in $\Enthigh$ and show that
\begin{equation}\label{extra-diff}
\begin{split}
\frac{d}{dt} \int  \Im (R_{\alpha}  \bar w   \bar r) \, d\alpha = & \
 \int  \Re \left((R_\alpha   \bar w -  \W_\alpha \bar r)  \bar w \right)
  + \Im\left(   R_\alpha  \bar r
  \bar r_\alpha \right) \, d\alpha
\\ & \ \ \,+  O_A\left( \uA\uB\right) \|(w,r)\|_{L^2 \times \dot H^\frac12}^2.
\end{split}
\end{equation}
Denote by $G_n$, respectively $K_n$, the two right hand sides in
\eqref{energy(n,3)}. By the definition of error terms and by
\eqref{bierror-n}-\eqref{bierror-nc} they satisfy the bounds
\[
\|(G_n,K_n)\|_{L^2 \times \dot H^{\frac12}} \lesssim_A \uB  \norm_n, \qquad \|K_n\|_{L^2} \lesssim_A 
\uA  \norm_n.
\]
Then their contribution in the above time derivative is estimated
\[
\left| \int  \Im (R_{\alpha} \bar \P \bar G_n   \bar r + R_{\alpha}  \bar w  \bar \P \bar K_n)\ d\alpha\right| = \left| \int  \Im (R_{\alpha} \bar \P \bar F_n   \bar r + P[R_{\alpha}  \bar w]  \bar \P \bar K_n)\ d\alpha\right| \lesssim_A \uA \uB \norm_n,
 \]
by using H\"older's inequality for the first term and the Coifman-Meyer commutator
estimate in Lemma~\ref{l:com} for the second.

We now consider the contributions due to the terms on the left of \eqref{energy(n,3)}.
These were already estimated in \cite{HIT} when $c=0$, so we only need to estimate the 
$c$ terms, i.e., the contributions of $b_1$ and $a_1$.

The contributions of the $b_1$ terms  are
collected together in the real part of the expression
\[
\begin{split}
I = & \  \frac{c}2 \int  \partial_\alpha (b_1 R_{\alpha})  \bar w   \bar r
+ R_\alpha  \bar \P(b_1 \bar w_\alpha)   \bar r  +  R_\alpha
  \bar w \bar \P(b_1 \bar r_\alpha )\, d\alpha \\  = & \
\frac{c}2 \int
 R_\alpha   \left([b_1,\P]( \bar w_\alpha) \, \bar r  +
  \bar w [b_1,\P] ( \bar r_\alpha )\right)\, d\alpha.
\end{split}
\]
Since $\|b_{1,\alpha}\|_{BMO} \lesssim A$, we can bound this using Lemma~\ref{l:com}, and
 then use H\"older's inequality for all terms.

The contribution of $a_1$ is  
\[
\frac{c}2 \Re \int R_\alpha \bar w  P \left[\frac{a_1 \bar w}{1+\W}\right] d\alpha
\]
for which it suffices to use the pointwise bound for $a_1$ in Lemma~\ref{l:a1} 
and the classical Coifman-Meyer bound.

Finally,  we consider the remaining contribution of the $c$ terms in the time
derivative of $R_\alpha$, for which we use the equation \eqref{ww2d-diff}
\[
\frac{c}{2} \Re \int   \bar w   \bar r \partial_\alpha  \frac{R\W + \bar R \W +N}{1+\W}  \, d\alpha .
\]
But this is a quartic expression for which we only need Sobolev embeddings 
and interpolation. 
\end{proof}

\subsection{Cubic  small data energy estimates: 
\texorpdfstring{$n \geq 1$} \  }
\label{s:ee3+} 

In this section we construct an $n$-th order energy functional $\Ent$ with cubic
estimates, for $n \geq 1$. While in essence the argument does not depend on $n$,
there are nevertheless a few small differences between the case $n=1$ 
and $n \geq 2$. These differences will be pointed out  at the appropriate 
places.

One ingredient in the construction of $\Ent$ is the high frequency (nearly) cubic energy
$\Enthigh$ in Lemma~\ref{l:hf-cubic} ($n \geq 2$), respectively $\Elint$ in Proposition~\ref{plin-short3}
($n = 1$).  However, this does not
suffice, as the right hand side of the energy relation
\eqref{en3-evolve} still contains lower order cubic terms.
To remedy this, here  we use  normal forms in order to add a lower order correction
to $\Enthigh$, which removes the above mentioned cubic terms.
For this we follow the method introduced in \cite{HIT}, though the computations here
are significantly more involved.

We first provide an outline of the argument, and then return to each
step in detail. The three main steps are as follows:

\bigskip

{\bf 1. Construct the normal form energy.}
The normal form variables $\tW,\tQ$ solve an equation of the form 
\begin{equation}
\left\{
\begin{aligned}
&\tW_t + \tQ_\alpha = \text{cubic and higher}
\\
&\tQ_t - i \tW =  \text{cubic and higher}.
\end{aligned}
\right.
\label{nft2eq}
\end{equation}
Thus, its associated linear energies satisfy 
\[
\frac{d}{dt} \Ez ( \tW^{(n+1)},\tQ^{(n+1)}) = \text{quartic and higher}.
\]
Here the quartic part of the left hand side  is highly unbounded,
and it is not uniquely  determined, as one can add arbitrary cubic terms to the 
normal form. Fortunately, it is also irrelevant for the above relation.

Thus, we eliminate it, and we define the normal form energy as
\[
\Ennf = \text{ (quadratic + cubic)}  [ \Ez ( \tW^{(n+1)},\tQ^{(n+1)})] .
\]
Here, we carefully make the choice of choosing the quadratic 
and cubic terms in the expansion of $\Ez ( \tW^{(n)},\tQ^{(n)})$ in terms 
of the diagonal variables $(\W,R)$, rather than $(W,Q)$.

The expression $\Ennf$  has only quadratic and cubic terms
in $(\W,R)$, and still satisfies 
\[
\frac{d}{dt} \Ennf = \text{unbounded quartic and higher}.
\]
Its disadvantage, due to the fact that our problem is quasilinear, is 
that the terms on the right are highly unbounded.

\bigskip

{\bf 2. Construct the  quasilinear modified energy.}
Here we construct \emph{the quasilinear modified energy} $\Ent$, using the 
normal form energy $\Ennf$ and the high frequency modified 
energy $\Enthigh(w,r)$. The first one has quartic estimates,
which are unbounded. The second one has good (quartic) high 
frequency estimates, but with cubic lower order errors,
\[
\frac{d}{dt} \Enthigh(w,r) = \text{bounded quartic and higher} + \text{lower order cubic}.
\]
We seek to combine the two into a quasilinear modified energy $\Ent$ which satisfies
\[
\frac{d}{dt} \Ent = \text{bounded quartic and higher}.
\]
To achieve this, it is natural to separate $\Ennf$ into a leading part and a lower order part,
\[
\Ennf = \Ennfhigh + \Ennflow,
\]
so that the leading part coincides with the high frequency energy 
up to quartic terms,
\[
 \Enthigh = \Ennfhigh + \text{quartic and higher},
\]
and then define 
\[
\Ent =  \Enthigh+ \Ennflow.
\]
This was the argument in \cite{HIT}; here it is slightly more
complicated, as additional terms in $\Ennflow$ need to be treated
as leading order, and thus corrected in a quasilinear fashion.
\bigskip

{\bf 3. Show that  $ \Ent $ is a good quasilinear cubic energy}, i.e., prove that
the bound \eqref{en-cubic} holds. Here, by construction, we know that 
\[
\frac{d}{dt} \Ent = \text{quartic and higher}; 
\]
therefore we can write 
\[
\frac{d}{dt} \Ent = \text{quartic and higher}(\frac{d}{dt} \Enthigh(w,r) )
+ \text{quartic and higher}( \frac{d}{dt} \Ennflow).
\]
The favorable bound for the first part is already a direct consequence 
of Lemma~\ref{l:hf-cubic}, so it remains to use the equations in order to bound the 
second term. But this is technically simple, and the argument is more a matter
of bookkeeping.

\bigskip

Now we  proceed to implement the above strategy.

\bigskip

{\bf STEP 1.} 
The first step described above is implemented in the following proposition:

\begin{proposition}
There exists a modified normal form energy $\Ennf=\Ennf(W,R) $ 
with the following properties:

A. $\Ennf$ has the form
\[
\Ennf =  \Ennfhigh +  \Ennfhighc   + \Ennflow,
\]
where the three components are defined as follows:
\begin{align}\label{ennf-high}
\begin{split}
\Ennfhigh = & \ \int  (1- 4(n+1) \Re \W)   \left(g |\W^{(n)}|^2
+  \Im[\bar R^{(n+1)} R^{(n)}]\right) 
d\alpha,
\\ & \ + 2 \int \Im[\bar{R}{\W}^{(n)} R^{(n+1)}] - g \Re[\bar {\W} (\W^{(n)})^2]  + (n+1) \Im[R_{\alpha} \bar \W^{(n)} \bar R^{(n)}]\, d\alpha
\\& \  \int c \Re R |\W^{(n)}|^2 + 2 \Im (\W \bar R^{(n+1)} R^{(n)}) \, d \alpha,
\end{split}
\end{align}
where the last term in the second integral is dropped if $n=1$,
\begin{align}
\label{ennf-high+}
\begin{split}
\Ennfhighc =  -  \int  \left [ c (2n+3)  \Re R  + c^2(2n+\frac52) \Im W\right] 
 (|\W^{(n)}|^2 - i \bar R^{(n+1)} R^{(n)})  d\alpha .
\end{split}
\end{align}
Finally, $\Ennflow$ is a trilinear integral form of order $2n$, with the following 
restrictions: 
\begin{enumerate}[label=(\roman*)]
\item The leading order for the terms without $c$ is at most $2n-1$.

\item The leading order for the terms with $c$ is at most $2n-\frac12$.

\item The highest power of $c$ is $4$. 

\item At most one factor of the form $c^2 W$ or $cR$ is present.
\end{enumerate}

B.  The energy functional $\Ennf$ is cubically accurate, 
\begin{equation}\label{enlow-diff}
\Lambda^{\leq 3} \frac{d}{dt} \Ennf = 0.
\end{equation}

C. Its components satisfy the following estimates: 
\begin{equation}\label{enlow-equiv}
\begin{aligned}
\Ennfhigh =& [1+ O(A)] \Ez (\W^{(n)}, R^{(n)}),  \qquad  \Ennfhighc = O(\uA) \Ez (\W^{(n)}, R^{(n)}), \\
 &\Ennflow = O(\uA)( \Ez (\W^{(n)}, R^{(n)})  + c^4\Ez (\W^{(n-1)}, R^{(n-1)})) .
\end{aligned}
\end{equation}
\end{proposition}

\begin{proof}
  Here we use the algebraic expression for the normal form
  transformation, which is given below. We call lower order terms
  (``lot'') any terms that can be included in $\Ennflow$.  Up to
  quartic terms we seek to have, at least formally, the relation
\begin{equation}\label{start-nf}
\begin{split}
E^n_{NF} = & \ \| (\tW^{(n+1)},\tQ^{(n+1)}) \|_{L^2 \times \dot H^\frac12}^2 + \text{quartic}
\\
= & \ \Ez(W^{(n+1)},Q^{(n+1)}) - 2 \Re \int \bar  W^{(n+1)} \partial^{n+1}W^{[2]} - i \bar Q^{(n+1)} 
 \partial^{n+1}Q^{[2]} d \alpha + \text{quartic}.
\end{split}
\end{equation}
In the first term we introduce the diagonal variable $R$ by using $Q_\alpha = R(1+\W)$,
which allows us to write it as 
\[
I_{main} =  \  \Ez (\W^{(n)},R^{(n)}) + 2 \Im \int  \bar R^{(n+1)}  \partial^{n} (R \W) d \alpha
+ \text{quartic}.
\]
Differentiating we obtain its principal part, containing all terms with leading order 
above $2n-1$:
\[
\begin{split}
 I_{main,high} =   \Ez (\W^{(n)},R^{(n)}) + 2 \Im \int   \W \bar R^{(n+1)} R^{(n)} 
+   R  \bar R^{(n+1)} \W^{(n)} - (n+2)   R_\alpha \bar R^{(n)} \W^{(n)} 
d \alpha .
\end{split}
\]
Here the last term on the right no longer appears if $n=1$, which accounts
for the similar modification in the proposition.

In the integral in \eqref{start-nf}, on the other hand, we have two concerns. First, we
want to eliminate its unbounded part of low frequency, i.e., all terms
with inverse derivatives of $W$, as well as all terms with
undifferentiated $Q$; we remark that once this is done, the switch to
diagonal variables is achieved simply by substituting $Q_\alpha$ with
$R$.  Secondly, we want to keep track of its highest order terms,
i.e., those terms which are at least energy strength  (i.e., $\W^{(n)} \bar \W^{(n)}$, 
or $\bar R^{(n+1)} R^{(n)}$, and also $\bar R^{(n)} \W^{(n)}$ for terms without $c$).
First we recall the expression for the normal form, see \eqref{nf-eq}:
\begin{equation*}
\left\{
\begin{aligned}
W^{[2]} = & \ - (W+\bar W) W_\alpha - \frac{c}{2g} \left[ (Q+\bar Q) W_\alpha + (W + \bar W) Q_\alpha\right] 
\\ & \ + \frac{ic^2}{2g} \left[ (\partial^{-1} W - \partial^{-1} \bar W) W_\alpha + W^2 +\frac12 |W|^2 \right]
- \frac{c^2}{4g^2} (Q+\bar Q) Q_\alpha 
\\ & \ + \frac{ic^3}{4g^2} \left[ (Q+\bar Q) W + (\partial^{-1} W - \partial^{-1} \bar W) Q_\alpha \right] + 
 \frac{c^4}{4g^2} (\partial^{-1} W - \partial^{-1} \bar W) W
\\
Q^{[2]} = & \ - (W+\bar W) Q_\alpha - \frac{c}{2g} (Q+\bar Q) Q_\alpha + \frac{ic}{4} (W^2+ 2|W|^2)
\\ & \ + \frac{ic^2}{2g} \left[ (\partial^{-1} W - \partial^{-1} \bar W) Q_\alpha + \frac12  (Q+\bar Q) W \right]
+  \frac{c^3}{4g} (\partial^{-1} W - \partial^{-1} \bar W) W .
\end{aligned}
\right.
\end{equation*}

Next, we successively consider the terms in the integral in \eqref{start-nf} organized by the power of $c$.

\medskip
{\em 0. Terms with $c^0$.}
These are
\[
I_0 = 2 \Re \int  - \bar \W^{(n)} \partial^{n+1} [ (W + \bar W) W_\alpha] 
+ i \bar Q^{(n+2)}  \partial^{n+1} [ (W + \bar W) Q_\alpha] d \alpha .
\]
Here there are no low frequency issues. Its high frequency part, on the other  
hand, has the form 
\[
\begin{split}
I_{0,high} = &  2 \Re \int - (2n+2) \Re \W ( |\W^{(n)}|^2 - i \bar Q^{(n+2)} Q^{(n+1)}) + i Q_\alpha \bar Q^{(n+2)} 
( W^{(n+1)} + \bar W^{(n+1)})  
\\ & \  \ \ \ - (n+2) Q_\alpha \bar Q^{(n+1)} \W^{(n)}   d \alpha  .
\end{split}
\]

\medskip
{\em 1. Terms with $c$.}
These are
\[
\begin{split}
I_1 = & \  c  \Re \int - \bar  \W^{(n)} \partial^{n+1}\left(\W (Q+\bar Q) +  Q_\alpha(W+\bar W)\right)
+  \bar Q^{(n+2)} \partial^{n+1}(\frac12 W^2 + |W|^2) 
\\ & \ + \frac{i}g \bar Q^{(n+2)}\partial^{n+1}(Q Q_\alpha + \bar Q Q_\alpha) d \alpha .
\\
\end{split}
\]
We first verify that the terms with undifferentiated $Q$ disappear
when integrating by parts. These are 
\[
c \Re \int -  2 \Re Q  \bar \W^{(n)} \W^{(n+1)} d \alpha = c \int \Re Q_\alpha |W^{(n)}|^2 d\alpha,
\]
respectively 
\[
\frac{c}{g} \Re \int i \Re Q  \bar Q^{(n+2)}  Q^{(n+2)} d \alpha = 0.
\]
Secondly, we compute the high frequency component. Taking into account the above integration 
by parts, this is 
\[
I_{1,high} = c \int - (2n+2) |\W^{(n)}|^2  \Re Q_\alpha   
+ \frac{i}{g} (2n+3) \Re Q_\alpha \bar Q^{(n+2)} Q^{(n+1)}  d\alpha.
\]

\medskip
{\em 2. Terms with $c^2$.}
These are 
\[
\begin{split}
I_2 = & \  c^2  \Re \int i \bar \W^{(n)} \partial^{n+1}\left( \W (\partial^{-1}W -  \partial^{-1}\bar W)
+ W^2  + \frac12|W|^2\right)
 - \frac1{2g} \bar \W^{(n)} \partial^{n+1}(Q Q_\alpha   + \bar Q Q_\alpha) 
\\ & \   + \frac{1}{g}  \bar Q^{(n+2)}  \partial^{n+1}( Q_\alpha (\partial^{-1} W  - \partial^{-1} \bar W))
+ \frac{1}{2g} \bar Q^{(n+2)}\partial^{n+1} (WQ + W \bar Q)  \,
 d\alpha.
\end{split}
\]
For the low frequency analysis
we compute
\[
c^2 \Re \int  - i \bar \W^{(n)} \W^{(n+1)}  (\partial^{-1}  W - \partial^{-1} \bar W) d \alpha = 
- c^2 \int |\W^{(n)}|^2 \Im W d \alpha,
\]
while the remaining $\partial^{-1} W$ terms and the $Q$ terms directly cancel.

Now we look at the high frequency part. This is 
\[
\begin{split}
I_{2,high} = & \ c^2 \int - (2 n+ \frac{5}2)  \Im W 
  \left( |\W^{(n)}|^2 - i \bar Q^{(n+2)} Q^{(n+1)} \right)   -  
(2n+4) (Q_\alpha + \bar Q_\alpha) 
\bar \W^{(n)} Q^{(n+1)} \,
d\alpha.
\end{split}
\]

\medskip
{\em 3. Terms with $c^3$.}
These are 
\[
\begin{split}
I_3 = - \frac{c^3}{2g} \Im \int & \  \bar \W^{(n)} \partial^{n+1} \left[ W(Q+\bar Q) + 
(\partial^{-1} W -\partial^{-1} \bar W) Q_\alpha  \right]
\\ & -  \bar Q^{(n+2)}  \partial^{n+1}\left[  (\partial^{-1} W -\partial^{-1} \bar W) W \right]  d\alpha .
\end{split}
\]
For the low frequency analysis we compute
\[
\Re \int  i \bar \W^{(n)} (\partial^{-1} W -\partial^{-1} \bar W) Q^{(n+2)} - i \bar Q^{(n+2)}  (\partial^{-1} W -\partial^{-1} \bar W) \W^{(n)} d\alpha = 0,
\]
and 
\[
\Re \int i  \bar \W^{(n)}  \W^{(n)} (Q + \bar Q) d\alpha = 0 ,
\]
while for the high frequency analysis, taking the above into account,  we are left with
\[
I_{3,high} =   \frac{c^3}{2g} \int i (2n+4)(W-\bar W) \bar \W^{(n)} Q^{(n+1)}  d \alpha .
\]

\medskip
{\em 4. Terms with $c^4$.}
This is 
\[
I_4 = \frac{c^4}{2g^2} \Re \int \W^{(n)} \partial^{(n+1)}[ (\partial^{-1} W- \partial^{-1} \bar W) W] d \alpha,
\]
which exhibits cancellation at low frequency and has no high frequency component.

\bigskip

The normal form energy $\Ennf$ is obtained by substituting
$Q_\alpha$ by $R$ in the all the above terms, after eliminating the
$Q$ and $\partial^{-1} W$. The expressions of $\Ennfhigh$ and
$\Ennfhighc$ are simply obtained by adding up all the high frequency
contributions above.

Since $\Ennf$ is obtained from the normal form, the relation
\eqref{enlow-diff} immediately follows.  We note that establishing this
property is a purely algebraic computation, which does not require the
direct use of the normal form (which is unbounded).

\bigskip

{\bf STEP 2.}
We begin with the  quasilinear energy $\Enthigh (w,r) $, which was defined in 
\eqref{ent-high} if $n \geq 2$, and is set to be equal to $\Elint(w,r)$ if $n=1$. 
We first compare it with the cubic energy  $\Ennfhigh$.  Precisely, we have
 \begin{equation}
\Lambda^{\leq 3}\Elint (w,r)  = \Lambda^{\leq 3} \Ennfhigh, \qquad n = 1,
\end{equation}
respectively
\begin{equation}
\Lambda^{\leq 3} \Enthigh (w,r)  = \Lambda^{\leq 3} \Ennfhigh, \qquad n \geq 2.
\end{equation}
This computation was already done in \cite{HIT} and is not repeated 
here. The only difference is in the $c$ term in $\Ennfhigh$,
which arises due to the  $a_1$ in  $\Elint (w,r) $, respectively  $\Enthigh (w,r)$.

Next, we consider the term $\Ennfhighc$. This contains the leading order 
energy density, so we have to treat it in a quasilinear manner.  Precisely, we replace it by 
\[
 \Enthighc(w,r) =  -  \int  \left [ c (2n+3)  \Re R  + c^2(2n+\frac52) \Im W\right] 
 ((g+ \ua) |w|^2 - i \bar r_\alpha r)  d\alpha, 
\]
which admits a straightforward bound
\begin{equation}
|  \Enthighc(w,r) | \lesssim c A_{-\frac12} \norm_n + c^2 A_{-1} \norm_n.
\end{equation}
Now we are in a position to define our quasilinear modified cubic energy
\begin{equation}
\Ent (\W,R) = \Enthigh(w,r) +  \Enthighc(w,r) + \Ennflow (\W, R).
\end{equation}
By construction this satisfies 
\begin{equation}\label{eq=3}
\Lambda^{<3} \Ent (W,R)  = \Lambda^{<3} \Ennf(W,R).
\end{equation}

\bigskip

{\bf STEP 3.}  Now we proceed to prove Proposition~\ref{t:en=small}. The norm equivalence
is already known from \eqref{en3-equiv} and \eqref{enlow-equiv}, so we
still need the energy estimate.   For this we use \eqref{eq=3} and \eqref{enlow-diff} 
to see that the cubic terms vanish,
\[
\Lambda^{\leq 3} \frac{d}{dt} \Ent  = \Lambda^{\leq 3} \frac{d}{dt} \Ennf = 0.
\]
Hence, it remains to look at the quartic terms and higher. For this we can split
the terms,
\[
\Lambda^{\geq 4} \frac{d}{dt} \Ent 
= \Lambda^{\geq 4} \frac{d}{dt} E^{n,(3)}(w,r) +   
\Lambda^{\geq 4} \frac{d}{dt} \Enthighc(w,r) +  \Lambda^{\geq 4} \frac{d}{dt} \Ennflow.
\]
The bound for the first term on the right is a direct consequence of 
Proposition~\ref{plin-short} for $n=1$, respectively 
Lemma~\ref{l:hf-cubic} for $n \geq 2$. From here on, the cases $n=1$ and $n \geq 2$
are identical.

 For the middle term $\Lambda^{\geq 4} \frac{d}{dt} \Enthighc(w,r) $ 
we will use Lemma~\ref{weighted-en},  with 
either $f = W$ or $f= R$, and $(w,r)$ as in this section. We claim that this 
yields
\[
\left| \frac{d}{dt} \Enthighc(w,r) \right| \lesssim_{\uA} \uB N_n ,
\]
which by homogeneity considerations yields
\[
\left| \Lambda^{\geq 4} \frac{d}{dt} \Enthighc(w,r) \right| \lesssim_{\uA} \uA \uB N_n ,
\]
as desired. 

a) If $f = W$ then the bounds
\[
\| f\|_{L^{\infty}} \lesssim_A A_{-1}, \qquad \| D^\frac12 f\|_{L^{\infty}} \lesssim A_{-1/2}
\]
allow us to estimate the difference $I$ in the Lemma~\ref{weighted-en}. The remaining terms
\[
c^2 \int \Lambda^{\geq 2} (\partial_t+ \ub \partial_\alpha) W ((g+ \ua) |w|^2 - i \bar r_\alpha r) d \alpha,
\]
and
\[
c^2  \int   f ((g+ \ua) \bar w F  - i \bar r_\alpha  G) d\alpha
\]
are $c$ times forms which are at least quartic, with order
$2n+\frac12$ and leading order $2n$ (this is because we have exactly
one $R$ factor). Hence we can bound them using  H\"older's inequality  and interpolation
by $c A \uA  N_n$.

b) If $f = R$ then the bounds
\[
\| f\|_{L^{\infty}} \lesssim_A A_{-1/2}, \qquad \| D^\frac12 f\|_{L^{\infty}} \lesssim A
\]
allow us to estimate the difference $D$ in the lemma. The remaining terms are 
exactly as above, still bounded by $c A \uA  N_n$.

Finally, for the lower order terms $\Lambda^{\geq 4} \frac{d}{dt} \Ennflow$ 
it suffices to have the  following property:

\begin{lemma}
For any term $I$ in $\Ennflow$ we have
\begin{equation}
\label{dt-low}
\left |\Lambda^{\geq 4} \frac{d}{dt} I \right| \lesssim  \uA \uB (N_n+c^4 N_{n-1}) .
\end{equation}
\end{lemma}
\begin{proof}
The terms without $c$ in this computation were already estimated in \cite{HIT}, so we need 
to consider only $c$ terms. A simple order analysis for these terms leads to the following 
properties:
\begin{enumerate}
\item[(a)] Their order is $2n + \frac12$, 
\item[(b)] Their leading order is at most $2n$,
 \item[(c)] The highest power of $c$ is $4$.
\end{enumerate}
Here part (b) is a obtained by using the relation~\eqref{dt-tri} to compute time derivatives.
This guarantees that the leading order does not increase by a full unit.

We also take a closer look at the $c^2W$ and $c R$ factors that can
arise. In the trilinear form $\Ennflow$ they can arise just once, so
the question is what happens if we differentiate with respect to
time. If $W$, respectively $R$ are differentiated in time then at most
we obtain another $W$ factor, respectively an $R$ or $cW$ factor, so
this pattern is preserved. If instead we differentiate a higher order
$W^{(k)}$ or $R^{(k)}$ factor using the equations \eqref{ww2d-diff},
we can further produce undifferentiated $W$ and $R$ factors, as
follows:

(i) Arising from $cM_1$ and its derivatives. There $cW$ will appear in combinations 
of the form $\P[W \bar Y^{(k)}]$ and its conjugate. However, the frequency localization 
enforced by $\P$ guarantees that $W$ is the higher frequency factor, which, for estimates,
ensures that we can freely move derivatives from $\bar Y$ to $W$.

(ii) Arising from derivatives of $cN$, where, for the same reasons as above, we can disregard
the $W$ factors, but we can still produce a $cR$ factor.    

(iii) Arising from $b_1$ and its derivatives. In the case of
derivatives of $b_1$ the discussion concerning $W$ factors is again
identical to case (i), so we can neglect those. Hence the only
potential $W$ factor can arise from undifferentiated $b_1$. However,
these are  avoided due to the use of \eqref{dt-tri}, where the transport is
fully included in the time differentiation and the undifferentiated
advection coefficient $\ub$ never appears.

\medskip

To summarize the above discussion, in the time derivative of $\Ennflow$ we can have 
at most one $c^2 W$ factor and  at most two $c^2W,cR$ factors. 
A further simplification arises from the constraint (b) above. 
Precisely, there  we can integrate by parts to rebalance derivatives in such a way  
that either 

a) both leading order terms have order at most $n$, or

b) we have exactly the product  $R^{(n+1)} \bar R^{(n)}$ as part of our 
multilinear expression.

In case (a), our estimates follow directly by H\"older's inequality and interpolation.
In case (b), the remaining factors must have order zero, except for one which has 
order $\frac12$. There we consider two scenarios.

b1) If the factors $ R^{(n+1)}$ and  $\bar R^{(n)}$ are not frequency balanced, then another
factor has the highest frequency and we can rebalance derivatives and use 
again H\"older's inequality and interpolation. 

b2)  If the factors $ R^{(n+1)}$ and  $\bar R^{(n)}$ are  frequency balanced, then all we need 
is to bound the remaining factors in $L^\infty$. The only factor of order $\frac12$ which we do not control in 
in $L^\infty$ is $R_\alpha$, so  H\"older's inequality and interpolation work unless 
our multilinear expression exactly contains 
\[
R^{(n+1)} \bar R^{(n)} R_\alpha,
\]
and possibly other zero order factors. Backtracking, the only way to
produce such a term is by differentiating a cubic expression which has a
leading order exactly $2n - \frac12$ as well as an $R_\alpha $
factor. Then this expression must be exactly
\[
\W^{(n)} \bar R^{(n)} R_\alpha,
\]
which cannot contain any $c$ factors so it is not within our purview here. 
\end{proof}

\end{proof}

\section{ Proofs of the main theorems}
\label{s:proofs}

The results in Theorem~\ref{baiatul} and Theorem~\ref{t:cubic} are a consequence of the estimates 
for the linearized equation in Section~\ref{s:linearized}, and of the energy estimates for the 
higher Sobolev norms of the solutions in Section~\ref{s:ee}. The arguments here closely 
match their counterparts in our previous gravity waves paper \cite{HIT}. Because of this,
we only provide an outline here, and refer the reader to \cite{HIT} for more details. 

\begin{proof}[Outline of proof of Theorem~\ref{baiatul}]
The steps of the proof are as follows:
\medskip

{\bf Step 1 : Existence of regular solutions.} A standard approach here
is to obtain solutions as limits of solutions to frequency truncated
equation.  As always, 
this truncation needs to happen in a symmetric way, so that uniform
energy estimates survive. In addition, a specific difficulty we
encounter in water wave evolutions is the fact that the evolution we
are trying to truncate is degenerate hyperbolic, a condition which
might not survive in a naive direct truncation procedure. In
\cite{HIT} we bypass this difficulty by solving directly the
differentiated system for the diagonal variables $(\W,R)$.  In our
case this is not entirely possible, as the system for $(\W,R)$ is not
fully self-contained. Precisely, it contains $W$ as a part of $\ub$.
Fortunately, this does not cause extra difficulties, as we can make the
frequency truncation consistent with differentiation in the $W$
equation, and thus be able to freely use $W$ in the $(\W,R)$ system.
We note that formally $\W$ also appears undifferentiated in $M_1$ 
and $N$, but the commutator structure of those expressions
implies that they actually can be rewritten (and estimated) in terms of $\W$
instead. 

Precisely, our main truncated system for $(\W,R)$ has the form
\begin{equation} \label{ww2d-diff2}
\left\{
\begin{aligned}
 & \W_{ t} + \P_{<\NN} \underline{b_\NN}\P_{<\NN} \W_{ \alpha} + \P_{<\NN}\frac{(1+\P_{<\NN}\W)\P_{<\NN} R_\alpha}{1+\overline{\P_{<\NN} \W}}  =\P_{<\NN} G(\W_{<\NN},R_{<\NN})
\\
&R_t  +\P_{<\NN} \underline{b_\NN}\P_{<\NN} R_{\alpha} +ic\P_{<\NN}R- i \P_{<\NN}\frac{g\P_{<\NN}\W-a_\NN}{1+\P_{<\NN}\W} =
\P_{<\NN} K(\W_{<\NN},R_{<\NN}),
\end{aligned}
\right.
\end{equation}
where $G$ and $K$ represent the right hand side terms in \eqref{ww2d-diff}, and 
\[
\ub_{\NN} = \ub(W_{<\NN},R_{<\NN}), \qquad  a_{\NN} = a(R_{<\NN}). 
\]
Here $n$ is a dyadic frequency parameter, and the multiplier $\P_{<n}$  selects the frequencies less than $2^{n}$.
The expression for $\W_t$ in the first equation above retains the structure from the original equations, i.e., it is 
an exact derivative. Thus we can recover the undifferentianted variable $W$ in a way consistent 
with the $\W$ equation above, simply by integrating in time the relation
\[
W_t   + \P_{<\NN} [ \underline{F_\NN} (1+ \P_{<\NN}\W)] + i \frac{c}2 \P_{<\NN} W = 0.
\]
The above set of equations is consistent, and for each dyadic frequency scale $\NN$ it generates 
an ode in the space $\H_k$, $k \geq 2$  for $(\W,R)$, with the additional property that $W$ in $L^2$. 
A-priori, the lifespan for this system depends on $\NN$. However, the same type of estimates 
used for the linearized system yield  uniform bounds and a uniform lifespan depending only on the 
Sobolev norm of the data, and not on $n$.  Then the regular solutions are obtained as weak limits of these 
approximate solutions as $n \to \infty$.

\bigskip

{\bf Step 2 : Uniqueness of regular  solutions ($(\W,R) \in \dH_2$)}. This is achieved in a more standard fashion.
Given two solutions $(W_1,Q_1)$ and $(W_2,Q_2)$ we subtract the equations and do energy estimates 
for the difference $(\W_1-\W_2,R_1-R_2)$ in $\dH_0$, as well as simpler integrated bounds
for $W_1-W_2$ in $L^2$. Then close the argument and prove uniqueness via Gronwall's inequality.

\bigskip

 {\bf Step 3 : $\dH_1$ bounds for $(\W,R)$}. Here we use the uniform bounds for the $\dH_2$ norm
of $(\W,R)$ which depend only on the $\dH_1$ norm  of $(\W,R)$ (via the control parameters $\uA$ and $\uB$),
as in Theorem~\ref{t:en=large}, to conclude that the regular solutions can be continued up to a time 
which depends only on the $\dH_1$ norm of the data.
 
\bigskip {\bf Step 4 : Construction of rough solutions, $(\W,R) \in
  \dH_1$ }.  We regularize the data, $(W_{<k}(0), Q_{<k}(0))$, by
truncating at frequencies $< 2^k$. The corresponding solutions will be
regular, with a uniform lifespan bound. Thinking of $k$ as a
continuous parameter, we associate the $k$ derivative of the solutions
$(W_{<k}, Q_{<k})$ to solutions for the linearized equation around
$(W_{<k}, Q_{<k})$. Using both the energy bounds for the solutions
$(\W_{<k}, R_{<k})$ in $\dH_0$ and $\dH_1$, and the $\dH_0$ bounds for
the linearized equation, we show that these solutions inherit not only
uniform bounds in $\dH_1$, but also uniform frequency envelope bounds
in terms of the initial data frequency envelope. This suffices in order 
to prove strong convergence of $(\W_{<k}, R_{<k})$ in $\dH_1$,
and of $(W_{<k},Q_{<k})$ in $\dH_0$.

\bigskip {\bf Step 4 : Weak Lipschitz dependence for rough solutions }.
Here we show that for rough solutions, i.e., $(W,Q) \in \dH_0$ and $(\W,R) \in \dH_1$,
we have Lipschitz dependence on the initial data in the $\dH_0$ topology 
for both $(W,Q) \in \dH_0$ and $(\W,R)$. This is a direct consequence of the $\dH_0$ 
bounds for the linearized equation.

\bigskip {\bf Step 5 : Continuous dependence on the data for rough solutions }.
This is a standard consequence of the frequency envelope bounds in \textbf{Step 3} and the weak Lipschitz bounds
in \textbf{Step 4}.

\bigskip {\bf Step 6 : Continuation of solutions }. Here we show that
the solutions extend with uniform bounds for as long as our control
norms $\uA$, $\uB$ remain bounded. This is a consequence of the above
local well-posedness in $\dH_1$ and the energy estimates in Theorem~\ref{t:en=large}.

\end{proof}

\begin{proof}[Outline of proof of Theorem~\ref{t:cubic}]
This is an easy consequence of the $\dH_1$ well-posedness and the $\dH_1$ uniform bounds in Theorem~\ref{t:en=small},
as the control  norms $\uA$ and $\uB$ can be estimated in terms if the $\dH_1$ norm of $(\W,R)$ and the $\dH_0$ 
norm of $(W,Q)$.

\end{proof}

\appendix

\section{Water Waves Related Estimates}
\label{s:est}
In this section we gather a number of bilinear and multilinear
estimates which are used throughout the paper. Some are from \cite{HIT},
and are just recalled here.  The rest are connected to the new 
structure induced by the vorticity. We begin with some commutator bounds  
from \cite{HIT}. 
\begin{lemma}\label{l:com}
 The following commutator estimates hold:
\begin{equation}  \label{first-com}
\Vert |D|^s \left[ \P,R\right] |D|^\sigma w \Vert _{L^2}\lesssim \Vert
  |D|^{\sigma+s} R\Vert_{BMO} \Vert w\Vert_{L^2}, \qquad \sigma \geq 0, \ \ s \geq 0,
\end{equation}
\begin{equation}\label{second-com}
\Vert |D|^s  \left[ \P,R\right] |D|^\sigma w \Vert _{L^2}\lesssim \Vert
  |D|^{\sigma+s}  R\Vert_{L^2} \Vert w\Vert_{BMO}, \qquad  \sigma > 0, \ \ \, s \geq 0.
\end{equation}
\end{lemma}
We remark that later this is applied to functions which are holomorphic/antiholomorphic,
but that no such assumption is made above.  Next, we have several bilinear estimates:

\begin{lemma}\label{l:N}
The functions $N$ and $M_1$  satisfy the pointwise bounds
\begin{equation}\label{N-infty}
\| N\|_{L^\infty} \lesssim_A AA_{-1/2},\quad \Vert M_1\Vert_{L^{\infty}}\lesssim_A A^2,
\end{equation}
as well as the Sobolev bounds 
\[
\Vert N\Vert_{\dot{H}^{n+\frac{1}{2}}}\lesssim_A A\mathbf{N}_{n} \quad \Vert M_1\Vert_{\dot{H}^{n}}\lesssim_A A\mathbf{N}_{n}.
\]
\end{lemma}
\begin{proof}
We begin with the bounds for $N$, where $N$ was defined in \eqref{defN} as follows
\[
N=\P\left[W\bar{R}_{\alpha}-\bar \W R\left]+\bar{\P}\right[\bar{W}R_{\alpha}- \W \bar R\right].
\] 
For the pointwise bound we claim that 
\begin{equation}\label{N-infty1}
\| N\|_{L^\infty} \lesssim \|W\|_{B^{\frac34,\infty}_2} \|R\|_{B^{\frac14,\infty}_2}.
\end{equation}
This suffices since each of the the right hand side factors is bounded
by $\sqrt{AA_{-1/2}}$ by interpolation.  To achieve this we observe that we can
alternatively write $N$ in the form
\[
N = \ \bar{\P}[\bar W R_\alpha - \W\bar R]+\P[W \bar R_\alpha  - \bar \W R ] = 
\partial_\alpha  (\bar{\P} [ \bar{W} R ] + \P [W \bar{R}]) -(\bar{\W}R+ \W\bar{R}).
\]
We apply a bilinear Littlewood-Paley decomposition and use the first expression 
above for the high-low interactions, and the second for the high-high interactions,
to write $N = N_1+N_2$, where
\[
\begin{split}
N_1 = & \ \sum_k [\bar W_{k}  R_{<k,\alpha} -W_{<k,\alpha} \bar{R}_k ]+  [W_{k} \bar R_{<k,\alpha} - \bar{W}_{<k,\alpha} R_k], 
\\
N_2 = & \ \sum_k\partial_\alpha(\bar{\P} [ \bar W_k R_k] + \P [W_k \bar R_k])-(\bar{W}_{k, \alpha}R_{k}+W_{k, \alpha}\bar{R}_{k}).
\end{split}
\]
We estimate the terms in $N_1$ separately; we show the argument for the first term:
\[
\|   \sum_k \bar W_{k}  R_{<k,\alpha}\|_{L^\infty} \lesssim 
\sum_{j \leq k} 2^{\frac34(j-k)} \|W_j\|_{L^{\infty}} \|R_k\|_{L^{\infty}}
\lesssim \|W\|_{B^{\frac34,\infty}_2} \|R\|_{B^{\frac14,\infty}_2}.
\]
For the first term in $N_2$ we note that the multiplier $\partial_\alpha \P_{<k+4}  \P$ 
has an $O(2^k)$ $L^\infty$ bound. Hence, we can estimate
\[
\|N_2\|_{L^\infty} \lesssim \sum_k 2^k \|W_k\|_{L^\infty} \|R_k\|_{L^\infty} 
\lesssim \|W\|_{B^{\frac34,\infty}_2} \|R\|_{B^{\frac14,\infty}_2}.
\]

For simplicity, we only prove the $\dot{H}^{\frac{1}{2}}$ bound for
$N$. the rest can be done in a very similar way. To obtain the bound
we apply Coifman-Meyer type commutator estimate:
\[
\Vert N \Vert_{\dot{H}^{\frac{1}{2}}}\lesssim \Vert \vert D\vert^{\frac{1}{2}} \left[ \bar{\P}, \bar{W} \right]\!R_{\alpha}\Vert_{L^2}+ \Vert  \vert D\vert^{\frac{1}{2}}  \left[ \bar{\P}, \bar{R} \right]\! \W\Vert_{L^2} +\Vert  \vert D\vert^{\frac{1}{2}} \left[ \P, W \right] \! \bar R_\alpha \Vert_{L^2}+\Vert  \vert D\vert^{\frac{1}{2}}  \left[ \P, \W\right]R\Vert_{L^2}.
 \]
 Each commutator can be bounded by $\lesssim A\norm_0$.  For pointwise
 and Sobolev bounds of $M_1$ we refer to Lemma~2.8 (\emph{Appendix} in
 \cite{HIT}); the exact same approach applies in the current case.
\end{proof}

Essential in the article are the bounds for $a$, which we established in \cite{HIT}:
\begin{proposition}\label{l:a}
  Assume that $R \in \dot H^{\frac12} \cap \dot H^{\frac32}$.  Then the real frequency-shift $a$
 is nonnegative and satisfies the $BMO$ bound
\begin{equation}\label{a-bmo}
\| a\|_{BMO}    \lesssim \Vert R\Vert_{BMO^\frac12}^2,
\end{equation}
and the uniform bound
\begin{equation}\label{a-point}
\| a\|_{L^\infty}    \lesssim \Vert R\Vert_{\dot B^{\frac12}_{\infty,2}}^2.
\end{equation}
 Moreover,
\begin{equation}\label{a-bmo+}
\Vert |D|^{\frac{1}{2}} a\Vert_{BMO}\lesssim \Vert R_\alpha \Vert_{BMO}
\Vert |D|^\frac12 R\Vert_{L^\infty}, \qquad
\Vert  a\Vert_{B^{\frac12,\infty}_2}\lesssim \Vert R_\alpha \Vert_{B^{\frac12,\infty}_2}
\Vert |D|^\frac12 R\Vert_{L^\infty},
\end{equation}
\begin{equation}\label{aflow}
\Vert (\partial_t+b\partial_{\alpha})a \Vert_{L^{\infty}}\lesssim AB,
\end{equation}
and
\begin{equation}\label{a-Hs}
\Vert  a\Vert_{\dot H^s}\lesssim \Vert  R \Vert_{\dot H^{s+\frac12}}
\Vert R\Vert_{BMO^\frac12}, \qquad s > 0.
\end{equation}
\end{proposition}

Here we need to supplement this with bounds for $a_1$. One notable difference 
between the two is that $a_1$ has a linear component, whereas $a$ is purely quadratic.
For various estimates we need to separate the two components of $a_1$:

\begin{proposition}\label{l:a1}
  Assume that $R \in \dot H^{\frac12} \cap \dot H^{\frac32}\cap L^{\infty}$. Then 
  \[
  \Vert a_1\Vert_{L^{\infty}} \lesssim _A A_{-1/2} (1+A).
  \]
 Moreover, the following estimate holds
\[
\Vert (\partial_t +\underline{b} \partial_{\alpha})a_1  +2 g\Im \W  - 2c \Im R \Vert_{L^{\infty}}\lesssim AB +cA^2.
\]
\end{proposition}
\begin{proof}
We first recall the expression for $a_1$:
\[
a_1=R+\bar{R}-N.
\]
Using the equation for $R$, we only need to prove the pointwise bound for $(\partial_t +\underline{b} \partial_{\alpha})N$. We begin with the following computation
\[
(\partial_t +\underline{b} \partial_{\alpha})a_1  +2 g\Im \W  - 2c \Im R=-2a\Im Y -c\Im \frac{R\W+\bar{R}\W +N}{1+\W} -(\partial_t +\underline{b} \partial_{\alpha})N,
\]
Based on the previously established pointwise bounds for $a$ and $N$ we can estimate all but the 
last term. This is considerably more delicate. However, as seen in the proof of Lemma~\ref{l:N},
the expression for $N$ exhibits exactly the same cancellation structure for $high \times high 
\to low$ interactions as we have in the bilinear expression for $a$. Hence, the argument 
in the proof of \eqref{aflow} in \cite{HIT} immediately adapts here.

 \end{proof}

\section{Holomorphic coordinates}
\label{s:eq}

Our goal here is to introduce the holomorphic coordinates and the
holomorphic set of variables $(W,Q)$ describing the free surface,
respectively the holomorphic velocity potential restricted to the free
surface, and to derive the evolution equations for $(W,Q)$.

We proceed as in \cite{HIT,IT3}. Let $\mathbb{H}$ be the lower
halfspace, with complex coordinates denoted by $\alpha+i\beta$.  Let
$\mathcal{F}:\mathbb{H}\rightarrow \Omega (t)$ to be the conformal
transformation that maps the $\alpha$-axis into $\Gamma (t)$, which is
unique up to $\alpha$ translations.  The $x$ and $y$ coordinates are given by
\[
x=x(\alpha, \beta, t), \quad y=y(\alpha, \beta, t),
\]
 and satisfy the Cauchy-Riemann equations
\begin{equation*}
x_{\alpha}=y_{\beta}, \quad x_{\beta}=-y_{\alpha}.
\end{equation*}
We fix the conformal map by assuming that 
\[
x(\alpha, \beta, t)+iy(\alpha, \beta, t) - \alpha + i\beta \to 0
\qquad \alpha, \beta \rightarrow \infty.
\]
To switch between the two sets of coordinates we will use
\begin{equation*}
\begin{aligned}
&\partial_{\alpha}=x_{\alpha}\partial_{x}+y_{\alpha}\partial_{y}, \qquad \quad  \partial_{\beta}=x_{\beta}\partial_{x}+y_{\beta}\partial_{y}, \\
&\partial_x=\frac{1}{j}(x_{\alpha}\partial_{\alpha}+ x_{\beta}\partial_{\beta}), \quad \ \partial_y=\frac{1}{j}(y_{\alpha}\partial_{\alpha}+ y_{\beta}\partial_{\beta}),
\end{aligned}
\end{equation*}
where
\[
j=x_{\alpha}y_{\beta}-x_{\beta}y_{\alpha}=x^2_{\alpha}+x_{\beta}^2=y^2_{\alpha}+y_{\beta}^2.
\]

Thus,
\[
\mathbf{u}=(u,v)=(cy+\varphi_x, \varphi_y).
\]

To switch between the two sets of coordinates we will use
\begin{equation*}
\begin{aligned}
&\partial_{\alpha}=x_{\alpha}\partial_{x}+y_{\alpha}\partial_{y}, \qquad \quad  \partial_{\beta}=x_{\beta}\partial_{x}+y_{\beta}\partial_{y}, \\
&\partial_x=\frac{1}{j}(x_{\alpha}\partial_{\alpha}+ x_{\beta}\partial_{\beta}), \quad \ \partial_y=\frac{1}{j}(y_{\alpha}\partial_{\alpha}+ y_{\beta}\partial_{\beta}),
\end{aligned}
\end{equation*}
where
\[
j=x_{\alpha}y_{\beta}-x_{\beta}y_{\alpha}=x^2_{\alpha}+x_{\beta}^2=y^2_{\alpha}+y_{\beta}^2.
\]
We also need a similar relation for the time derivative.
 Assume we have
\[
\mathbb{H}\xrightarrow{\mathcal{F}}\Omega (t)\xrightarrow{f} \mathbb{C}
\]
for an arbitrary function $f$. Let 
\[
g(\alpha,\beta, t)=f(x(\alpha, \beta, t), y(\alpha, \beta, t),t),
\]
 then
\[
g_t=f_t+x_tf_x+f_yy_t.
\]
So 
\begin{equation}
\label{chain}
\begin{aligned}
f_t&=g_t-x_tf_x-f_yy_t\\
&=g_t-x_t\frac{1}{j}(x_{\alpha}g_{\alpha}+x_{\beta}g_{\beta})-y_t\frac{1}{j}(y_{\alpha}g_{\alpha}+y_{\beta}g_{\beta})\\
&=g_t-\frac{1}{j}(x_t x_{\alpha}+y_ty_{\alpha})g_{\alpha}-\frac{1}{j}(x_t x_{\beta}+y_ty_{\beta})g_{\beta}.
\end{aligned}
\end{equation}

Define $\psi $ to be the (harmonic) composition of $\mathcal{F}$ to
the velocity potential $\varphi$,
\[
\psi =\varphi \circ \mathcal{F}:\mathbb{H}\rightarrow \mathbb{R}, \qquad \psi (\alpha, \beta, t)=\varphi (x(\alpha,\beta, t), y(\alpha, \beta, t), t).
\]
The velocity components $(u,v)$ can now be expressed in terms of the
velocity potential $\psi$ by
\begin{equation}
\label{rel1}
\left\{
\begin{aligned}
&u=\frac{1}{j}(x_{\alpha}\psi_{\alpha}+x_{\beta}\psi_{\beta})+cy\\
&v=\frac{1}{j}(y_{\alpha}\psi_{\alpha}+y_{\beta}\psi_{\beta}).
\end{aligned}
\right.
\end{equation}

It follows that 
\[
u^2+v^2=\varphi_{x}^2+\varphi_y^2=\frac{1}{j}(\psi_{\alpha}^2+\psi_{\beta^2}).
\]

\subsection{Boundary values}

Setting $\beta =0$ gives the boundary values of the holomorphic
functions defined in the lower half plane. In particular, we introduce
the notation $x(\alpha, 0, t)=:X(\alpha, t)$ and $y(\alpha, 0,
t)=:Y(\alpha, t)$, so that $\alpha \to X+iY$ parametrizes
 $\Gamma (t)$. The function $(x-\alpha)+i(y-\beta)$  is holomorphic in $\mathbb{H}$
and decays at infinity, which implies that on the boundary we have
\[
\left\{
\begin{aligned}
&Y=H(X-\alpha)\\
&X=\alpha +HY.
\end{aligned}
\right.
\]
Also set  
\[
z(\alpha,\beta,t)=x(\alpha, \beta,t)+iy(\alpha, \beta,t), \qquad Z(\alpha,t) = z(\alpha,0,t).
\]
Then our ``holomorphic'' variable $W$, which describes the surface $\Gamma (t)$, will be 
\begin{equation*}
W = Z - \alpha.
\end{equation*}
 As $z$ is holomorphic in $\mathbb{H}$, so is $ \dfrac{1}{z_\alpha}-1$; further, it decays
at infinity. Its boundary values on the real axis are given by
\[
\frac{1}{Z_\alpha -1} =\left( \frac{X_{\alpha}}{J}-1\right)-i\frac{Y_{\alpha}}{J},
\]
which leads to the following relations:
\begin{equation}
\label{rel2}
\left\{
\begin{aligned}
&\frac{X_{\alpha}}{J}-1=-H\left[\frac{Y_{\alpha}}{J} \right],\\
&\frac{Y_{\alpha}}{J}=H\left[ \frac{X_{\alpha}}{J}-1\right]  .
\end{aligned}
\right.
\end{equation}

We also introduce the notation $\Psi (\alpha, t):=\psi(\alpha, \beta, t)$
for the real velocity potential restricted to $\Gamma (t)$ and expressed in 
holomorphic coordinates, and at the same time define $\Theta (\alpha, t)$ by
\[
\left\{
\begin{aligned}
\Psi &=H\Theta\\
\Theta &=-H\Psi.
\end{aligned}
\right.
\]
Up to a constant this is the trace of the stream function $\theta$ on the free boundary.
Since $\psi$ is harmonic in the lower half plane, we have 
\[
\psi_{\beta}|_{\beta =0}-H\Psi_{\alpha}=-\Theta_{\alpha}.
\]
Our holomorphic velocity potential will be the function 
\[
Q = \Psi + i \Theta.
\]
Further we need to focus on the two boundary conditions: kinematic and dynamic. 

\subsection{The kinematic boundary condition.} 
The kinematic boundary condition states that the normal component of
the velocity of the boundary is equal to the normal component of the
fluid velocity, meaning that
\[
(X_t, Y_t)\cdot (-Y_{\alpha}, X_{\alpha})=(u,v)\cdot  (-Y_{\alpha}, X_{\alpha}),
\]
where $ (-Y_{\alpha}, X_{\alpha})$ is a normal to $\Gamma (t)$. Expanding the expression above and using \eqref{rel1} re-expresses the kinematic boundary condition in holomorphic coordinates:
\begin{equation}
\label{eq1}
X_{\alpha}Y_{t}-Y_{\alpha}X_{t}=H\Psi_{\alpha}-cYY_{\alpha}=-\Theta_{\alpha}-cYY_{\alpha}.
\end{equation}
The goal now is to obtain a second equations for $X_t$ and $Y_t$, and then solve for an explicit form of those boundary values. Divide \eqref{eq1} by $J$
\begin{equation*}
\begin{aligned}
\left( \frac{X_{\alpha}}{J}-1\right) Y_{t}-\frac{Y_{\alpha}}{J}X_{t}+Y_t&=-\frac{\Theta_{\alpha}}{J}-c\frac{YY_{\alpha}}{J},
\end{aligned}
\end{equation*}
 and use \eqref{rel2} to obtain 
 \begin{equation*}
 -H\left[ \frac{Y_{\alpha}}{J}\right] Y_{t}-\frac{Y_{\alpha}}{J}H\left[Y_t \right] +Y_t=-\frac{\Theta_{\alpha}}{J}-c\frac{YY_{\alpha}}{J},
 \end{equation*}
which further simplifies to
\begin{equation*}
 \frac{Y_{\alpha}}{J}Y_{t}-H\left[ \frac{Y_{\alpha}}{J}\right] H\left[Y_t \right] +X_t=-H\left[ \frac{\Theta_{\alpha}}{J}\right] -cH\left[ \frac{YY_{\alpha}}{J}\right].
 \end{equation*}
Thus, a second equation for $\left( X_t, Y_t\right) $ is
\begin{equation}
\label{eq2}
X_{\alpha}X_t+Y_{\alpha}Y_t=-JH\left[ \frac{\Theta_{\alpha}}{J}\right] -cJH\left[ \frac{YY_{\alpha}}{J}\right].
 \end{equation}
From \eqref{eq1} and \eqref{eq2} we have
\begin{equation}
\label{eq3}
\left\{
\begin{aligned}
X_{t}=&-H\left[\frac{\Theta_{\alpha}}{J}\right] X_{\alpha} -cH\left[ \frac{YY_{\alpha}}{J}\right]X_{\alpha}+\frac{\Theta_{\alpha}}{J}Y_{\alpha}+c\frac{YY_{\alpha}}{J}Y_{\alpha} \\
Y_{t}=&-H\left[ \frac{\Theta_{\alpha}}{J}\right]Y_{\alpha}-cH\left[ \frac{YY_{\alpha}}{J}\right]Y_{\alpha}-\frac{\Theta_{\alpha}}{J}X_{\alpha}-c\frac{YY_{\alpha}}{J}X_{\alpha}.  
\end{aligned}
\right.
\end{equation}

\subsection{The dynamic boundary condition.} We have already
determined the the spatial form of the dynamic boundary condition in
\eqref{imp}. From \eqref{chain} and the kinematic boundary conditions
\eqref{eq1}-\eqref{eq2} it follows that, on on the boundary, $\varphi_t
$ is given as
\[
\varphi_t|_{\beta =0}= \Psi _t -\frac{1}{J} \left(  X_{\alpha}X_t +Y_{\alpha}Y_t\right)\Psi_{\alpha}+\frac{1}{J}\left( X_{\alpha}Y_t-Y_{\alpha}X_t\right)\Theta_{\alpha}.  
\]
Substituting $X_t$ and $Y_t$ from \eqref{eq3} yields
\begin{equation}
\varphi_{t}|_{\beta=0}=\Psi_{t}+H\left[\frac{\Theta_{\alpha}}{J} \right]\Psi_{\alpha}+c H\left[\frac{YY_{\alpha}}{J} \right]\Psi_{\alpha}-\frac{1}{J}\Theta ^2_{\alpha}-\frac{1}{J}cYY_{\alpha}\Theta_{\alpha}. 
\end{equation}
We  express all the terms on the right in  \eqref{imp}
terms of the traces on the boundary of the corresponding functions. 
Doing so, after some simplifications  we arrive at the following equation:
\begin{equation}\label{eq-Psi}
\Psi_{t}+H\left[\frac{\Theta_{\alpha}}{J} \right]\Psi_{\alpha}-\frac{\Theta ^2_{\alpha}}{J}+\frac{1}{2J}\left( \Psi_{\alpha}^2 +\Theta_{\alpha}^2 \right)+gY+c H\left[\frac{YY_{\alpha}}{J} \right]\Psi_{\alpha}-c\Theta +c\frac{Y}{J}  X_{\alpha}\Psi_{\alpha}=0.
\end{equation}

\subsection{The real form of the equations.}
The equations \eqref{eq3} and \eqref{eq-Psi} provide us with a system
describing the evolution of the free boundary and the velocity
potential restricted to the free boundary, as follows:
\[
\left\{
\begin{aligned}
&Y_{t}=-H\left[ \frac{\Theta_{\alpha}}{J}\right]Y_{\alpha}-cH\left[ \frac{YY_{\alpha}}{J}\right]Y_{\alpha}-\frac{\Theta_{\alpha}}{J}X_{\alpha}-c\frac{YY_{\alpha}}{J}X_{\alpha}\\  
&\Psi_{t}=-H\left[\frac{\Theta_{\alpha}}{J} \right]\Psi_{\alpha}+\frac{\Theta ^2_{\alpha}}{J}-\frac{1}{2J}\left( \Psi_{\alpha}^2 +\Theta_{\alpha}^2 \right)-gY-c H\left[\frac{YY_{\alpha}}{J} \right]\Psi_{\alpha}+c\Theta -c\frac{Y}{J}  X_{\alpha}\Psi_{\alpha}.
\end{aligned}
\right.
\]
Here $X$ and $\Theta $ are dependent variables.

The Hamiltonian associated to the system is
\[
\mathcal{E}(Y,\Psi)=\frac{1}{2}\int\left\lbrace \Psi \vert \partial_{\alpha} \vert \Psi +gY^2X_{\alpha}+c\Psi_{\alpha}Y^2+\frac{1}{3}c^2Y^3X_{\alpha}\right\rbrace  \, d\alpha,
\]
where $\vert \partial_{\alpha}\vert =H\partial_{\alpha}$. Thus
\[
\left\{
\begin{aligned}
&\frac{\delta \mathcal{E}}{\delta Y}=gYX_{\alpha}+g\frac{1}{2}\vert \partial_{\alpha}\vert (Y^2)+c\Psi_{\alpha}Y+\frac{1}{2}c^2Y^2X_{\alpha}+\frac{1}{6}c^2\vert \partial_{\alpha}\vert (Y^3)\\
&\frac{\delta \mathcal{E}}{\delta \Psi}=\vert \partial_{\alpha}\vert \Psi -cYY_{\alpha}.
\end{aligned}
\right.
\]

We can write the above  equations for $(Y,\Psi)$ in a skew-symmetric form
\begin{equation}
\label{skew}
\left( \begin{array}{c}
Y   \\
\Psi   \\
 \end{array} \right)_t =\left( \begin{array}{cc}
0 & \mathbf{A}  \\
-\mathbf{A}^{*} & \mathbf{B} \\
 \end{array} \right)  \left( \begin{array}{c}
\delta \mathcal{E}/\delta Y   \\
\delta \mathcal{E}/\delta \Psi \\
 \end{array} \right)
 \end{equation}
 where 
 \begin{equation*}
 \left\{
 \begin{aligned}
 &\mathbf{A}=\frac{X_{\alpha}}{J}+Y_{\alpha}H\frac{1}{J}, \quad \mathbf{A}^{*}=\frac{X_{\alpha}}{J}-\frac{1}{J}HY_{\alpha},\\
 &\mathbf{B}=\Psi_{\alpha}H\frac{1}{J}+\frac{1}{J}H\Psi_{\alpha}-c\partial_{\alpha}^{-1}, \qquad \mathbf{B}^{*}=-\mathbf{B}.
 \end{aligned}
 \right.
 \end{equation*}
 This form immediately implies that $\mathcal{E}$ is conserved along
 the flow.
There are several Hamiltonian symmetries, which correspond to conservation
laws of the system in accordance with Noether's principle. The horizontal translation
invariance
\[
Y(\alpha, t)\rightarrow Y(\alpha +a, t),\quad  \Psi(\alpha, t)\rightarrow \Psi(\alpha +a, t),
\]
is generated by the functional $\mathcal{P}$ which we will derive below. This functional  corresponds to total horizontal momentum, and it is given by the following system
\begin{equation}
\label{momentum}
\left( \begin{array}{c}
Y   \\
\Psi   \\
 \end{array} \right)_\alpha =\left( \begin{array}{cc}
0 & \mathbf{A}  \\
-\mathbf{A}^{*} & \mathbf{B} \\
 \end{array} \right)  \left( \begin{array}{c}
\delta \mathcal{P} /\delta Y  \\
\delta \mathcal{P} /\delta \Psi , \\
 \end{array} \right)
 \end{equation}
where 
\[
\mathcal{P}(Y, \Psi)=\int \left\lbrace\Psi Y_{\alpha}-\frac{c}{2}Y^2X_{\alpha}\right\rbrace  \, d\alpha.
\]

The variational derivatives of $\mathcal{P}$ are
\[
\left\{
\begin{aligned}
&\delta \mathcal{P} /\delta Y=-\Psi_{\alpha} -cYX_{\alpha}-\frac{c}{2}\vert \partial_{\alpha}\vert (Y^2)\\
&\delta \mathcal{P} /\delta \Psi =Y_{\alpha}.
\end{aligned}
\right.
\]

Thus, our conserved energies are
\begin{equation}
\label{energies}
\left\{
\begin{aligned}
&\mathcal{E}(Y,\Psi)=\frac{1}{2}\int\left\lbrace \Psi \vert \partial_{\alpha} \vert \Psi +gY^2X_{\alpha}+c\Psi_{\alpha}Y^2+\frac{1}{3}c^2Y^3X_{\alpha}\right\rbrace  \, d\alpha\\
&\mathcal{P}(Y, \Psi)=\int \left\lbrace\Psi Y_{\alpha}-\frac{c}{2}Y^2X_{\alpha}\right\rbrace  \, d\alpha ,
\end{aligned}
\right.
\end{equation}
where $\vert \partial_{\alpha}\vert =H\partial_{\alpha}$.

\subsection{The complex form of the equations} 
Recall that our holomorphic variables are
\[
Z=X+iY, \quad Q=\Psi+i\Theta.
\]
We also introduce the notation
\begin{equation}
\label{def-F}
F=H\left[ \frac{\Theta_{\alpha}}{J}\right]+i\frac{\Theta_{\alpha}}{J} =\mathbf{P}\left[ \frac{Q_{\alpha}-\bar{Q}_{\alpha}}{J}\right], \qquad J=\vert Z_{\alpha}\vert ^2,
\end{equation}
noting that $F$ is also the boundary value of a function which is holomorphic in
the lower half plane. 

Using these notations, a straightforward computation allows us to
re-expresses the kinematic boundary conditions \eqref{eq3} in one
single equation for the motion of the boundary
 \begin{equation*}
Z_t+FZ_{\alpha}+2ic\mathbf{P}\left[ \frac{YY_{\alpha}}{J}\right]Z_{\alpha}=0,
\end{equation*}
which is equivalent to 
\begin{equation}
\label{kbc}
Z_t+FZ_{\alpha}-i\frac{c}{4}\mathbf{P}\left[ \frac{\partial_{\alpha}(Z-\bar{Z})^2}{J}\right]Z_{\alpha}=0. 
\end{equation}

Next we rewrite the dynamic boundary condition. For this we apply the
operator $2\P = \mathbf{I}-iH$ in the equation
\eqref{eq-Psi}.  We do not repeat the computations in \cite{HIT} for the case $c=0$.
Instead, we use them directly to obtain, from  \eqref{eq-Psi}, the following 
equation:
\begin{equation}
\label{dbc}
Q_t-ig(Z-\alpha)+FQ_{\alpha}+icQ+\mathbf{P}\left[\frac{\vert Q _{\alpha}\vert ^2}{J}\right]+cK =0,
\end{equation}
where 
\[
K = (\mathbf{I}-iH)\left[ H\left[\frac{YY_{\alpha}}{J} \right]\Psi_{\alpha}+\frac{Y}{J}  X_{\alpha}\Psi_{\alpha}   \right].
\]
To simplify $K$ we use the relation $H(fg - H fH g) = f H g + H f g$ to rewrite it as 
\[
\begin{split}
K = & \ 
H\left[\frac{YY_{\alpha}}{J} \right]\Psi_{\alpha}
- i  H\left[\frac{YY_{\alpha}}{J} \Theta_\alpha \right] + i   \frac{YY_{\alpha}}{J} \Psi_\alpha
+  i  H\left[\frac{YY_{\alpha}}{J} \right] \Theta_\alpha + (\mathbf{I} - i H)\left[\frac{Y}{J}  X_{\alpha}\Psi_{\alpha}   \right]
\\ = & \ i(\mathbf{I}-iH) \left[\frac{YY_{\alpha}}{J} \right] Q_\alpha + 
 (\mathbf{I} - i H)\left[\frac{Y}{J}  (X_{\alpha}\Psi_{\alpha}+Y_\alpha \Theta_\alpha )  \right]
\\ = & \ \P \left[\frac{Y}{\bar Z_{\alpha}}  - \frac{Y}{Z_{\alpha}}  \right] Q_\alpha+  \P \left[ \frac{Y Q_\alpha}{Z_\alpha} + \frac{Y \bar Q_\alpha}{\bar Z_\alpha}\right] 
\end{split}
\]
Further, we write $Y = - i((Z-\alpha) - (\bar Z-\alpha))$, eliminate the projected antiholomorphic terms
and remove the projection in front of holomorphic terms to obtain
\begin{equation}\label{K}
K = - \frac{i}{2} \P \left[\frac{Z-\alpha}{\bar Z_{\alpha}} +
  \frac{\bar Z-\alpha}{Z_{\alpha}} \right] Q_\alpha- \frac{i}2 \P \left[  \frac{(Z-\alpha) \bar
    Q_\alpha}{\bar Z_\alpha}-
  \frac{(\bar Z-\alpha) Q_\alpha}{Z_\alpha} \right].
\end{equation}
Thus our set of holomorphic equations consists of \eqref{kbc} and \eqref{dbc}-\eqref{K}.

For the last step we use the relation $Z=W+\alpha$ to substitute $Z$ by $W$. 
For the last expression in  \eqref{kbc} we have 
\[
\begin{aligned}
\mathbf{P}\left[ \frac{\partial_{\alpha}(Z-\bar{Z})^2}{J}\right]
= & \ 2\mathbf{P}\left[  \frac{Z-\alpha}{\bar Z_\alpha}  -  \frac{Z-\alpha}{Z_\alpha} + 
\frac{\bar Z-\alpha}{Z_\alpha} \right] \\
=&  2\mathbf{P}\left[  \frac{W}{1+\bar W_\alpha}  -  
\frac{W}{1+W_\alpha} + \frac{\bar W}{1+W_\alpha} \right] 
\\   = &2F_1 - \frac{2W}{1+W_\alpha} ,
\end{aligned}
\]
where 
\[
F_1 = \mathbf{P}\left[  \frac{W}{1+\bar W_\alpha}  + \frac{\bar W}{1+W_\alpha} \right] .
\]
For $K$ on the other hand we have 
\[
K = -\frac{i}2 F_1 Q_\alpha - \frac{i}2 \P\left[ \frac{W  \bar Q_\alpha}{1+\bar W_\alpha} - 
\frac{ \bar W Q_\alpha}{ 1+W_\alpha}\right] .
\]
Hence, setting 
\[
\underline{F} = F - i \frac{c}2 F_1,
\]
our equations become  
\begin{equation}\label{whole1}
\left\{
\begin{aligned}
&W_t+\underline{F} (W_{\alpha}+1)+ i\frac{c}{2}W=0\\
&Q_t+\underline{F}Q_{\alpha}+icQ-igW+\mathbf{P}\left[\frac{\vert Q _{\alpha}\vert ^2}{J}\right]+
i\frac{c}{2}\left\lbrace  \P\left[  \frac{W\bar{Q}_{\alpha}}{1+\bar{W}_{\alpha}} \right]   -\P\left[ \frac{\bar{W}Q_{\alpha}}{1+W_{\alpha}}\right] \right\rbrace=0  .
\\
\end{aligned}
\right.
\end{equation}
We can also re-express the Hamiltonian and horizontal momentum in terms of the holomorphic variables 
$(W,Q)$. This gives
\begin{equation}\label{ww-energy}
\begin{aligned}
\E(W,Q) =& \Re \int  g\vert W\vert^2(1+W_\alpha) - iQ\bar{Q}_{\alpha}
+ c  Q_{\alpha} (\Im W)^2 -
\frac{c^3}{2i} |W|^2W(1+W_\alpha)  \, d\alpha .\\
\end{aligned}
\end{equation}
A second conserved quantity is the horizontal momentum,
\begin{equation}
\mathcal{P}=\int \left\lbrace\frac{1}{i}\left( \bar{Q}W_{\alpha}-Q\bar{W}_{\alpha}\right)-c\vert W\vert^2+\frac{c}{2}\left( W^2\bar{W}_{\alpha}+\bar{W}^2W_{\alpha}\right)\right\rbrace\, d\alpha .
\label{horizontal-m}
\end{equation}

\end{document}